\documentclass[12pt]{article}
\usepackage{authblk}
\usepackage[letterpaper,top=2cm,bottom=2cm,left=2.5cm,right=2.5cm,marginparwidth=1.75cm]{geometry}
\usepackage{setspace}
\usepackage{amsmath,amsfonts,amssymb,amsthm,stmaryrd}
\usepackage[T1]{fontenc}
\usepackage{textcomp}
\usepackage{graphicx}
\usepackage{epstopdf}
\usepackage{caption}
\usepackage{subfig}
\usepackage{float}
\usepackage{booktabs}
\usepackage{adjustbox}
\usepackage{multirow}
\usepackage{xcolor}
\usepackage[colorlinks=true, allcolors=blue]{hyperref}
\usepackage{tikz}
\usetikzlibrary{arrows.meta,positioning,calc,fit,backgrounds,matrix}
\usepackage{algorithm}
\usepackage{algpseudocode}
\usepackage[numbers,sort&compress]{natbib}

\newtheorem{theorem}{Theorem}
\newtheorem{lemma}{Lemma}
\newtheorem{proposition}{Proposition}
\newtheorem{corollary}{Corollary}
\theoremstyle{definition}

\theoremstyle{remark}
\newtheorem{remark}{Remark}
\newcommand{\R}{\mathbb{R}}
\newcommand{\Q}{Q}
\newcommand{\bxi}{\boldsymbol{\xi}}
\newcommand{\norm}[1]{\lVert #1\rVert}

\title{CS-PIKAN: Curvilinear Separable Physics-Informed Kolmogorov-Arnold Networks with Autograd-Based Metric Evaluation and Certified Error Bounds\\ for PDEs on Complex Domains}
\author[1]{Mohammad E. Heravifard}
\author[1]{Kazem Hejranfar\thanks{Corresponding author: khejran@sharif.edu (K. Hejranfar)}}
\affil[1]{\it{Department of Aerospace Engineering, Sharif University of Technology, Tehran, Iran}}
\date{}

\begin{document}
\maketitle
\begin{abstract}
Physics-Informed Kolmogorov-Arnold Networks (PIKANs) offer effective neural solvers for partial differential equations (PDEs), but their conventional formulations rely on Cartesian domains and scattered collocation points, often resulting in high computational costs, approximate boundary enforcement, and no rigorous error guarantees. We introduce CT-PIKAN, which maps curvilinear domains to the unit square and computes all geometric quantities, including those required by the Laplace-Beltrami operator, through automatic differentiation. Building on this transformation, Curvilinear Separable PIKAN (CS-PIKAN) represents solutions as rank-$R$ sums of products of one-dimensional KANs with topology-adapted bases (Fourier for periodic and Chebyshev for bounded coordinates and time). A separable envelope and transfinite lift enforce boundary and initial conditions exactly, while a single $\sqrt{g}$-weighted residual is minimized on a stratified tensor grid, reducing network cost from the product to the sum of grid sizes. We derive a metric-explicit Poincar\'e constant, establish computable a-posteriori error bounds for elliptic, parabolic, and hyperbolic problems, and prove that the loss provides an unbiased estimator of the squared residual norm underlying these bounds. Nine benchmarks cover advection, Poisson, and heat conduction on Cartesian, polar, wavy, annular, star-shaped, and sampled-boundary domains. Across three seeds, CS-PIKAN reduces the relative $L_2$ error of CT-PIKAN by factors of 7-126 at equal or lower training cost (up to 10$\times$ faster), with all certified bounds satisfied. Ablations highlight the importance of exact periodicity and boundary constraints; KAN factors outperform MLPs for transport-dominated problems but not smooth elliptic ones. Nevertheless, classical body-fitted solvers remain substantially cheaper at comparable accuracy for these two-dimensional linear problems.
\end{abstract}


\newpage

\section{Background}
\label{sec:intro}

Partial differential equations (PDEs) underpin the modeling of fluid dynamics, heat transfer, wave propagation, electromagnetics, and material deformation \cite{f1,p1}. Finite difference, finite volume, and finite element methods solve them reliably on structured and unstructured grids \cite{g1,g2,g3}, but complex or curvilinear geometries require mesh generation and problem-dependent discretizations \cite{g3,g4}. Body-fitted grids handle such geometries through a smooth map from a reference square, after which the equations are rewritten with metric coefficients \cite{g4,Thompson1985,Gordon1973}; deriving and discretizing these terms, however, is problem dependent and error prone. Physics-Informed Neural Networks (PINNs) instead embed the governing equations and boundary and initial conditions in the loss of a neural network \cite{PINN1,PINN11}, using automatic differentiation to avoid explicit meshing and discretization \cite{PINN3,AD}. They may nevertheless struggle with stiff and multi-scale problems \cite{PINN4,PINN5}, and irregular geometries cause sampling distortion, gradient imbalance, and inaccurate boundary conditions \cite{PINN6,PINN7}. More generally, physics-informed solvers impose boundary and initial data through penalty terms that compete with the residual \cite{PINN3,PINN6}, pay a per-step cost that grows with the number of collocation points and the depth of nested differentiation, and provide no guarantee that a small training loss implies a small solution error \cite{Mishra2023,Hillebrecht2022}.\\ 

Kolmogorov-Arnold Networks (KANs) replace the fixed activations of multilayer perceptrons (MLPs) with learnable univariate edge functions, improving interpretability and the representation of nonlinearity \cite{KAN1}. Their physics-informed form (PIKAN) \cite{KAN2} has been applied to a wide range of PDEs with gains in accuracy, efficiency, and scalability \cite{KAN3,KAN33}, including solid mechanics, multiscale systems, and dynamical processes \cite{KAN2,PIKAN3}, and variants such as ChebPIKAN, SPIKAN, and Scaled-cPIKAN target fluid dynamics and high-dimensional problems \cite{PIKAN13,PIKAN14,Mos1,Mos2}. Most PINN and PIKAN solvers, however, assume rectangular or easily parameterized domains, and curvilinear geometries remain challenging \cite{XPINN,PhyGeoNet}: direct training on them suffers from poor convergence, distorted sampling, and boundary conditions that are hard to represent. Geometry-aware approaches address this in several ways. Coordinate-transformed PINNs require analytical Jacobians and metrics for each geometry \cite{Hwang2022}. PhyGeoNet combines elliptic grid generation with convolutional networks and finite-difference metrics \cite{PhyGeoNet}. JacobiNet learns the map and differentiates it automatically \cite{JacobiNet}, and networks can also generate structured grids \cite{Khairullin2025}. Other methods encode geometry through Laplace--Beltrami eigenfunctions \cite{DeltaPINN}, learned deformations \cite{GeoFNO}, or learned input metrics for KANs \cite{GeoKAN}. All of them still train a point-wise network on scattered points, using the transformation, if at all, only to rewrite the operator. Separable PINNs and PIKANs gain order-of-magnitude speed-ups by representing the solution as a sum of products of one-dimensional networks on tensor grids \cite{SPINN,KAN3}, but they are restricted to box-shaped domains; for irregular domains, the SPINN authors mask exterior points and mention coordinate transformations only as a possible alternative \cite{SPINN}.\\

This work develops the coordinate transformation in two steps. We first formulate a Coordinate-Transformed PIKAN (CT-PIKAN) for general curvilinear domains. Its core idea is a smooth and invertible map from a user-defined physical domain onto the unit square, under which the governing equations are reformulated and enforced within the physics-informed loss. All geometric quantities (the Jacobian matrix, inverse Jacobian, metric tensor and transformed differential operators) are computed by automatic differentiation \cite{AD,AD2,AD4} directly from the map, which removes the need for manual derivation of geometry-specific terms. CT-PIKAN serves as the point-wise baseline throughout the paper; we use it to identify precisely what a coordinate transformation does and does not achieve when the network is still trained on scattered points. The second and main step is based on the observation that a coordinate transformation turns \emph{every} curvilinear domain into a box. A box is exactly the setting in which separable representations, tensor-product quadrature, exact imposition of boundary data and coordinate-adapted bases become available. The resulting Curvilinear Separable PIKAN (CS-PIKAN) exploits all four. These ingredients also give exactly the structure needed for computable error bounds. The main contributions of this work are:
\begin{enumerate}
\item \textbf{Autograd-based curvilinear formulation (CT-PIKAN).} A geometry-aware physics-informed KAN formulation in which the metric tensor, its inverse and determinant, and the derivatives that enter the Laplace--Beltrami operator are obtained by automatic differentiation of the coordinate map. Closed-form, numerically constructed and sampled-boundary maps are all handled by the same code.
\item \textbf{A new architecture (CS-PIKAN).} The pulled-back solution is a rank-$R$ separable sum of one-dimensional Chebyshev-KANs with a coordinate-adapted first layer (Fourier for closed coordinates, Chebyshev for bounded coordinates and time). It is combined with an exact-constraint ansatz consisting of a separable envelope and a transfinite (Boolean-sum) lift, so the training loss is a single, weight-free term. The curvilinear Laplace--Beltrami and transport operators are assembled exactly from one-dimensional derivatives and metric fields (Proposition~\ref{prop:sep}), so the network cost scales with $\sum_k N_k$ instead of $\prod_k N_k$.
\item \textbf{Certified error bounds.} We prove a Poincar\'e inequality on the physical domain whose constant is explicit in the metric of the map (Theorem~\ref{thm:poincare}). From it we derive a-posteriori bounds for elliptic, parabolic and hyperbolic problems (Theorem~\ref{thm:apost}), and we show that the $\sqrt g$-weighted CS-PIKAN loss is an unbiased estimator of the squared residual norm that enters these bounds (Proposition~\ref{prop:unbiased}). Exact constraints and exact periodicity are precisely the hypotheses the bounds need, which gives a mathematical reason for the architectural choices.
\item \textbf{Systematic evaluation.} Nine benchmarks covering hyperbolic, elliptic and parabolic problems on seven maps are studied, including one map constructed numerically from sampled boundary points. Each benchmark is compared with CT-PIKAN under the configurations of its public code, over three seeds, against analytical or Richardson-extrapolated finite-difference references. An ablation isolates each component, a comparison of KAN and MLP factors inside the same framework identifies when the KAN factors pay off, and a cost comparison includes a classical body-fitted solver.
\end{enumerate}

The remainder of this paper is organized as follows. Section~\ref{sec:intro} reviews KANs and PIKANs, including the univariate basis families used later. Section~\ref{sec:transform} presents the coordinate transformation, the autograd-based metric evaluation and the transformed model equations. Section~\ref{sec:ctpikan} formulates the point-wise CT-PIKAN and discusses its limitations. Section~\ref{sec:method} introduces CS-PIKAN, and Section~\ref{sec:theory} develops the error analysis. Section~\ref{sec:results} reports the numerical experiments, Section~\ref{sec:discussion} discusses the findings, Section~\ref{sec:limits} lists the limitations, and Section~\ref{sec:conc} concludes. The appendices give the closed-form metrics of all maps, the finite-difference reference scheme, and the results of every individual run.

\subsection{Kolmogorov-Arnold representation and KAN layers}
Physics-Informed Kolmogorov-Arnold Networks (PIKANs) have recently demonstrated their efficiency as a PDE-solving paradigm by leveraging the Kolmogorov-Arnold representation theorem (KART). The theorem states that any continuous multivariate function $f:[0,1]^n\to\R$ can be written as Eq.~\eqref{eq:1}:
\begin{equation}
	f(\mathbf{X}) = f(x_1, x_2, \dots, x_n) = \sum_{j=1}^{2n+1} \psi_j \left( \sum_{i=1}^{n} \phi_{ij} (x_i) \right),
	\label{eq:1}
\end{equation}
where $\psi_j$ and $\phi_{ij}$ are continuous univariate functions. This theorem suggests that high-dimensional function learning can be handled by decomposing the target function into a hierarchy of univariate transformations. For the $k$-th layer of a Kolmogorov-Arnold Network (KAN), the univariate functions are collected in the matrix of Eq.~\eqref{eq:2}:
\begin{equation}
	\Phi^{k} =
	\begin{bmatrix}
		\phi_{11} & \phi_{12} & \cdots & \phi_{1n} \\
		\phi_{21} & \phi_{22} & \cdots & \phi_{2n} \\
		\vdots & \vdots & \ddots & \vdots \\
		\phi_{m1} & \phi_{m2} & \cdots & \phi_{mn}
	\end{bmatrix},
	\qquad
	x^{k+1}_j=\sum_{i=1}^{n}\phi_{ji}\big(x^{k}_i\big).
	\label{eq:2}
\end{equation}
Therefore, a full single-output KAN with $L$ layers is expressed as Eq.~\eqref{eq:3}:
\begin{equation}
	f(x_1, x_2, \dots, x_n) =
	\sum_{i_{L}=1}^{n_{L}} \phi^{L}_{i_{L+1},i_{L}}
	\left( \sum_{i_{L-1}=1}^{n_{L-1}} \dots
	\left( \sum_{i_2=1}^{n_2} \phi^{2}_{i_3,i_2}
	\left( \sum_{i_1=1}^{n_1} \phi^{1}_{i_2,i_1}(x_{i_1}) \right) \right) \right).
	\label{eq:3}
\end{equation}
In conventional networks such as multi-layer perceptrons (MLPs), fixed activation functions are applied at the nodes. In KANs, the learnable functions sit on the edges: each edge, which corresponds to a weight connection in a traditional network, carries its own univariate function, which allows a more adaptive representation.

\subsection{Univariate basis families}
\label{sec:bases}
The edge functions are expanded in a basis. Four families are used in this paper.

\paragraph{B-splines.} In the original KAN \cite{KAN1}, each edge function is a weighted combination of a fixed base function and a B-spline expansion, Eq.~\eqref{eq:4}:
\begin{equation}
	\phi(x) = w_b\, b(x) + w_s \sum_{i} c_i B_i(x),
	\label{eq:4}
\end{equation}
where $b(x)$ is a predefined base function such as SiLU, $B_i(x)$ are B-spline basis functions on a grid, and $w_b, w_s, c_i$ are learnable parameters. The SiLU (Sigmoid Linear Unit), also known as the Swish activation, is defined in Eq.~\eqref{eq:5}:
\begin{equation}
	\text{SiLU}(x) = \frac{x}{1 + e^{-x}}.
	\label{eq:5}
\end{equation}
B-splines are local and well suited to function fitting, but cubic B-splines have piecewise-constant third derivatives, and their second derivatives vanish outside the grid support. Grid extension or grid updates are also needed when the inputs leave the grid range.

\paragraph{Gaussian radial basis functions.} A Gaussian RBF layer \cite{FastKAN} expands each input $z_i\in[0,1]$, rescaled to $\tilde z_i=2z_i-1$, on $K$ Gaussian bumps with learnable centres $\mu_{ik}$ and widths $s_{ik}=e^{\ell_{ik}}$:
\begin{equation}
\phi_{ik}(\tilde z_i)=\exp\!\left(-s_{ik}\,(\tilde z_i-\mu_{ik})^2\right),\qquad
\mathbf{h}^{(1)}=\tanh\!\Big(\textstyle\sum_{i}\sum_{k=1}^{K} \mathbf{W}_{ik}\,\phi_{ik}(\tilde z_i)+\mathbf{b}\Big).
\label{eq:rbf}
\end{equation}
The basis is $C^\infty$ and needs no grid extension. The CT-PIKAN baseline of Section~\ref{sec:ctpikan} uses this layer, with $K=16$, as its input layer, followed by a $\tanh$ trunk; only its first layer carries learnable univariate functions, so it is a \emph{hybrid} KAN.

\paragraph{Chebyshev polynomials.} Chebyshev-KANs \cite{ChebyKAN,Mos2} expand each edge function in Chebyshev polynomials $T_k$ of a bounded input, usually $\tanh$-normalised to $[-1,1]$:
\begin{equation}
\phi(z)=\sum_{k=0}^{K-1}c_k\,T_k\!\big(\tanh z\big),\qquad T_0=1,\quad T_1(s)=s,\quad T_{k+1}(s)=2sT_k(s)-T_{k-1}(s).
\label{eq:cheb}
\end{equation}
They are global, spectrally accurate for smooth functions, and their derivatives of any order are available by recursion or automatic differentiation.

\paragraph{Fourier modes.} For a coordinate $z$ that is periodic with period $1$, a natural basis is $\{1,\cos2\pi kz,\sin2\pi kz\}_{k=1}^{K_f}$. A function expanded in it is exactly periodic together with all of its derivatives. CS-PIKAN uses this basis in the first layer of every periodic coordinate (Section~\ref{sec:method}).

\subsection{Physics-informed training}
Given a PDE $\mathcal N[u]=0$ in $\Omega$ with boundary operator $\mathcal B[u]=g$ on $\partial\Omega$ and, for transient problems, an initial condition $u(\cdot,0)=u_0$, a PIKAN $u_\theta$ is trained by minimizing
\begin{equation}
\mathcal L(\theta)=\frac{\lambda_r}{N_r}\sum_{i=1}^{N_r}\big|\mathcal N[u_\theta](\mathbf x_i)\big|^2+\frac{\lambda_b}{N_b}\sum_{j=1}^{N_b}\big|\mathcal B[u_\theta](\mathbf x_j)-g(\mathbf x_j)\big|^2+\frac{\lambda_0}{N_0}\sum_{k=1}^{N_0}\big|u_\theta(\mathbf x_k,0)-u_0(\mathbf x_k)\big|^2,
\label{eq:pinnloss}
\end{equation}
where the derivatives in $\mathcal N$ and $\mathcal B$ are evaluated by automatic differentiation, and $\lambda_r,\lambda_b,\lambda_0$ are user-chosen weights. Problem-dependent choice of these weights, and the competition between the terms, is a well-known source of training difficulties \cite{PINN3}.

\subsection{Preliminary validation of the PIKAN building block}
\label{sec:prelim}
\begin{figure}[H]
	\centering
	\subfloat [PIKAN ($t=0$)]{\includegraphics[trim={0 0 0 1cm},clip,width=0.22\textwidth]{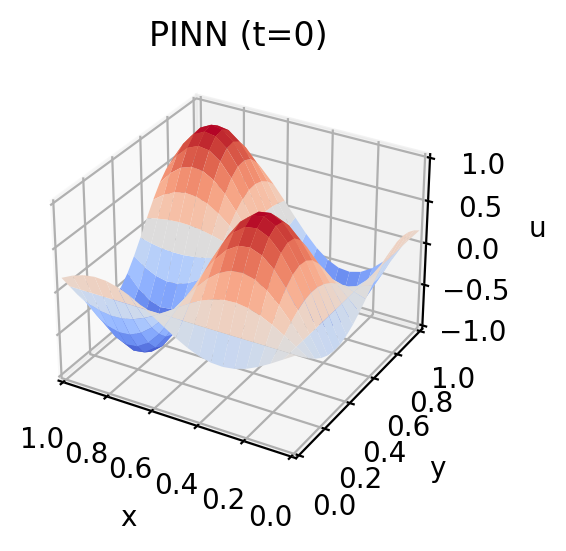}}
	\subfloat [PIKAN ($t=0.2$)]{\includegraphics[trim={0 0 0 1cm},clip,width=0.22\textwidth]{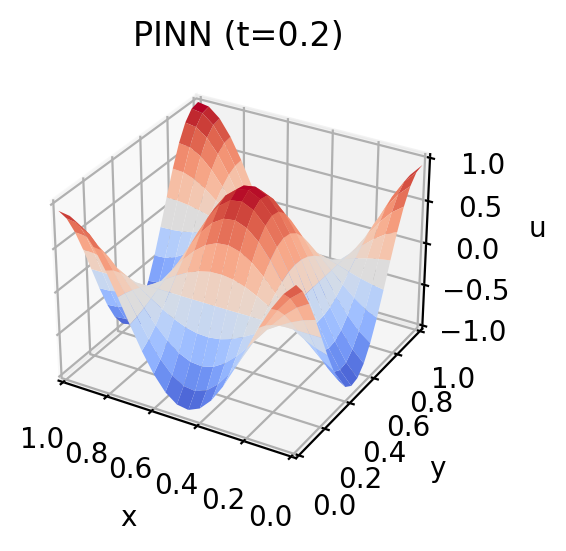}}
	\subfloat [Analytical ($t=0$)]{\includegraphics[trim={0 0 0 1cm},clip,width=0.22\textwidth]{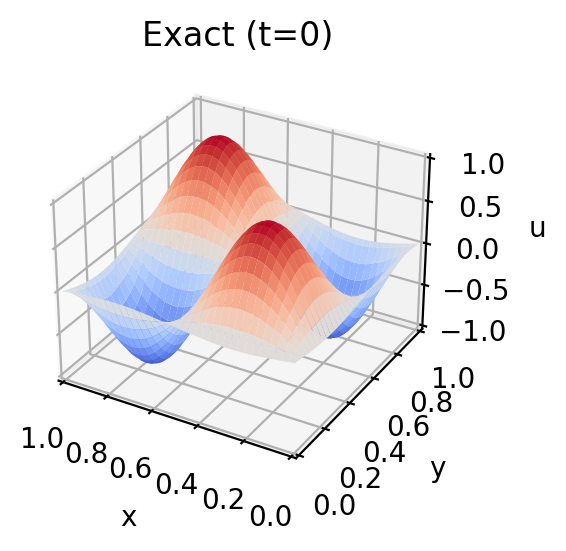}}
	\subfloat [Analytical ($t=0.2$)]{\includegraphics[trim={0 0 0 1cm},clip,width=0.22\textwidth]{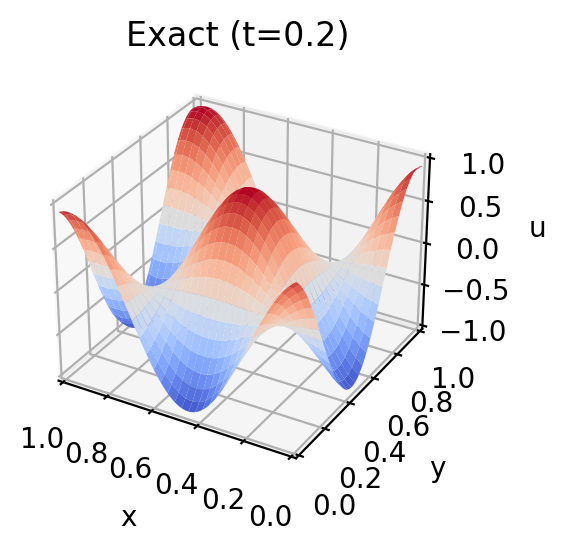}}
	\\
	\subfloat [PIKAN ($t=0$)]{\includegraphics[trim={0 0 0 0.65cm},clip,width=0.22\textwidth]{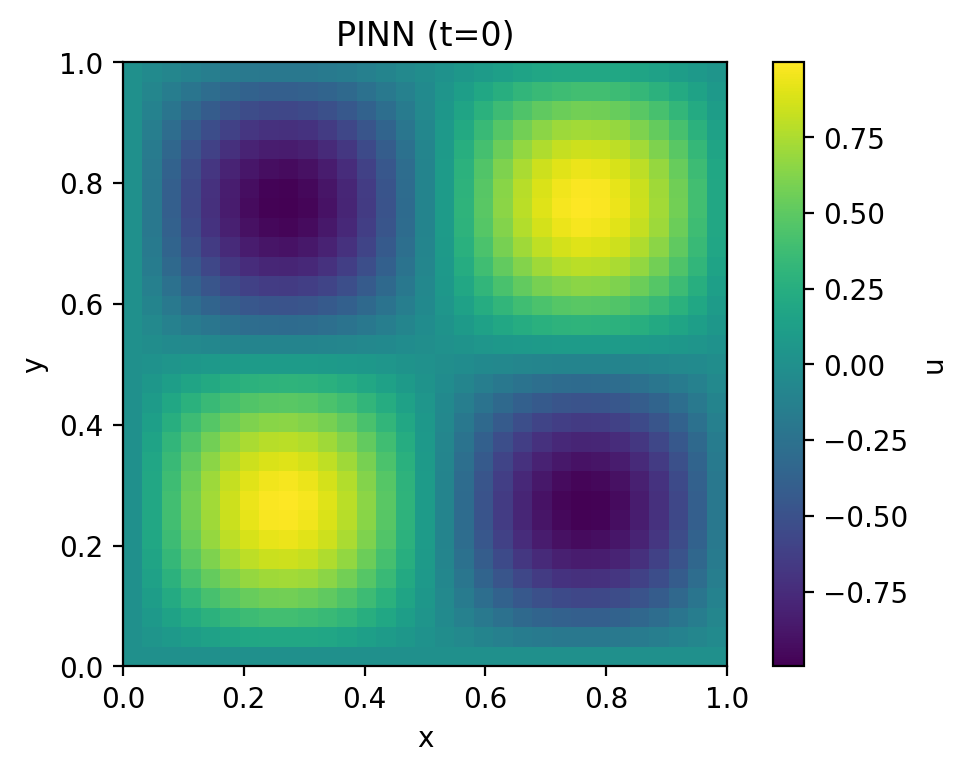}}
	\subfloat [PIKAN ($t=0.2$)]{\includegraphics[trim={0 0 0 0.65cm},clip,width=0.22\textwidth]{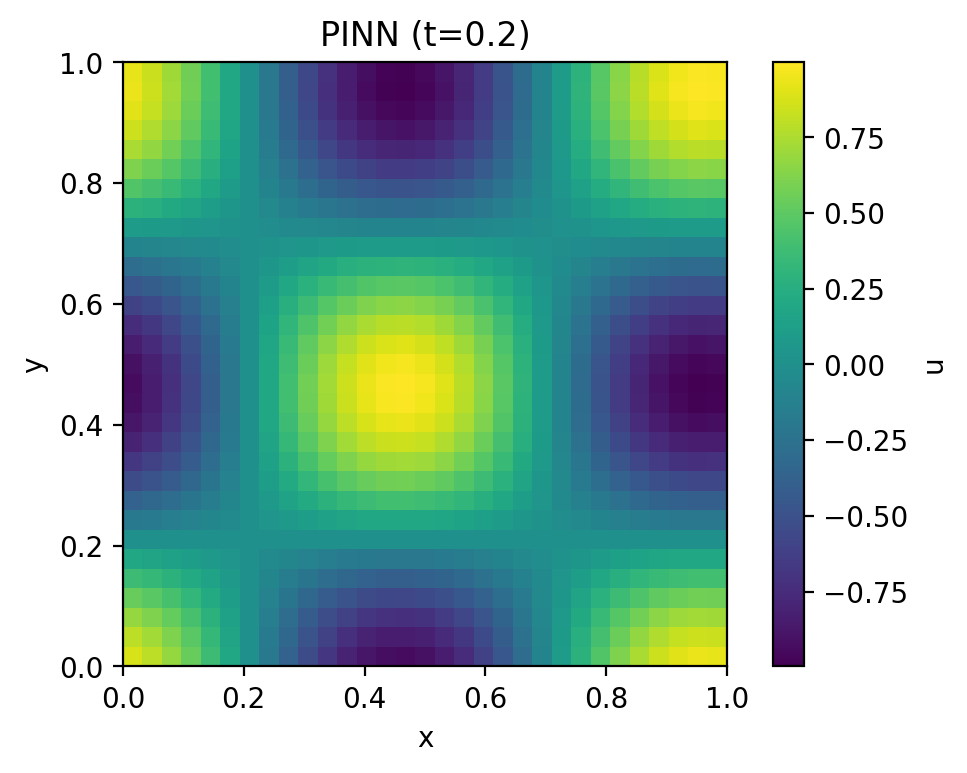}}
	\subfloat [Analytical ($t=0$)]{\includegraphics[trim={0 0 0 0.65cm},clip,width=0.22\textwidth]{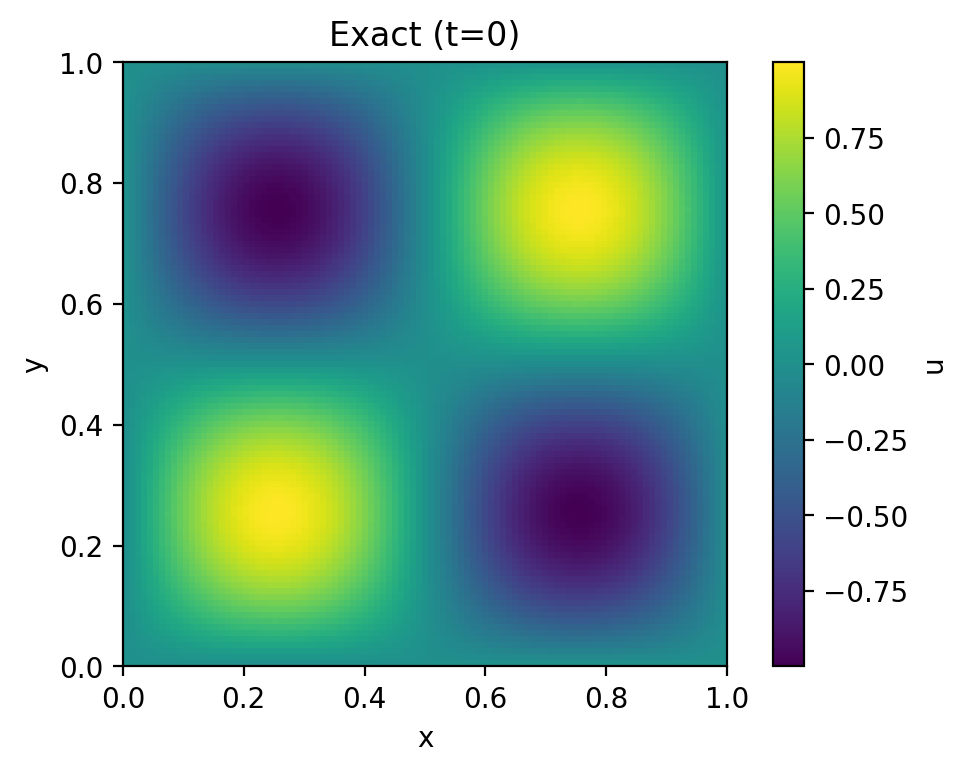}}
	\subfloat [Analytical ($t=0.2$)]{\includegraphics[trim={0 0 0 0.65cm},clip,width=0.22\textwidth]{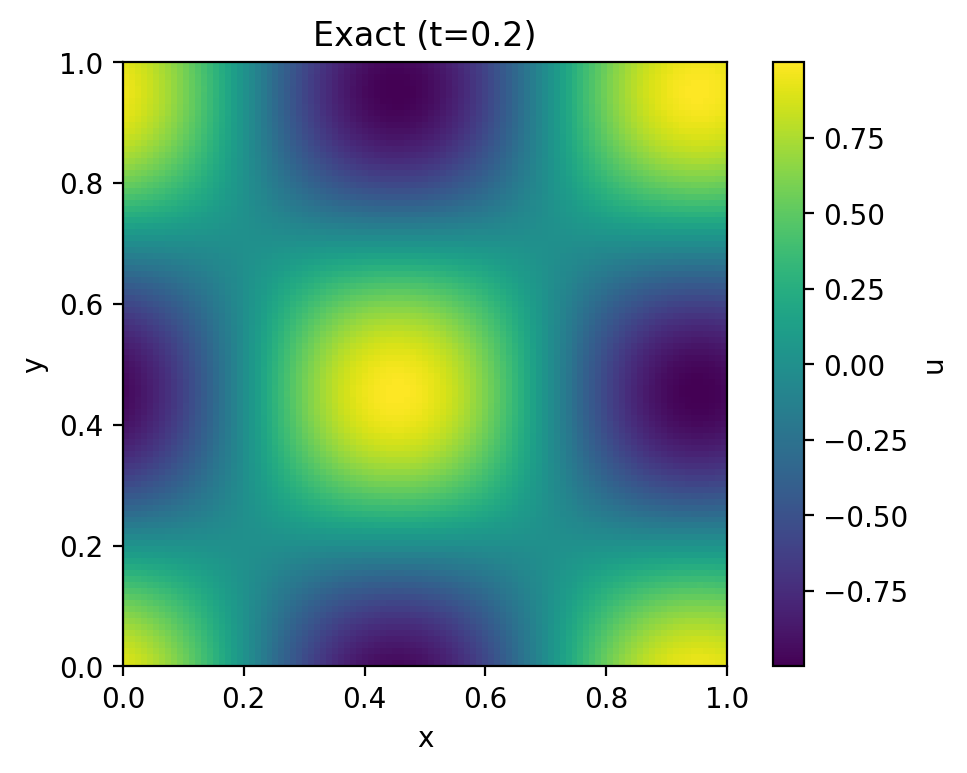}}
	\caption{Preliminary validation of the data-free B-spline PIKAN against the analytical solution of the Cartesian advection problem \eqref{eq:24}-\eqref{eq:25}.}
	\label{fig:3}
\end{figure}

Before introducing coordinate transformations, the data-free B-spline PIKAN of Eqs.~\eqref{eq:3}--\eqref{eq:5} was validated on the two-dimensional linear advection equation \cite{2Dadvection1,2Dadvection2,2Dadvection3,2Dadvection4}
\begin{equation}
u_{t} + c(u_{x} + u_{y}) = 0,
\label{eq:24}
\end{equation}
on the unit square with the sinusoidal initial condition
\begin{equation}
	u(x,y,0) = \sin(2\pi x)\sin(2\pi y).
	\label{eq:25}
\end{equation}
Fig.~\ref{fig:3} shows the predicted and analytical fields. This preliminary test only confirms that the KAN building block can represent a transported field. The quantitative comparison on this problem, with exact inflow data, three seeds and both CT-PIKAN and CS-PIKAN, is given in Section~\ref{sec:res_adv}.

\section{Coordinate transformation and geometric formulation}
\label{sec:transform}

\subsection{Curvilinear coordinates}
In many scientific and engineering applications, PDEs are posed on domains whose boundaries are curved or geometrically complex. Solving such equations directly in the physical domain can be challenging. A common approach, adopted in this work, is to introduce a smooth bijective coordinate transformation. Let $\Omega\subset\R^2$ be a physical domain with curvilinear geometry, and let $\Q=(0,1)^2$ be the computational square with coordinates $\bxi=(\xi,\eta)$. If $\Omega$ is doubly connected (an annulus or a star-shaped domain with a hole), $\eta$ is taken $1$-periodic, and $\Q$ is understood as $(0,1)\times\mathbb T$. We consider a smooth and invertible map
\begin{equation}
	\mathbf{x} = \Phi(\bxi) = \big(x(\xi,\eta),\, y(\xi,\eta)\big).
	\label{eq:6}
\end{equation}
The inverse map is assumed to exist and to be differentiable, although only the forward map is needed in practice. The Jacobian matrix of the transformation and its determinant are \cite{Transform1}
\begin{equation}
	\mathbf{J} = \frac{\partial \mathbf{x}}{\partial \bxi} =
	\begin{bmatrix}
		x_{\xi} & x_{\eta} \\
		y_{\xi} & y_{\eta}
	\end{bmatrix},
	\qquad
	J = \det (\mathbf{J}) = x_{\xi} y_{\eta} - x_{\eta} y_{\xi} > 0 .
	\label{eq:7}
\end{equation}
By the chain rule, the Cartesian derivatives are expressed in computational coordinates as
\begin{equation}
\frac{\partial}{\partial x} = \frac{1}{J} \left( y_{\eta}\frac{\partial}{\partial \xi} - y_{\xi}\frac{\partial}{\partial \eta} \right),
\qquad
\frac{\partial}{\partial y} = \frac{1}{J} \left( - x_{\eta}\frac{\partial}{\partial \xi} + x_{\xi}\frac{\partial}{\partial \eta} \right),
\label{eq:8}
\end{equation}
so that the inverse metrics are $\xi_x=y_\eta/J$, $\xi_y=-x_\eta/J$, $\eta_x=-y_\xi/J$ and $\eta_y=x_\xi/J$. The covariant metric tensor, the gradient transformation and the Laplacian in general curvilinear coordinates are given by Eqs.~\eqref{eq:9}--\eqref{eq:11}:
\begin{equation}
	G = \mathbf J^{\top} \mathbf J, \qquad g_{ij}=\frac{\partial\mathbf x}{\partial\xi_i}\cdot\frac{\partial\mathbf x}{\partial\xi_j},\qquad g=\det G = J^2,
	\label{eq:9}
\end{equation}
\begin{equation}
	\nabla_x u = \mathbf J^{-\top} \nabla_\xi \hat u,
	\label{eq:10}
\end{equation}
\begin{equation}
	\Delta_x u = \frac{1}{\sqrt{g}} \nabla_\xi \cdot \left( \sqrt{g}\, G^{-1} \nabla_\xi \hat u \right),
	\label{eq:11}
\end{equation}
where $\hat u=u\circ\Phi$ is the pulled-back field. Eq.~\eqref{eq:11} is the Laplace--Beltrami operator of the metric $G$; in index notation,
\begin{equation}
\Delta_x u = \frac{1}{\sqrt{g}} \frac{\partial}{\partial \xi_i} \left( \sqrt{g}\, g^{ij} \frac{\partial \hat u}{\partial \xi_j} \right)
= g^{ij}\,\frac{\partial^2\hat u}{\partial\xi_i\partial\xi_j}+b^j\,\frac{\partial\hat u}{\partial\xi_j},\qquad b^j=\frac{1}{\sqrt g}\frac{\partial}{\partial\xi_i}\big(\sqrt g\,g^{ij}\big),
\label{eq:12}
\end{equation}
with summation over $i,j\in\{1,2\}$, $(\xi_1,\xi_2)=(\xi,\eta)$, and $g^{ij}$ the entries of $G^{-1}$. The first form is conservative (divergence form); the second, non-divergence form, is the one used by CS-PIKAN. The coefficients $b^j$ contain second derivatives of the map. First-order (advection-type) terms transform as
\begin{equation}
\nabla_{x} u = \left( \hat u_{\xi}\, \xi_x + \hat u_{\eta}\, \eta_x ,\; \hat u_{\xi}\, \xi_y + \hat u_{\eta}\, \eta_y \right),
\qquad
\mathbf v\cdot\nabla_x u=\mathbf w\cdot\nabla_\xi\hat u,\qquad \mathbf w=\mathbf J^{-1}\mathbf v,
\label{eq:13}
\end{equation}
where $\mathbf w$ is the contravariant velocity. Two further quantities are used in the analysis of Section~\ref{sec:theory}:
\begin{equation}
M=\sqrt g\,G^{-1},\qquad \kappa=\frac{\sigma_1}{\sigma_2},
\label{eq:Mkappa}
\end{equation}
where $\sigma_1\ge\sigma_2>0$ are the singular values of $\mathbf J$. The matrix $M$ appears in the pulled-back Dirichlet energy, and $\kappa$ measures the local anisotropy of the map.

\subsection{Automatic differentiation for metric computation}
\label{sec:autograd}
Traditional coordinate-transformed solvers require manual analytical computation of the metric terms $x_{\xi}, x_{\eta}, y_{\xi}, y_{\eta}, g_{ij}, g^{ij}, J$ and $\sqrt{g}$. In contrast, the proposed method computes all these quantities automatically by differentiating the map:
\begin{equation}
x_{\xi} = \frac{\partial x}{\partial \xi},\qquad
x_{\eta} = \frac{\partial x}{\partial \eta},\qquad
y_{\xi} = \frac{\partial y}{\partial \xi},\qquad
y_{\eta} = \frac{\partial y}{\partial \eta}.
\label{eq:14}
\end{equation}
The map is defined once as a differentiable function, $\mathbf{x} = \Phi(\xi,\eta)$ (Eq.~\eqref{eq:6}), so that any first- or higher-order Cartesian derivative can be written in terms of $(\xi,\eta)$. The metric tensor, its determinant, the inverse metric and the Jacobian entries are
\begin{equation}
g_{ij} =
\begin{bmatrix}
	x_{\xi}^{2} + y_{\xi}^{2} & x_{\xi}x_{\eta} + y_{\xi}y_{\eta} \\
	x_{\xi}x_{\eta} + y_{\xi}y_{\eta} & x_{\eta}^{2} + y_{\eta}^{2}
\end{bmatrix},
\qquad
g = \det(g_{ij}), \qquad \sqrt{g} = \sqrt{\det(g_{ij})},
\label{eq:16}
\end{equation}
\begin{equation}
g^{ij} = (g_{ij})^{-1}, \qquad  J_{ij} = \frac{\partial x_i}{\partial \xi_j},
\label{eq:17}
\end{equation}
all obtained directly by automatic differentiation; the metric tensor itself is $G=\mathbf J^\top\mathbf J$. The coefficients $b^j$ of Eq.~\eqref{eq:12} require one further level of differentiation: the function $\bxi\mapsto M(\bxi)=\sqrt g\,G^{-1}$ is itself differentiated, and $b^j=\frac{1}{\sqrt g}\sum_i\partial_i M_{ij}$. Algorithm~\ref{al:1} summarises the procedure.

\begin{algorithm}[H]
	\caption{Autograd-based coordinate transformation and metric computation}
	\begin{algorithmic}[1]
		\Require Differentiable coordinate map $\Phi(\xi,\eta)$; evaluation points $\{(\xi_i,\eta_i)\}_{i=1}^{N}$
		\For{each point $(\xi_i,\eta_i)$}
		\State Physical coordinates: $(x_i, y_i) = \Phi(\xi_i, \eta_i)$
		\State Jacobian by automatic differentiation: $\mathbf J_i = \partial (x,y)/\partial (\xi,\eta)$
		\State Inverse Jacobian $\mathbf J_i^{-1}$; metric tensor $G_i = \mathbf J_i^\top \mathbf J_i$; $\sqrt{g_i}=|\det\mathbf J_i|$; $G_i^{-1}$
		\State $M_i=\sqrt{g_i}\,G_i^{-1}$ as a function of $\bxi$; differentiate once more: $b^j_i=\frac{1}{\sqrt{g_i}}\sum_k\partial_k M_{kj}$
		\State (Advection) contravariant velocity $\mathbf w_i=\mathbf J_i^{-1}\mathbf v(x_i,y_i)$
		\EndFor
		\State \Return $\{\mathbf J_i, \mathbf J_i^{-1}, G_i, \sqrt{g_i}, G_i^{-1}, b_i, \mathbf w_i\}$
	\end{algorithmic}
	\label{al:1}
\end{algorithm}

\paragraph{Requirements on the map.} The map must be a diffeomorphism onto $\Omega$ ($J>0$ everywhere), and it must be twice differentiable, because the Laplace--Beltrami operator differentiates $\sqrt g\,g^{ij}$ once more. We check $\min J$ on a dense grid before training (Table~\ref{tab:geom}).

\paragraph{Maps without a closed form.} Nothing in Algorithm~\ref{al:1} requires $\Phi$ to be an analytical expression; it only has to be a differentiable tensor function. In Section~\ref{sec:res_poisson} we construct the map numerically from boundary curves that are known only as point samples. Each curve is represented by the trigonometric interpolant of its samples, and the curves are blended by linear transfinite interpolation \cite{Gordon1973,g4},
\begin{equation}
\Phi(\xi,\eta)=(1-\xi)\,\mathbf c_{\rm in}(\eta)+\xi\,\mathbf c_{\rm out}(\eta).
\label{eq:tfi_map}
\end{equation}
The metrics are again obtained by automatic differentiation. A map represented by a trained network, as in \cite{JacobiNet,Khairullin2025}, would enter in the same way.

\paragraph{Verification.} For every closed-form map used in this paper, Appendix~\ref{app:metrics} lists the analytical metrics. In double precision the autograd values agree with them to within $2.3\times10^{-13}$ on $2000$ random points per map.

\subsection{Transformed model equations}
\label{sec:modeleqs}
A general second-order PDE in the physical domain,
\begin{equation}
	\mathcal{F}\big(u, \partial_x u, \partial_y u, \partial_{xx}u, \partial_{xy}u, \partial_{yy}u \big) = 0,
	\qquad (x,y)\in\Omega,
	\label{eq:19}
\end{equation}
with suitable boundary conditions (Dirichlet, Neumann, or Robin) and, for time-dependent problems, an initial condition, becomes
\begin{equation}
\widetilde{\mathcal{F}}\big(\hat u, \hat u_{\xi}, \hat u_{\eta}, \hat u_{\xi\xi}, \hat u_{\xi\eta}, \hat u_{\eta\eta}; \Phi_{\xi},\Phi_{\eta},\Phi_{\xi\xi},\Phi_{\xi\eta},\Phi_{\eta\eta}\big)=0,
\qquad (\xi,\eta) \in \Q
\label{eq:22}
\end{equation}
after substitution of Eqs.~\eqref{eq:8}--\eqref{eq:13}. The three model equations of this paper take the following forms.

\paragraph{Poisson equation.} $-\Delta_x u=f$ becomes
\begin{equation}
-\Big(g^{11}\hat u_{\xi\xi}+2g^{12}\hat u_{\xi\eta}+g^{22}\hat u_{\eta\eta}+b^1\hat u_\xi+b^2\hat u_\eta\Big)=\hat f(\xi,\eta),\qquad \hat f=f\circ\Phi.
\label{eq:poisson_t}
\end{equation}

\paragraph{Heat equation.} $u_t=\alpha\Delta_x u+f$ becomes
\begin{equation}
\hat u_t=\alpha\Big(g^{11}\hat u_{\xi\xi}+2g^{12}\hat u_{\xi\eta}+g^{22}\hat u_{\eta\eta}+b^1\hat u_\xi+b^2\hat u_\eta\Big)+\hat f.
\label{eq:heat_t}
\end{equation}

\paragraph{Advection equation.} $u_t+\mathbf v\cdot\nabla_x u=0$ becomes the advection equation in computational coordinates with the contravariant velocity,
\begin{equation}
\hat u_t+w^1\hat u_\xi+w^2\hat u_\eta=0,\qquad (w^1,w^2)^\top=\mathbf J^{-1}\mathbf v .
\label{eq:adv_t}
\end{equation}
For the constant velocity $\mathbf v=c(1,1)$ of Eq.~\eqref{eq:24}, substitution of Eqs.~\eqref{eq:33}--\eqref{eq:34},
\begin{equation}
\frac{\partial}{\partial x} = \frac{\partial}{\partial\xi}\,\frac{\partial\xi}{\partial x} + \frac{\partial}{\partial\eta}\,\frac{\partial\eta}{\partial x},
\label{eq:33}
\end{equation}
\begin{equation}
\frac{\partial}{\partial y} = \frac{\partial}{\partial\xi}\,\frac{\partial\xi}{\partial y} + \frac{\partial}{\partial\eta}\,\frac{\partial\eta}{\partial y},
\label{eq:34}
\end{equation}
gives
\begin{equation}
u_{t} + c\left[ (\xi_{x} + \xi_{y})u_{\xi} + (\eta_{x} + \eta_{y})u_{\eta} \right] = 0,
\label{eq:39}
\end{equation}
that is, $w^1=c(\xi_x+\xi_y)$ and $w^2=c(\eta_x+\eta_y)$. For a general map these coefficients \emph{vary in space}. For the wavy map of Section~\ref{sec:res_adv},
\begin{equation}
x=\xi+A\sin(6\pi\eta),\qquad y=\eta+A\sin(6\pi\xi),
\label{eq:30}
\end{equation}
the metrics are
\begin{equation}
\begin{aligned}
&x_\xi=1,\quad x_\eta=6\pi A\cos(6\pi\eta),\quad y_\xi=6\pi A\cos(6\pi\xi),\quad y_\eta=1,\\
&J=1-36\pi^2A^2\cos(6\pi\xi)\cos(6\pi\eta),
\end{aligned}
\label{eq:wavy_metrics}
\end{equation}
such that Eq.~\eqref{eq:39} reads
\begin{equation}
u_{t} + c\left( a(\xi,\eta)\,u_{\xi} + b(\xi,\eta)\,u_{\eta} \right) = 0,\qquad
a=\frac{1-6\pi A\cos(6\pi\eta)}{J},\quad b=\frac{1-6\pi A\cos(6\pi\xi)}{J}.
\label{eq:40}
\end{equation}
For $A=0.02$, $J\in[0.858,1.142]$, and the map is a diffeomorphism. In both solvers of this paper, $a$ and $b$ are evaluated pointwise from the autograd metrics; Eq.~\eqref{eq:wavy_metrics} is used only for verification.

\paragraph{Advection in polar coordinates.} When the domain is a circular sector, polar coordinates $(r,\theta)$ give the natural form \cite{2Dpolar}
\begin{equation}
\frac{\partial u}{\partial t} + v_{r} \frac{\partial u}{\partial r} + \frac{v_{\theta}}{r} \frac{\partial u}{\partial \theta} = 0,
\label{eq:26}
\end{equation}
where $v_r$ and $v_\theta$ are the radial and angular velocity components. With $r=r_0+(r_1-r_0)\xi$ and $\theta=\theta_1\eta$ this is again Eq.~\eqref{eq:adv_t}, with $w^1=v_r/(r_1-r_0)$ and $w^2=v_\theta/(r\,\theta_1)$. For $v_r=v_\theta=1$ the characteristics satisfy $\dot r=1$ and $\dot\theta=1/r$, so $r-t$ and $\theta-\ln r$ are invariant along them. For the initial condition
\begin{equation}
	u(r,\theta,0) = \alpha \sin(\beta \pi r) \sin(\lambda \theta),
	\label{eq:27}
\end{equation}
this gives the exact solution
\begin{equation}
u(r,\theta,t)=\alpha\,\sin\!\big(\beta\pi(r-t)\big)\,\sin\!\Big(\lambda\big(\theta-\ln\tfrac{r}{r-t}\big)\Big),\qquad r-t>0,
\label{eq:exact_polar}
\end{equation}
which provides both the inflow data and the reference solution in Section~\ref{sec:res_adv}.

\section{The point-wise baseline: CT-PIKAN}
\label{sec:ctpikan}

\subsection{Formulation}
\begin{figure}[H]
	\centering
	\includegraphics[width=1\textwidth]{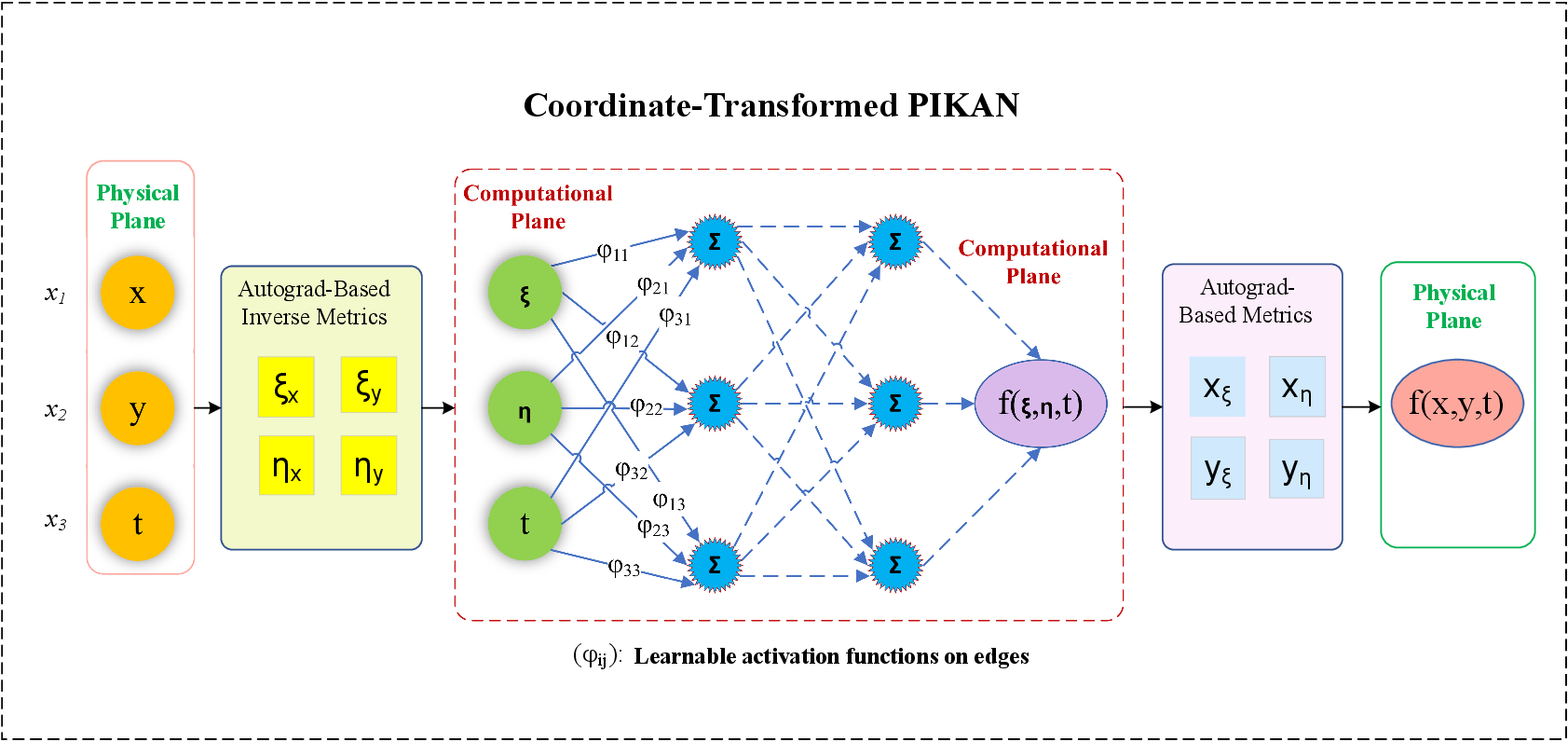}
	\caption{Schematic of the point-wise CT-PIKAN architecture used as baseline.}
	\label{fig:1}
\end{figure}

Physics-informed KANs are inherently mesh-free, yet their accuracy and optimization can degrade when the physical domain exhibits strong curvature or geometric irregularity. In the Coordinate-Transformed PIKAN (CT-PIKAN), the physical domain $\Omega$ is therefore mapped onto the computational domain $\Q$, the PDE is solved there, and the solution is evaluated on the physical domain through the forward map. The framework integrates three components: (i) the coordinate transformation of the curvilinear domain into the computational domain, (ii) automatic differentiation for geometry-aware metric computation, and (iii) a PIKAN for the solution. Fig.~\ref{fig:1} presents a schematic of the architecture.

\begin{figure}[H]
	\centering
	\includegraphics[width=0.86\textwidth]{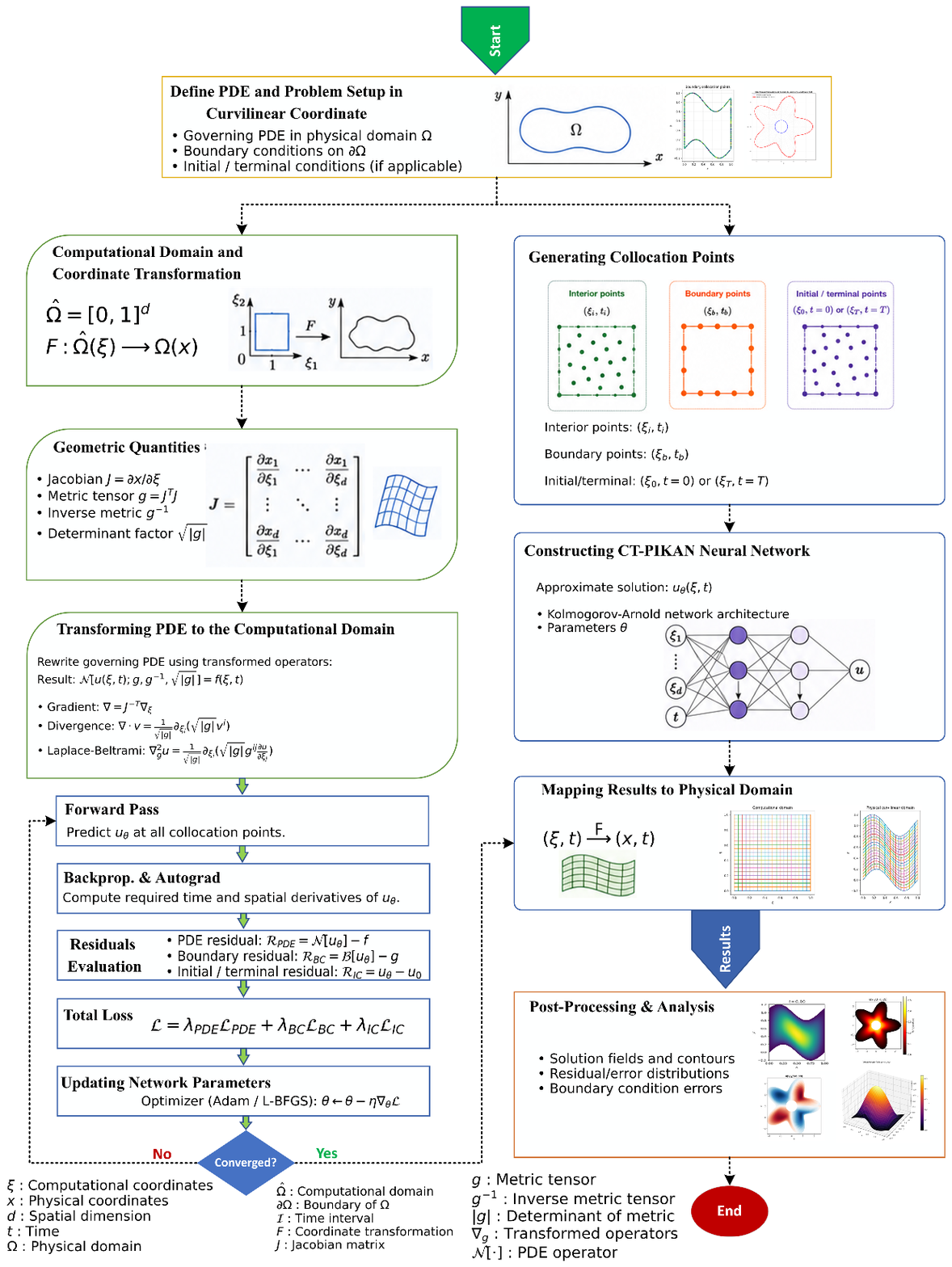}
	\caption{The training process of the point-wise CT-PIKAN.}
	\label{fig:2}
\end{figure}

Within CT-PIKAN, a network $\hat u_\theta(\xi,\eta[,t])$ approximates the pulled-back solution directly in the computational domain. The network receives $(\xi,\eta)$ (and $t$) as inputs. Automatic differentiation provides $\hat u_{\xi}$, $\hat u_{\eta}$ and the higher-order terms, which are substituted into the transformed operator $\widetilde{\mathcal{F}}$ of Eq.~\eqref{eq:22}. The PDE-residual loss, the boundary loss and (if needed) the initial-condition loss are aggregated into
\begin{equation}
	\mathcal{L}_{\text{CT-PIKAN}} = \lambda_r\mathcal{L}_{\text{PDE}} + \lambda_{b}\mathcal{L}_{\text{BC}} + \lambda_{0}\mathcal{L}_{\text{IC}},
	\label{eq:23}
\end{equation}
which is Eq.~\eqref{eq:pinnloss} evaluated at points of $\Q$. CT-PIKAN thus unifies (i) the coordinate transformation, (ii) the metric-tensor construction and (iii) the residual computation, all through automatic differentiation. Training proceeds entirely in the computational domain, where every physical boundary is a coordinate line. The training procedure is depicted in Fig.~\ref{fig:2}, and Algorithm~\ref{al:2} lists the steps.

\begin{algorithm}[H]
	\caption{CT-PIKAN physics-informed training procedure (point-wise baseline)}
	\begin{algorithmic}[1]
		\Require PIKAN model $u_\theta$; metric routine (Algorithm~\ref{al:1}); PDE operator $\mathcal{N}$; boundary operator $\mathcal{B}$; weights $\lambda_r, \lambda_b, \lambda_0$
		\State Initialize the network parameters $\theta$
		\For{each training iteration}
		\State Sample interior points $\{(\xi_i,\eta_i[,t_i])\}$, boundary points $\{(\xi_j,\eta_j[,t_j])\}$ and initial points $\{(\xi_k,\eta_k)\}$ at random
		\State Forward prediction $u_i = u_\theta(\xi_i,\eta_i[,t_i])$; metrics at the points by Algorithm~\ref{al:1}
		\State Gradients in computational space $\nabla_\xi u_i$; Hessian if needed; $\nabla_x u_i = \mathbf J_i^{-\top} \nabla_\xi u_i$; $\Delta_x u_i$ by Eq.~\eqref{eq:12}
		\State Residual $r_i = \mathcal{N}[u_\theta](\xi_i,\eta_i[,t_i])$; $\mathcal{L}_r = \frac{1}{N_r} \sum_i |r_i|^2$
		\State Boundary loss $\mathcal{L}_b = \frac{1}{N_b} \sum_j |\mathcal{B}(u_j) - g_j|^2$; initial loss $\mathcal L_0=\frac{1}{N_0}\sum_k|u_\theta(\xi_k,\eta_k,0)-\hat u_0(\xi_k,\eta_k)|^2$
		\State Total loss $\mathcal{L} = \lambda_r \mathcal{L}_r + \lambda_b \mathcal{L}_b+\lambda_0\mathcal L_0$; Adam update $\theta \leftarrow \theta - \alpha \nabla_\theta \mathcal{L}$
		\EndFor
		\State \Return trained CT-PIKAN model $u_\theta$
	\end{algorithmic}
	\label{al:2}
\end{algorithm}

\paragraph{Architecture and training.} The network used in this paper is the one of the public CT-PIKAN code: a Gaussian-RBF KAN input layer (Eq.~\eqref{eq:rbf}, $16$ learnable bases per input) followed by a $\tanh$ trunk of width $64$. Training uses Adam with learning rate $10^{-3}$, penalty weights $\lambda_b=\lambda_0=10$, and fresh random collocation points at every epoch. Table~\ref{tab:config} lists the per-case settings.

\subsection{Collocation strategies}
\begin{figure}[H]
	\centering
	\subfloat [$65 \times 10^{2}$ random c.p.]{\includegraphics[width=0.28\textwidth]{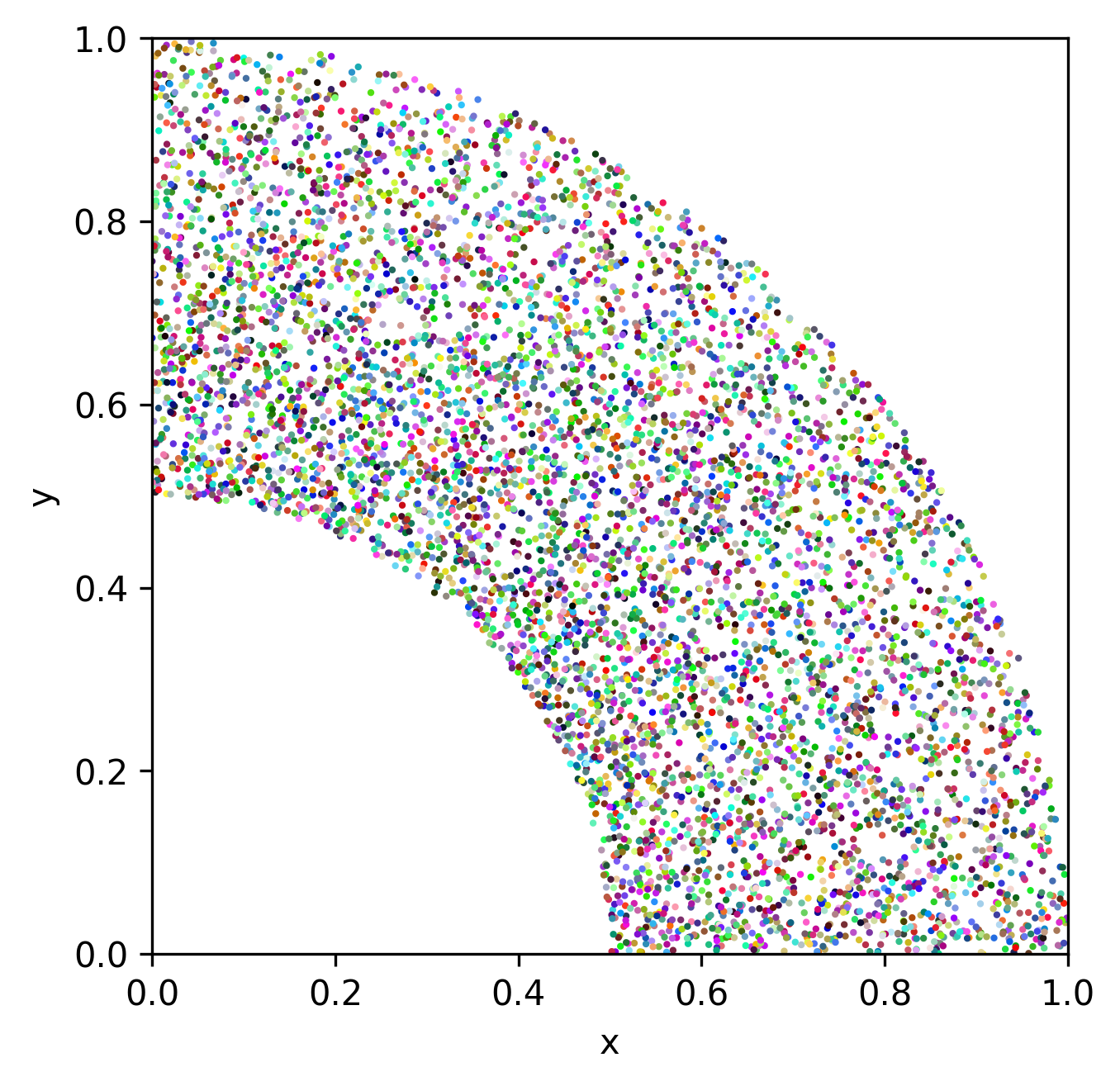}}
	\subfloat [$5 \times 10^{5}$ random c.p.]{\includegraphics[width=0.28\textwidth]{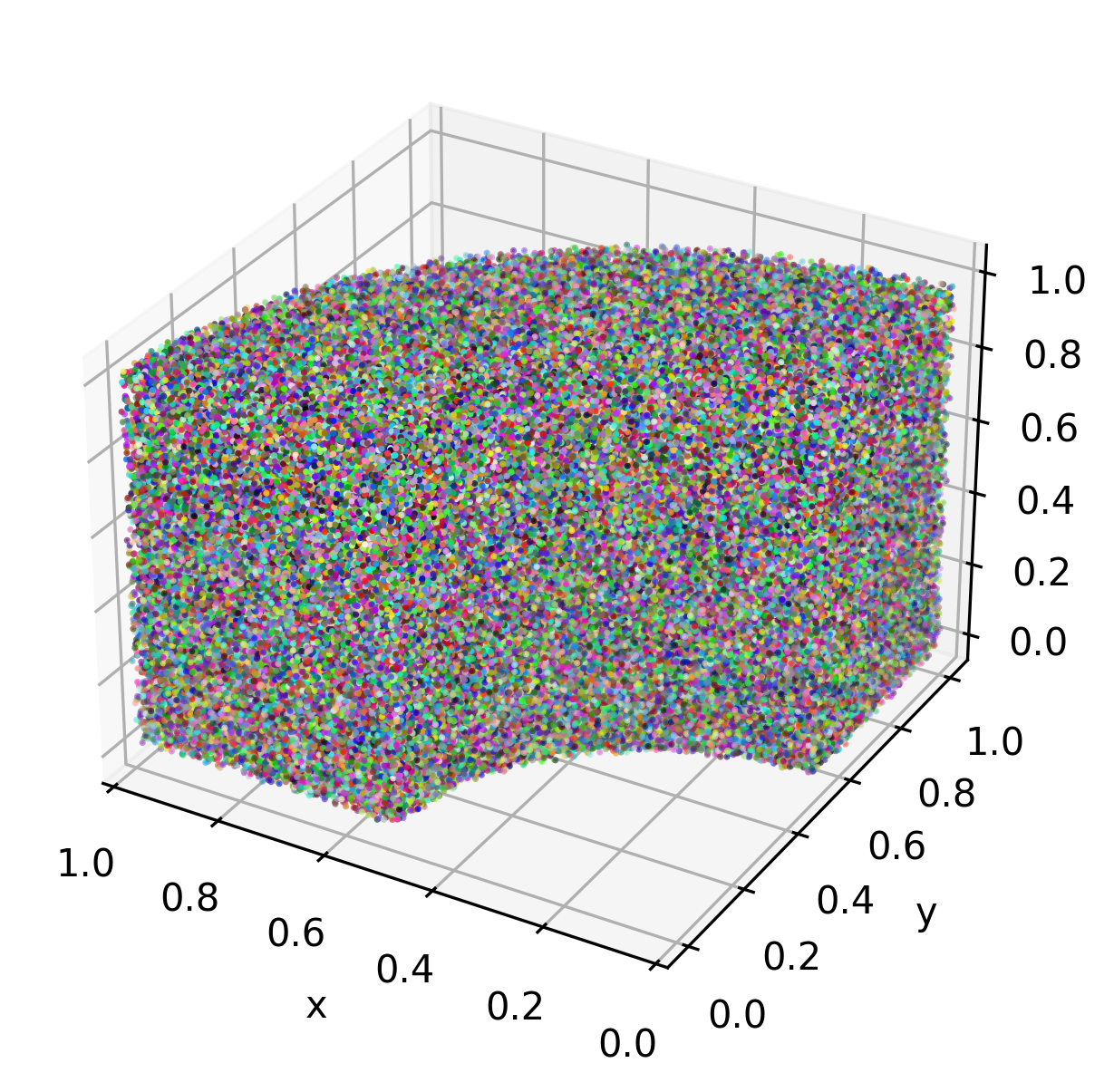}}
	\subfloat [$5 \times 10^{5}$ uniform c.p.]{\includegraphics[width=0.28\textwidth]{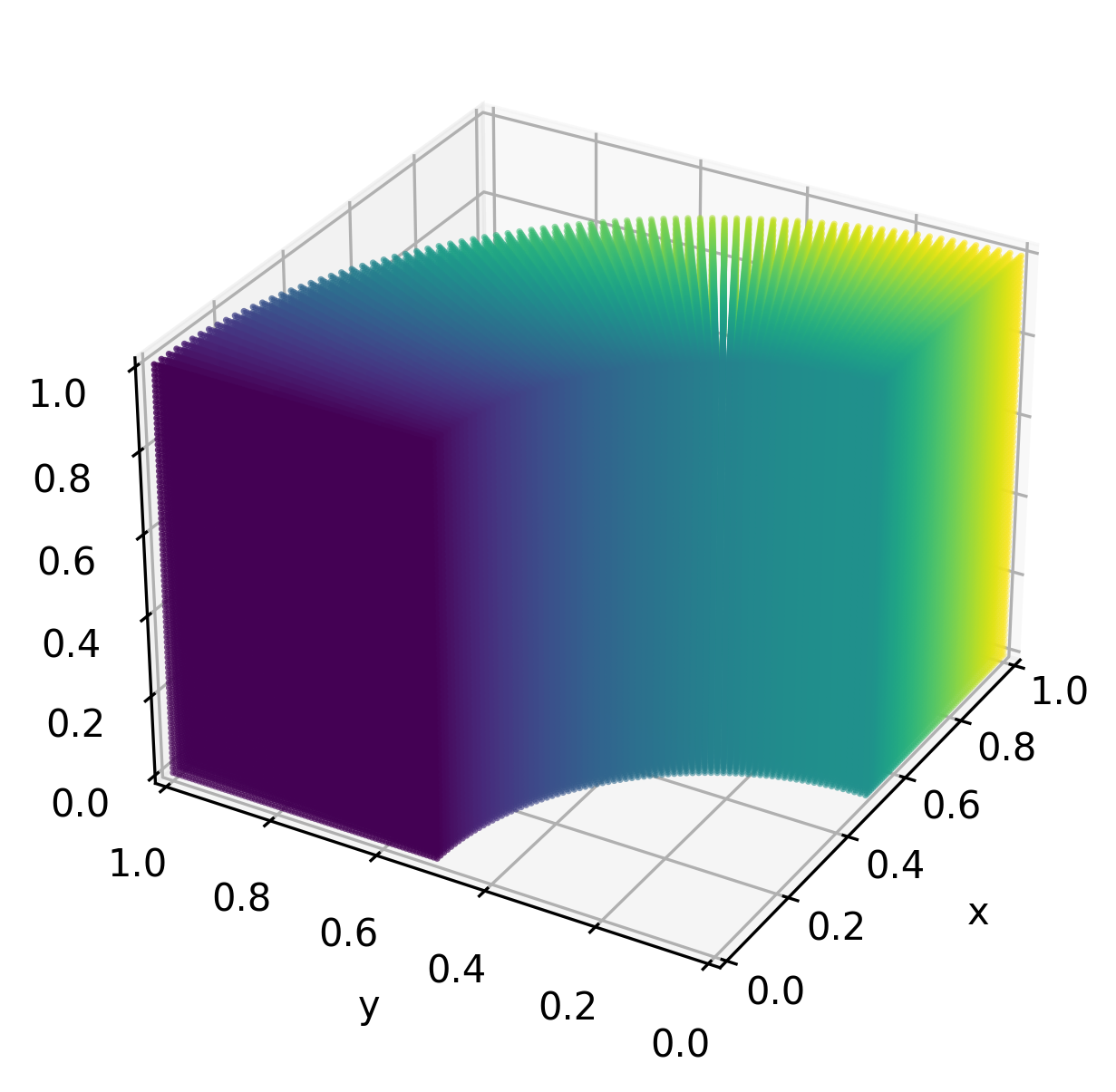}}
	\caption{Point-wise collocation strategies of the original study: spatial and spatio-temporal distribution of collocation points for the advection problem on the quarter-circle sector.}
	\label{fig:5}
\end{figure}

\begin{figure}[H]
	\centering
	\subfloat [Spatial collocation points ($10^{4}$)]{\includegraphics[width=0.35\textwidth, height=.35\textwidth]{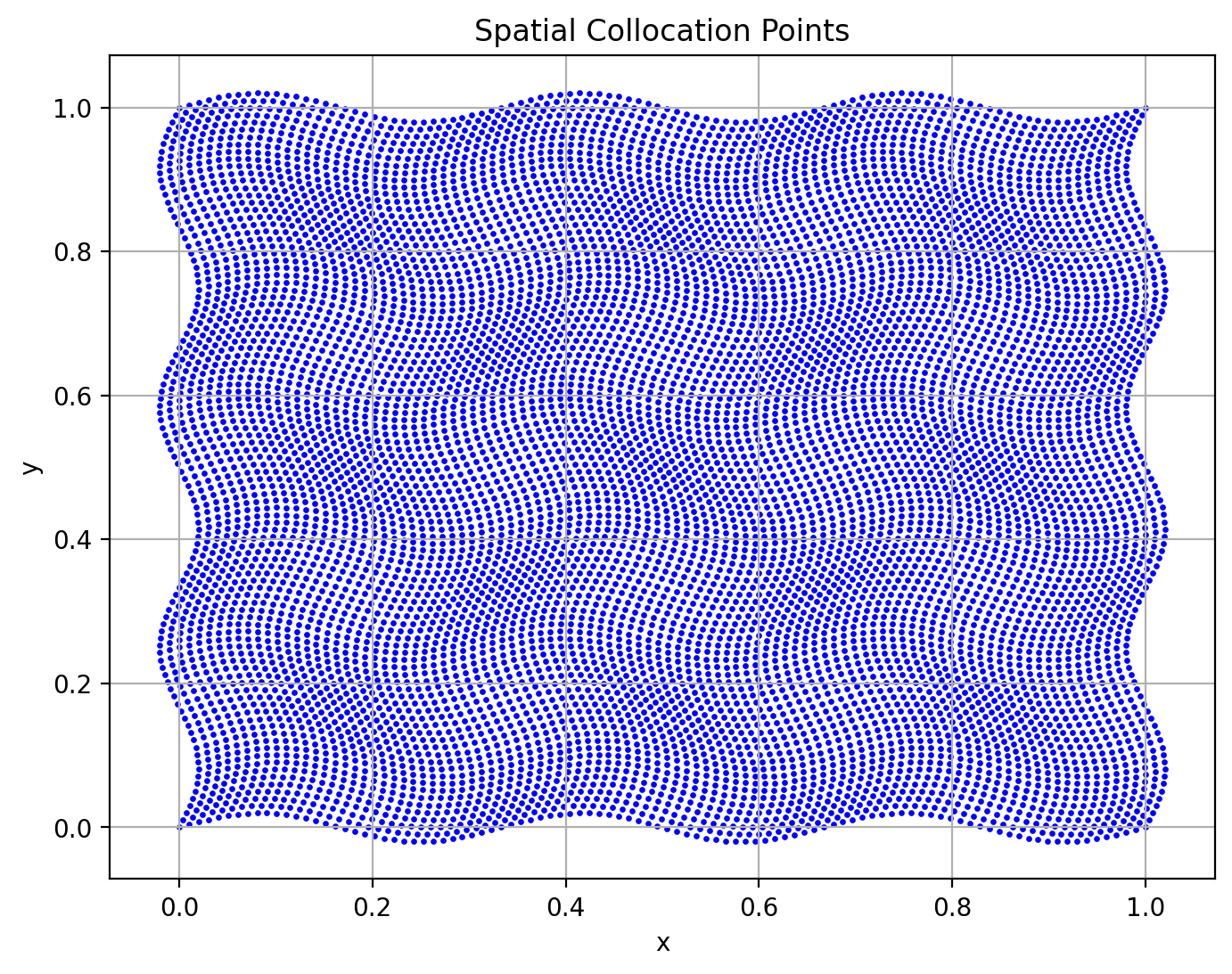}}
	\hspace{1cm}
	\subfloat [Spatio-temporal collocation points ($5\times10^{5}$)]{\includegraphics[width=0.35\textwidth, trim=0cm 0cm 2cm 2cm, clip]{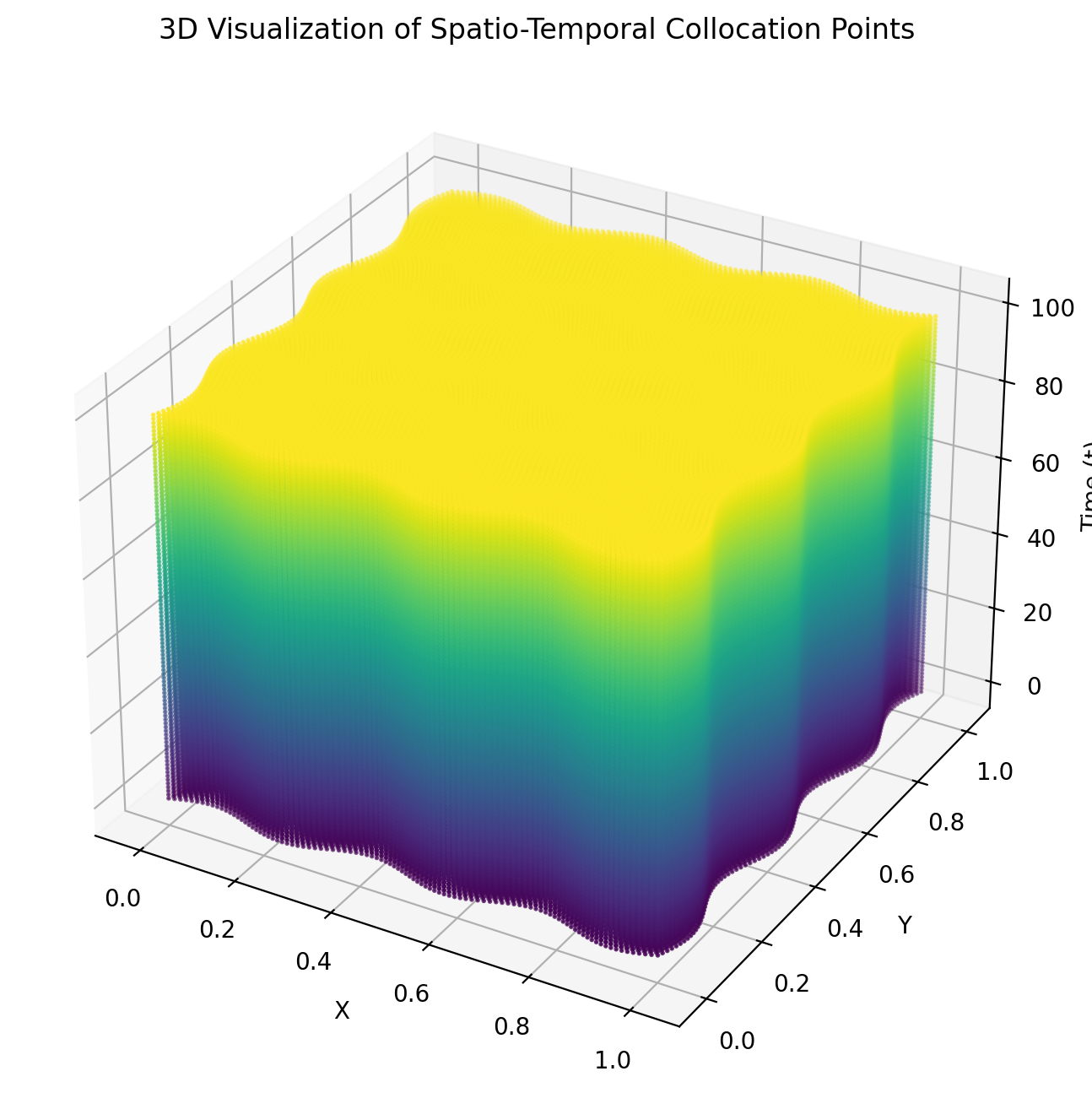}}
	\caption{Point-wise collocation points for the curvilinear (wavy) advection domain.}
	\label{fig:6}
\end{figure}

\begin{figure}[H]
	\centering
	\subfloat [Interior points]{\includegraphics[width=0.28\textwidth]{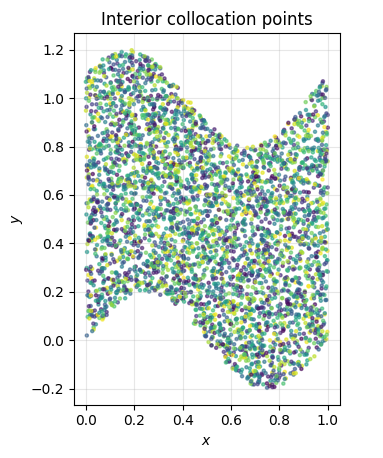}}
	\subfloat [B.C.~points]{\includegraphics[width=0.293\textwidth]{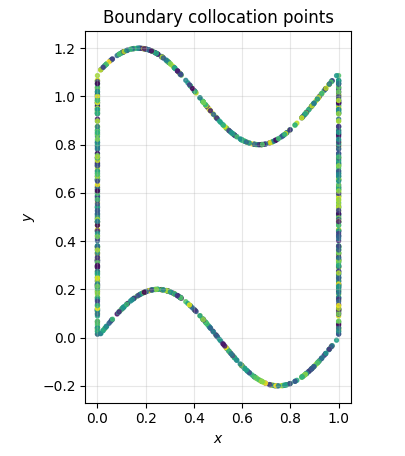}}
	\subfloat [I.C.~points]{\includegraphics[width=0.262\textwidth]{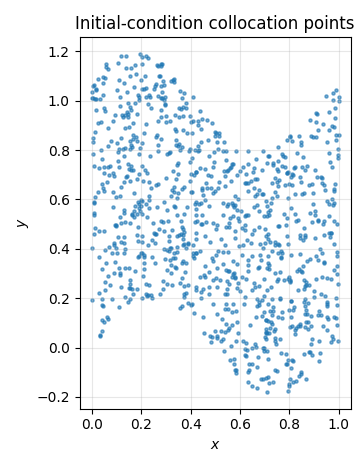}}
	\\
	\subfloat [Interior points]{\includegraphics[width=0.3\textwidth]{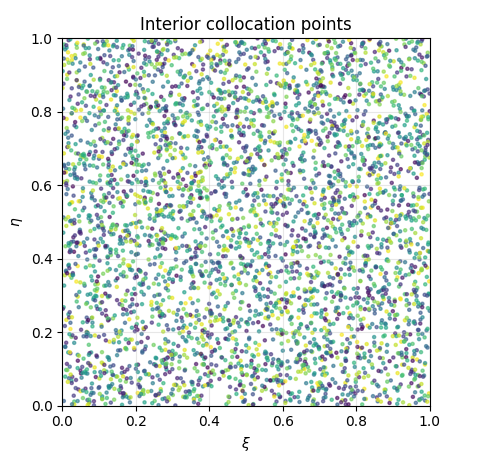}}
	\subfloat [B.C.~points]{\includegraphics[width=0.3\textwidth]{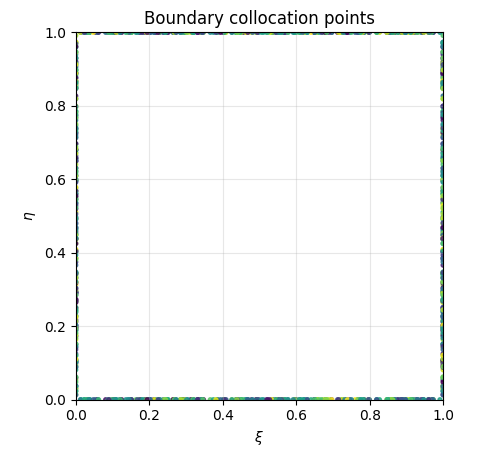}}
	\subfloat [I.C.~points]{\includegraphics[width=0.29\textwidth]{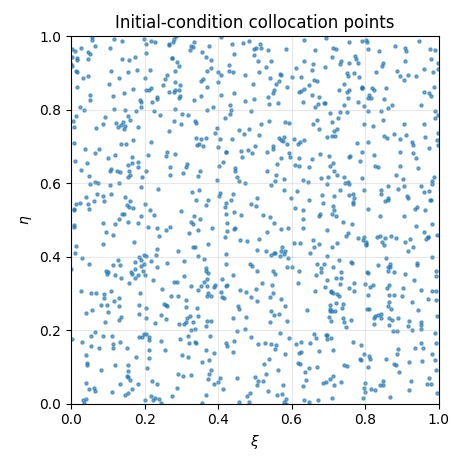}}
	\caption{Point-wise collocation points of CT-PIKAN in the physical curvilinear (wavy-channel) domain (top) and in the computational domain (bottom) for the transient heat equation.}
	\label{fig:11}
\end{figure}

Because CT-PIKAN is mesh-free, its collocation points can be generated in several ways. Figs.~\ref{fig:5}, \ref{fig:6} and \ref{fig:11} show the strategies used in the original study for the polar sector, the wavy advection domain and the wavy heat channel: random and uniform points in space--time, and the images of random computational points in the physical domain. Uniform sampling in $(\xi,\eta)$ is \emph{not} uniform in physical area: the physical density of the points is proportional to $1/\sqrt g$, which can be seen in the physical-domain panels.

\subsection{Limitations of point-wise training on transformed domains}
\label{sec:ct_limits}
CT-PIKAN removes the need to derive metric terms by hand, but it leaves three structural problems untouched, and our experiments (Section~\ref{sec:results}) show that they dominate its error.
\begin{enumerate}
\item \textbf{Penalised data.} The boundary and initial conditions are only approximately satisfied, and their weights compete with the residual. The diagnostic maps of Fig.~\ref{fig:15} (Section~\ref{sec:res_heat}) show boundary errors of the same order as the residual.
\item \textbf{The seam of closed coordinates.} For annular and star-shaped domains, $\theta=2\pi\eta$, so the edges $\eta=0$ and $\eta=1$ of the computational square are the same physical curve. A network that receives $\eta$ as a raw input is not periodic, and nothing in the loss forces $u$ and its derivatives to match across this seam. Section~\ref{sec:res_poisson} shows that the CT-PIKAN error concentrates exactly there.
\item \textbf{Cost.} Every collocation point requires nested automatic differentiation through the full network, including second derivatives for the Laplace--Beltrami operator, so the cost grows with the number of points times the network size.
\end{enumerate}
CS-PIKAN, introduced next, addresses all three by changing the representation of the solution rather than only the operator.

\section{CS-PIKAN: Curvilinear Separable Physics-Informed KAN}
\label{sec:method}

\subsection{Pull-back to the computational domain}
CS-PIKAN uses the same map $\Phi:\bar\Q\to\bar\Omega$ and the same autograd metrics as CT-PIKAN (Section~\ref{sec:transform}). We assume a $C^2$ diffeomorphism with $\det\mathbf J>0$, and use the quantities
\begin{equation}
G=\mathbf J^{\top}\mathbf J,\qquad \sqrt g=\det \mathbf J,\qquad M=\sqrt g\,G^{-1},\qquad b^j=\frac{1}{\sqrt g}\,\partial_{i}\!\left(\sqrt g\,g^{ij}\right),
\label{eq:metric}
\end{equation}
so that, for $\hat u=u\circ\Phi$,
\begin{align}
\nabla_x u &= \mathbf J^{-\top}\nabla_\xi \hat u,\qquad
\Delta_x u=g^{ij}\partial_{ij}\hat u+b^j\partial_j\hat u, \label{eq:LB}\\
\mathbf v\cdot\nabla_x u &= \mathbf w\cdot\nabla_\xi\hat u,\qquad \mathbf w=\mathbf J^{-1}\mathbf v. \label{eq:adv}
\end{align}
The difference from CT-PIKAN lies in how $\hat u$ is represented, how the data are imposed and how the loss is computed.

\subsection{Separable representation with coordinate-adapted KAN factors}
For a time-dependent problem, the solution on $\Q\times(0,T)$ is written as
\begin{equation}
\hat u_\theta(\xi,\eta,t)=\mathcal L(\xi,\eta,t)+\sum_{r=1}^{R} \underbrace{\big[e_\xi(\xi)\,a_r(\xi)\big]}_{A_r(\xi)}\,\underbrace{\big[e_\eta(\eta)\,b_r(\eta)\big]}_{B_r(\eta)}\,\underbrace{\big[e_t(t)\,c_r(t)\big]}_{C_r(t)},
\label{eq:ansatz}
\end{equation}
where $\mathbf a=(a_1,\dots,a_R)$, $\mathbf b$ and $\mathbf c$ are vector-valued one-dimensional KANs, $e_\xi,e_\eta,e_t$ are envelope functions, and $\mathcal L$ is a lift of the data (Section~\ref{sec:constraints}). For steady problems the factor $C_r$ is omitted. Each factor network $\mathbf f:\R\to\R^{R}$ has the form
\begin{equation}
\mathbf h^{(0)}=W_0^{\top}\boldsymbol\psi(z),\qquad
h^{(\ell+1)}_o=\sum_{i=1}^{H}\sum_{k=0}^{K-1} c^{(\ell)}_{iko}\,T_k\!\big(\tanh h^{(\ell)}_i\big),\quad \ell=0,1,
\label{eq:factor}
\end{equation}
where $T_k$ are the Chebyshev polynomials of Eq.~\eqref{eq:cheb}. The layers after the first are therefore Chebyshev-KAN layers \cite{ChebyKAN,Mos2} with learnable univariate edge functions. The first layer is a learnable univariate function in a basis $\boldsymbol\psi$ \emph{adapted to the topology of the coordinate}:
\begin{equation}
\boldsymbol\psi(z)=\begin{cases}
\big(1,\cos 2\pi z,\sin 2\pi z,\dots,\cos 2\pi K_f z,\sin 2\pi K_f z\big) & z \text{ periodic},\\
\big(T_0(2z-1),\dots,T_{K-1}(2z-1)\big) & z \text{ bounded or time}.
\end{cases}
\label{eq:psi}
\end{equation}
The Fourier layer makes $\hat u_\theta$ and all of its $\eta$-derivatives exactly $1$-periodic, so the representation has no seam along $\eta=0\equiv1$. Section~\ref{sec:theory} shows that this is needed for the error bounds.

\begin{remark}[Expressivity]
Finite sums of products of continuous univariate functions form a subalgebra of $C(\bar\Q\times[0,T])$ that contains the constants and separates points. By the Stone--Weierstrass theorem it is therefore dense, so a separable ansatz with sufficiently large $R$ can approximate any continuous pulled-back solution uniformly. In the periodic case, the same argument applies on $[0,1]\times\mathbb T\times[0,T]$ with trigonometric functions of $\eta$. The rank needed depends on how well the solution aligns with the computational coordinates; this is the main approximation-theoretic limitation of the approach (Section~\ref{sec:limits}).
\end{remark}

\paragraph{Initialisation and size.} The initial coefficients of mode $k$ in the first layer are scaled by $(1+k)^{-2}$, and those of degree $k$ in the Chebyshev layers by $(1+k)^{-1}$. This keeps the initial field and its second derivatives of order one, which matters because the operator contains second derivatives of high-frequency modes. Unless stated otherwise we use $R=32$, $H=16$, $K=8$, $K_f=8$ and two KAN layers. This gives $12{,}688$ parameters for steady problems, $18{,}816$ for transient problems on non-periodic maps, and $18{,}960$ for transient problems on periodic maps.

\subsection{Exact constraints and a single-term loss}
\label{sec:constraints}
Because every physical boundary is a coordinate line of $\Q$, boundary and initial data can be built into \eqref{eq:ansatz}, so that they are satisfied exactly rather than penalised.
\begin{itemize}
\item \emph{Homogeneous Dirichlet data} (Poisson, heat). Here $\mathcal L=0$ (Poisson) or $\mathcal L=\hat u_0(\xi,\eta)$ (heat). The envelopes are $e_\xi=\xi(1-\xi)$, $e_\eta=\eta(1-\eta)$ if $\eta=0,1$ are walls ($e_\eta=1$ if $\eta$ is periodic), and $e_t=t/T$. The initial data must vanish on the walls (compatibility).
\item \emph{Inflow data} for advection, prescribed on the inflow faces $\xi=0$, $\eta=0$ and at $t=0$. We use the Boolean sum of the face projectors \cite{Gordon1973},
\begin{equation}
\mathcal L=(P_\xi\oplus P_\eta\oplus P_t)\,d
= \big(P_\xi+P_\eta+P_t-P_\xi P_\eta-P_\xi P_t-P_\eta P_t+P_\xi P_\eta P_t\big)\,d,
\label{eq:coons}
\end{equation}
with $(P_\xi d)(\xi,\eta,t)=d(0,\eta,t)$, $(P_\eta d)=d(\xi,0,t)$ and $(P_t d)=d(\xi,\eta,0)$. Here $d$ denotes the prescribed boundary and initial data. The envelopes are $e_\xi=\xi$, $e_\eta=\eta$ and $e_t=t/T$.
\end{itemize}

\begin{lemma}[The Boolean-sum lift interpolates the data]\label{lem:coons}
Let $d$ be continuous and consistent at the edges where faces meet. Then $\mathcal L$ of Eq.~\eqref{eq:coons} satisfies $\mathcal L(0,\eta,t)=d(0,\eta,t)$, $\mathcal L(\xi,0,t)=d(\xi,0,t)$ and $\mathcal L(\xi,\eta,0)=d(\xi,\eta,0)$. Consequently $\hat u_\theta$ in Eq.~\eqref{eq:ansatz} satisfies the inflow and initial data exactly for all parameters $\theta$.
\end{lemma}
\begin{proof}
Each projector is idempotent and $P_\xi d$ evaluated at $\xi=0$ equals $d(0,\cdot,\cdot)$. Writing $\mathcal L=d-(I-P_\xi)(I-P_\eta)(I-P_t)d$, the product vanishes at $\xi=0$, because $(I-P_\xi)h$ vanishes at $\xi=0$ for any $h$ and the projectors commute. The same holds at $\eta=0$ and $t=0$. The envelopes vanish on the same faces, so the separable part does not change the traces.
\end{proof}

With all data built in, the loss is a single term, and there are no penalty weights to tune:
\begin{equation}
\mathcal J(\theta)=\frac{1}{N_\xi N_\eta N_t}\sum_{i,j,k}\sqrt{g}(\xi_i,\eta_j)\;\mathcal R_\theta(\xi_i,\eta_j,t_k)^2 ,
\label{eq:loss}
\end{equation}

\begin{figure}[H]
\centering
\includegraphics[width=0.85\textwidth]{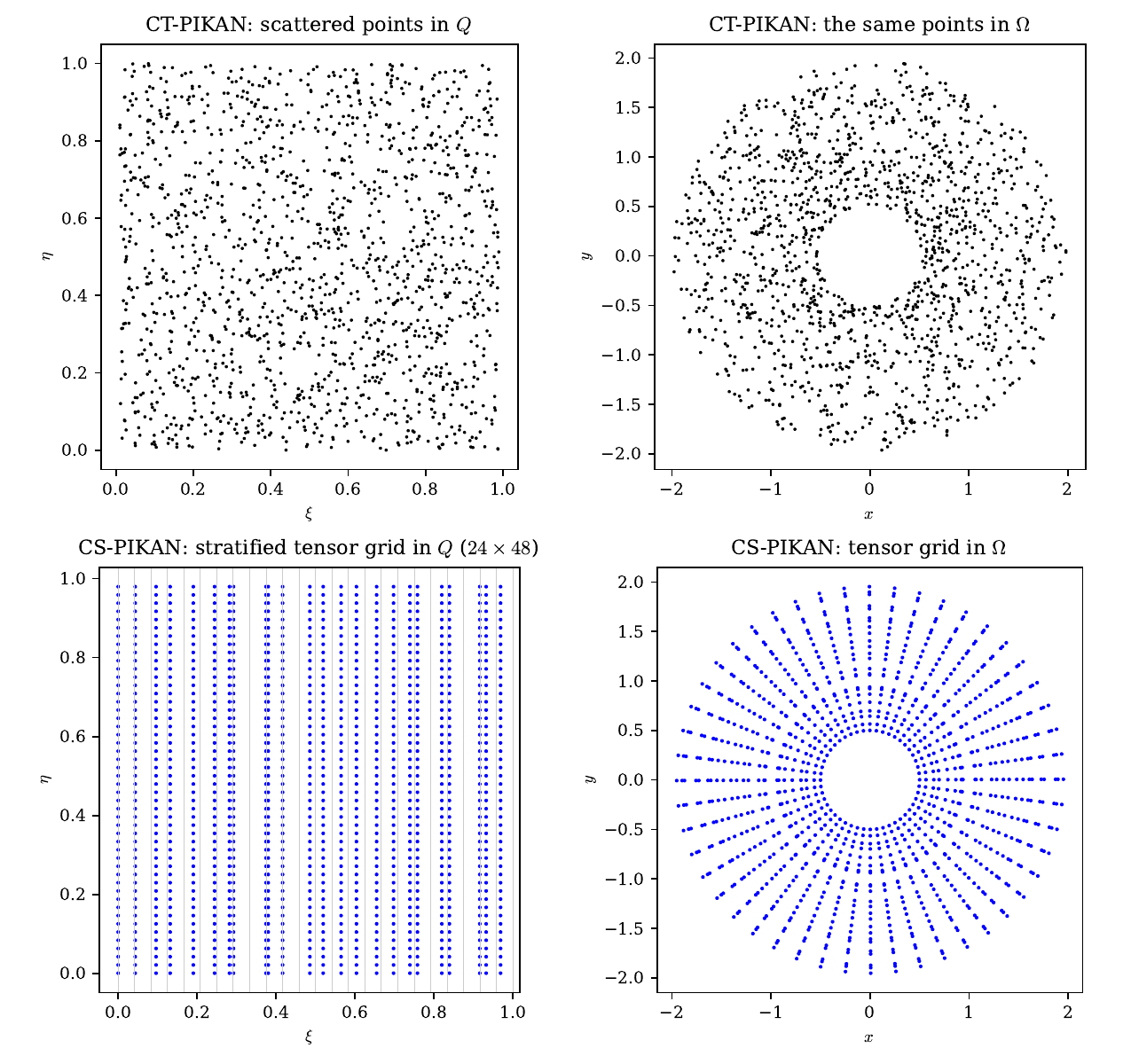}
\caption{Collocation in CT-PIKAN (top: scattered random points, shown in the computational and the physical domain) and in CS-PIKAN (bottom: one draw of the stratified tensor grid; grey lines mark the strata in $\xi$, and the periodic $\eta$-grid is shifted as a whole). The annulus map is used.}
\label{fig:sampling}
\end{figure}
\noindent where $\mathcal R_\theta$ is the PDE residual. The nodes are drawn afresh at every iteration by \emph{stratified tensor sampling}: $\xi_i=(i+U_i)/N_\xi$ and $t_k=T(k+V_k)/N_t$ with independent uniform $U_i,V_k$. A periodic $\eta$ uses a randomly shifted uniform grid $\eta_j=(j+S)/N_\eta$, and a bounded $\eta$ is stratified like $\xi$. The weight $\sqrt g$ makes \eqref{eq:loss} a quadrature of the residual over the \emph{physical} domain (Proposition~\ref{prop:unbiased}). Without it, the loss measures the residual in the computational rather than the physical domain and is no longer tied to the certificate. Fig.~\ref{fig:sampling} contrasts the scattered points of CT-PIKAN with the stratified tensor grid of CS-PIKAN.
\subsection{Separable assembly of curvilinear operators}
\begin{proposition}[Exact separable assembly]\label{prop:sep}
Let $A,A',A''\in\R^{N_\xi\times R}$, $B,B',B''\in\R^{N_\eta\times R}$ and $C,C'\in\R^{N_t\times R}$ collect the factors of \eqref{eq:ansatz} and their derivatives at the grid nodes, and write $\llbracket X,Y,Z\rrbracket_{ijk}=\sum_r X_{ir}Y_{jr}Z_{kr}$. On the tensor grid,
\begin{equation}
\Delta_x \hat u_\theta = g^{11}\llbracket A'',B,C\rrbracket+2g^{12}\llbracket A',B',C\rrbracket+g^{22}\llbracket A,B'',C\rrbracket+b^1\llbracket A',B,C\rrbracket+b^2\llbracket A,B',C\rrbracket+\Delta_x\mathcal L,
\label{eq:sepLB}
\end{equation}
where the metric fields $g^{ij},b^j$ are evaluated on the $N_\xi\times N_\eta$ spatial grid and broadcast in $t$. Analogous identities hold for $\partial_t\hat u_\theta=\llbracket A,B,C'\rrbracket+\partial_t\mathcal L$ and for $\mathbf w\cdot\nabla_\xi\hat u_\theta$. Evaluating the networks requires $N_\xi+N_\eta+N_t$ forward-mode passes, against $N_\xi N_\eta N_t$ second-order automatic-differentiation passes for a point-wise network on the same nodes. The metric fields depend only on the nodes, not on $\theta$.
\end{proposition}
\begin{proof}
Identity \eqref{eq:sepLB} follows by inserting \eqref{eq:ansatz} into \eqref{eq:LB}. Mixed derivatives of a product of univariate functions factorise ($\partial_{\xi\eta}(A_rB_rC_r)=A_r'B_r'C_r$), and all metric terms are pointwise multipliers. The cost statement follows because each factor depends on one variable only.
\end{proof}

We verified \eqref{eq:sepLB} against point-wise automatic differentiation of the assembled network for all problem classes; the two agree to $6$--$7$ significant digits in single precision. Table~\ref{tab:complexity} compares the leading costs of the two formulations, and Algorithm~\ref{alg:cs} summarises one training iteration. Fig.~\ref{fig:schematic} shows the architecture.

\begin{table}[H]
\centering
\caption{Leading cost of one training iteration for a transient problem with $N=N_\xi N_\eta N_t$ residual nodes. $P$ is the cost of one network evaluation, $c_2$ the overhead factor of nested second-order automatic differentiation, and $P_1$ the cost of a one-dimensional factor network.}
\label{tab:complexity}
\footnotesize
\begin{tabular}{lcc}
\toprule
Operation & CT-PIKAN (point-wise) & CS-PIKAN (separable)\\
\midrule
Network evaluations with derivatives & $c_2\,N\,P$ & $c_2\,(N_\xi+N_\eta+N_t)\,P_1$\\
Assembly of the residual & included above & $O(R\,N)$ tensor contractions\\
Metric fields & $O(N)$ (per point) & $O(N_\xi N_\eta)$ (spatial grid only)\\
Boundary and initial terms & $O((N_b+N_0)P)$ & none (exact constraints)\\
\bottomrule
\end{tabular}
\end{table}

\begin{figure}[H]
\centering
\resizebox{\textwidth}{!}{%
\begin{tikzpicture}[font=\small,>=Stealth,
  box/.style={draw,rounded corners=2pt,align=center,minimum height=8mm,inner sep=3pt},
  kan/.style={box,fill=blue!8,anchor=west,minimum width=7.4cm},
  geo/.style={box,fill=orange!12},
  resbox/.style={box,fill=green!10,minimum width=5.6cm}]
\draw[thick,fill=gray!10] plot[smooth cycle,tension=.8] coordinates {(0,0.2) (0.6,1.0) (1.6,1.1) (2.0,0.3) (1.5,-0.6) (0.5,-0.6)};
\draw[thick,fill=white] (1.0,0.25) circle (0.28);
\node at (1.0,-1.05) {$\Omega$ (physical)};
\draw[thick,fill=gray!10] (3.3,-0.6) rectangle (4.8,0.9);
\foreach \s in {0.3,0.6,0.9,1.2} {\draw[gray] (3.3+\s,-0.6) -- (3.3+\s,0.9); \draw[gray] (3.3,-0.6+\s) -- (4.8,-0.6+\s);}
\node at (4.05,-1.05) {$\Q=(0,1)^2$};
\draw[->,thick] (3.2,0.15) -- node[above]{$\Phi$} (2.1,0.15);
\node[kan] (fa) at (5.8,1.5) {$\xi\ \to$ Chebyshev layer $\to$ Cheb-KAN $\to\ \mathbf a(\xi)\in\R^R$};
\node[kan] (fb) at (5.8,0.15) {$\eta\ \to$ \textbf{Fourier} layer (periodic) $\to$ Cheb-KAN $\to\ \mathbf b(\eta)\in\R^R$};
\node[kan] (fc) at (5.8,-1.2) {$t\ \to$ Chebyshev layer $\to$ Cheb-KAN $\to\ \mathbf c(t)\in\R^R$};
\draw[->] (4.8,0.7) -- (fa.west); \draw[->] (4.8,0.15) -- (fb.west); \draw[->] (4.8,-0.4) -- (fc.west);
\node[resbox] (sum) at (19.9,0.15) {$\hat u_\theta=\mathcal L+\sum_{r=1}^R A_r(\xi)\,B_r(\eta)\,C_r(t)$\\[2pt]
  {\footnotesize $A_r=e_\xi a_r,\ \ B_r=e_\eta b_r,\ \ C_r=e_t c_r$}\\[1pt]{\footnotesize exact BC/IC: envelope + Boolean-sum lift}};
\draw[->] (fa.east) -- (sum.west); \draw[->] (fb.east) -- (sum.west); \draw[->] (fc.east) -- (sum.west);
\node[geo,minimum width=5.2cm] (metric) at (4.05,-2.9) {automatic differentiation of $\Phi$:\\ $\sqrt g,\ g^{ij},\ b^{j},\ \mathbf w=J^{-1}\mathbf v$ on the $\xi\times\eta$ grid};
\node[resbox] (op) at (19.9,-2.9) {separable operator (Proposition~\ref{prop:sep})\\ residual $\mathcal R_\theta$ on the $N_\xi\times N_\eta\times N_t$ grid};
\draw[->] (sum) -- (op); \draw[->] (metric.east) -- (op.west);
\node[resbox] (loss) at (19.9,-4.6) {$\mathcal J=\mathrm{mean}\big(\sqrt g\,\mathcal R_\theta^2\big)\ \approx\ \|\mathcal R_\theta\|^2_{L_2(\Omega)}$\\{\footnotesize squared certificate (Theorem~\ref{thm:apost})}};
\draw[->] (op) -- (loss);
\draw[->] (4.05,-0.6) -- (metric.north);
\end{tikzpicture}}

\caption{CS-PIKAN. The curvilinear domain is pulled back to the unit square. Each computational coordinate is processed by its own one-dimensional KAN with a topology-adapted first layer, and the factors are combined in a rank-$R$ separable sum together with the exact-constraint envelope and lift. Metric fields from automatic differentiation of the map enter only as pointwise multipliers in the separable operator (Proposition~\ref{prop:sep}). The weighted loss is the squared certificate of Theorem~\ref{thm:apost}.}
\label{fig:schematic}
\end{figure}

\begin{algorithm}[H]
\caption{One CS-PIKAN training iteration}\label{alg:cs}
\begin{algorithmic}[1]
\Require map $\Phi$, data lift $\mathcal L$, envelopes $e_\xi,e_\eta,e_t$, factor networks $\mathbf a,\mathbf b,\mathbf c$
\State Draw stratified nodes $\{\xi_i\}$, $\{\eta_j\}$, $\{t_k\}$
\State Metric fields on $\{\xi_i\}\times\{\eta_j\}$ by Algorithm~\ref{al:1}: $\sqrt g$, $g^{ij}$, $b^j$ (and $\mathbf w=\mathbf J^{-1}\mathbf v$)
\State Lift $\mathcal L$ and its derivatives on the grid (explicit function, automatic differentiation)
\State Factors and derivatives on the 1-D nodes by forward-mode AD: $A,A',A''$; $B,B',B''$; $C,C'$
\State Assemble the residual $\mathcal R_\theta$ on the tensor grid with Proposition~\ref{prop:sep}
\State $\mathcal J\gets$ mean of $\sqrt g\,\mathcal R_\theta^2$; take an Adam step on $\theta$
\end{algorithmic}
\end{algorithm}

\section{Analysis: geometry-explicit certified error bounds}
\label{sec:theory}

Throughout, $\norm{\cdot}$ denotes the $L_2(\Omega)$ norm, and $\Q$ carries Dirichlet conditions on $\xi\in\{0,1\}$ together with either periodicity in $\eta$ (\emph{case P}) or Dirichlet conditions on $\eta\in\{0,1\}$ (\emph{case D}). We set $\lambda_\Q=\pi^2$ in case P and $\lambda_\Q=2\pi^2$ in case D; these are the first Dirichlet(-periodic) eigenvalues of $-\Delta_\xi$ on $\Q$. As before, $\sigma_1\ge\sigma_2>0$ are the singular values of $\mathbf J$ and $\kappa=\sigma_1/\sigma_2$ is its condition number.

\begin{lemma}[Pull-back identities]\label{lem:pull}
For $v\in H^1(\Omega)$ with $\hat v=v\circ\Phi$,
\[
\norm{v}^2=\int_\Q \hat v^2\sqrt g\,d\bxi,\qquad \norm{\nabla v}^2=\int_\Q \nabla_\xi\hat v^{\top} M\,\nabla_\xi\hat v\,d\bxi .
\]
Moreover, $\det M=1$, the eigenvalues of $M$ are $\kappa^{-1}$ and $\kappa$, and $M^{-1}=G/\sqrt g$.
\end{lemma}
\begin{proof}
The identities follow from $dx=\sqrt g\,d\bxi$ and
\[|\nabla_x v|^2=\nabla_\xi\hat v^\top \mathbf J^{-1}\mathbf J^{-\top}\nabla_\xi\hat v=\nabla_\xi\hat v^\top G^{-1}\nabla_\xi\hat v.\] The eigenvalues of $G^{-1}$ are $\sigma_1^{-2}$ and $\sigma_2^{-2}$, and $\sqrt g=\sigma_1\sigma_2$, so the eigenvalues of $M$ are $\sigma_2/\sigma_1$ and $\sigma_1/\sigma_2$, and their product is $1$.
\end{proof}

\begin{theorem}[Geometry-explicit Poincar\'e inequality]\label{thm:poincare}
Let $v\in H^1(\Omega)$ be such that $\hat v$ vanishes on $\xi\in\{0,1\}$ (and on $\eta\in\{0,1\}$ in case D). Then $\norm{v}\le C_\Omega\norm{\nabla v}$, with
\begin{equation}
C_\Omega^2=\Big(\max_\Q\sqrt g\Big)\min\Big\{\frac{\max_\Q\kappa}{\lambda_\Q},\ \frac{1}{\pi^2}\max_\Q\frac{g_{11}}{\sqrt g},\ \Big[\frac{1}{\pi^2}\max_\Q\frac{g_{22}}{\sqrt g}\Big]_{\text{case D}}\Big\}.
\label{eq:CP}
\end{equation}
\end{theorem}
\begin{proof}
By Lemma~\ref{lem:pull}, $\norm{v}^2\le\max\sqrt g\int_\Q\hat v^2$. (i) By the Poincar\'e inequality on $\Q$, $\int_\Q\hat v^2\le\lambda_\Q^{-1}\int_\Q|\nabla_\xi\hat v|^2$. Since $\lambda_{\max}(M^{-1})=\kappa$,
\[|\nabla_\xi\hat v|^2\le\kappa\,\nabla_\xi\hat v^\top M\nabla_\xi\hat v,\]
and therefore $\int_\Q|\nabla_\xi\hat v|^2\le\max\kappa\,\norm{\nabla v}^2$. (ii) For fixed $\eta$, the one-dimensional inequality $\int_0^1\hat v^2d\xi\le\pi^{-2}\int_0^1\hat v_\xi^2d\xi$ holds. By the Cauchy--Schwarz inequality in the $M$-inner product,
\[\hat v_\xi^2=(\mathbf e_1\!\cdot\!\nabla_\xi\hat v)^2\le(\mathbf e_1^\top M^{-1}\mathbf e_1)(\nabla_\xi\hat v^\top M\nabla_\xi\hat v)=\frac{g_{11}}{\sqrt g}\,\nabla_\xi\hat v^\top M\nabla_\xi\hat v.\] Integrating over $\eta$ gives the second term. The third term follows in the same way with $\mathbf e_2$. Each bound is valid, so their minimum is too.
\end{proof}

All quantities in \eqref{eq:CP} are obtained from the autograd Jacobian of the map. Bound (ii) is direction-aware and much sharper than (i) for strongly anisotropic maps: for the annulus, (i) gives $C_\Omega^2=16.0$, whereas (ii) gives $0.912$ (Table~\ref{tab:geom}). For the annulus the constant is available in closed form: with $\sqrt g=3\pi r$ and $g_{11}=9/4$ (Appendix~\ref{app:metrics}), (ii) gives $C_\Omega^2=6\pi\cdot\frac{9/4}{3\pi\cdot 0.5}\cdot\pi^{-2}=9/\pi^2\approx0.912$.

\begin{theorem}[A-posteriori certificates]\label{thm:apost}
Let $\hat u_\theta$ be a CS-PIKAN approximation, i.e.\ smooth, exactly periodic in $\eta$ in case P, and satisfying the boundary and initial data exactly. Let $e=u-u_\theta$, where $u_\theta=\hat u_\theta\circ\Phi^{-1}$.
\begin{enumerate}
\item[(a)] \emph{Poisson} ($-\Delta u=f$ in $\Omega$, $u=0$ on $\partial\Omega$), with residual $\mathcal R=\Delta u_\theta+f$:
\[\norm{\nabla e}\le C_\Omega\norm{\mathcal R},\qquad \norm{e}\le C_\Omega^2\norm{\mathcal R}.\]
\item[(b)] \emph{Heat} ($u_t=\alpha\Delta u$, homogeneous Dirichlet data), with residual $\mathcal R=\partial_tu_\theta-\alpha\Delta u_\theta$:
\[\norm{e(t)}\le\int_0^t e^{-\alpha(t-s)/C_\Omega^2}\norm{\mathcal R(s)}\,ds,\qquad
\norm{e(T)}^2+\alpha\!\int_0^T\!\norm{\nabla e}^2dt\le\frac{C_\Omega^2}{\alpha}\int_0^T\!\norm{\mathcal R}^2dt.\]
\item[(c)] \emph{Advection} ($u_t+\mathbf v\cdot\nabla u=0$, data on the inflow boundary $\Gamma_-=\{\mathbf v\cdot\mathbf n<0\}$, and $\mathbf v\cdot\mathbf n\ge0$ on $\partial\Omega\setminus\Gamma_-$), with residual $\mathcal R=\partial_tu_\theta+\mathbf v\cdot\nabla u_\theta$ and $\beta=\max(0,\sup_\Omega\nabla\!\cdot\mathbf v)$:
\[\norm{e(t)}\le\int_0^t e^{\beta(t-s)/2}\norm{\mathcal R(s)}\,ds .\]
\end{enumerate}
\end{theorem}
\begin{proof}
Exact constraints give $e=0$ on the Dirichlet (respectively inflow) boundary and $e(\cdot,0)=0$. Exact periodicity in case P means that $\hat e$ has no jump across $\eta=0\equiv1$, so $e\in H^1(\Omega)$ and Green's formula applies without interface terms. (a) Since $-\Delta e=\mathcal R$, we have $\norm{\nabla e}^2=(\mathcal R,e)\le\norm{\mathcal R}\norm{e}\le C_\Omega\norm{\mathcal R}\norm{\nabla e}$ by Theorem~\ref{thm:poincare}, and $\norm{e}\le C_\Omega\norm{\nabla e}$ completes the argument. (b) Since $e_t-\alpha\Delta e=-\mathcal R$, testing with $e$ gives $\tfrac12\tfrac{d}{dt}\norm{e}^2+\alpha\norm{\nabla e}^2=-(\mathcal R,e)$. With $\alpha\norm{\nabla e}^2\ge\alpha C_\Omega^{-2}\norm{e}^2$ we get $\tfrac{d}{dt}\norm{e}\le-\alpha C_\Omega^{-2}\norm{e}+\norm{\mathcal R}$, and Gr\"onwall's lemma with $e(0)=0$ yields the first bound. For the second, $|(\mathcal R,e)|\le C_\Omega\norm{\mathcal R}\norm{\nabla e}\le\frac{C_\Omega^2}{2\alpha}\norm{\mathcal R}^2+\frac{\alpha}{2}\norm{\nabla e}^2$; integrate in time. (c) Since $e_t+\mathbf v\cdot\nabla e=-\mathcal R$, testing with $e$ and integrating by parts gives $\tfrac12\tfrac{d}{dt}\norm{e}^2+\tfrac12\int_{\partial\Omega}(\mathbf v\cdot\mathbf n)e^2-\tfrac12\int_\Omega(\nabla\!\cdot\mathbf v)e^2=-(\mathcal R,e)$. The boundary integral is non-negative because $e=0$ on $\Gamma_-$. Hence $\tfrac{d}{dt}\norm{e}\le\tfrac\beta2\norm{e}+\norm{\mathcal R}$, and Gr\"onwall's lemma concludes.
\end{proof}

\begin{proposition}[The loss estimates the certificate]\label{prop:unbiased}
For stratified tensor sampling, $\mathbb E[\mathcal J(\theta)]=\frac1T\int_0^T\norm{\mathcal R_\theta(t)}^2dt$, and $\mathbb E[\mathcal J(\theta)]=\norm{\mathcal R_\theta}^2$ for steady problems.
\end{proposition}
\begin{proof}
Each node is uniformly distributed in its grid cell, and every cell has volume $T/(N_\xi N_\eta N_t)$. This holds for the shifted periodic grid as well, because the marginal of $\eta_j$ is uniform on its cell. Therefore
\[\mathbb E\big[\sqrt g\,\mathcal R^2(\xi_i,\eta_j,t_k)\big]=\frac{N_\xi N_\eta N_t}{T}\int_{\text{cell}}\sqrt g\,\mathcal R^2.\] Summing over cells, dividing by $N_\xi N_\eta N_t$, and applying Lemma~\ref{lem:pull} gives the result.
\end{proof}

\begin{corollary}[Energy certificate for the heat equation]\label{cor:energy}
Under the hypotheses of Theorem~\ref{thm:apost}(b), for the expected training loss $\bar{\mathcal J}=\mathbb E[\mathcal J]$ of Proposition~\ref{prop:unbiased},
\[\norm{e(T)}^2+\alpha\int_0^T\norm{\nabla e}^2dt\le\frac{C_\Omega^2\,T}{\alpha}\,\bar{\mathcal J}(\theta).\]
Thus driving the expected training loss to zero drives a guaranteed space--time energy error to zero.
\end{corollary}
\begin{proof}
Combine the second inequality of Theorem~\ref{thm:apost}(b) with Proposition~\ref{prop:unbiased}.
\end{proof}

\begin{remark}
Theorem~\ref{thm:apost} and Proposition~\ref{prop:unbiased} together show that CS-PIKAN training minimises an unbiased estimate of the square of a guaranteed error bound. Two of its design choices are exactly the hypotheses needed. Penalty-based imposition of the boundary data, as in CT-PIKAN, leaves $e\ne0$ on $\partial\Omega$. A non-periodic representation of a closed coordinate leaves a jump of $\hat e$ across the seam, so that $e\notin H^1(\Omega)$. In both cases the certificate is lost, and the uncontrolled boundary or interface terms would have to be estimated separately.
\end{remark}

\begin{remark}[Sharpness]
The certificates can be pessimistic for three reasons: the maxima in \eqref{eq:CP} are taken over the whole domain; the Poisson bound uses the Poincar\'e constant twice; and the error is controlled by the $L_2$ norm of the residual, whereas the sharp quantity is its $H^{-1}$ norm, which is much smaller for oscillatory residuals. Table~\ref{tab:geom} quantifies the first effect by comparing $C_\Omega^2$ with the sharp value $1/\lambda_1(\Omega)$, and Table~\ref{tab:bounds} reports the resulting effectivity indices.
\end{remark}

In practice, $\norm{\mathcal R(t)}$ is computed with a $48$-point Gauss--Legendre rule in each bounded coordinate ($96$ uniform points in a periodic one) and $201$ time levels, and the maxima in \eqref{eq:CP} are taken on a $401\times401$ grid. The certificates are therefore exact up to quadrature error, which is negligible at these resolutions.

\section{Numerical experiments}
\label{sec:results}

In this section we evaluate CS-PIKAN and the CT-PIKAN baseline on nine benchmarks spanning hyperbolic, elliptic, and parabolic PDEs on Cartesian and curvilinear domains: the two-dimensional linear advection equation on Cartesian, polar and non-orthogonal wavy domains; the Poisson equation on an annulus and on a domain whose boundaries are known only through sample points; and the transient heat equation on a wavy channel and a star-shaped domain. The examples increase progressively in geometric and mathematical complexity. In all cases the metrics are obtained by automatic differentiation of the map.

\subsection{Protocol}
\label{sec:protocol}
\paragraph{Baseline.} CT-PIKAN (Section~\ref{sec:ctpikan}) is run with the configurations of its public code: depth $4$, $4096$ interior and $1024$ boundary points, $4000$ epochs, step decay (factor $0.5$ every $1000$ epochs) and gradient clipping for Poisson; depth $3$, $3000$ epochs and a constant learning rate for the wavy channel; depth $4$, $2500$ interior points, $3000$ epochs, step decay and clipping for the star; and the Poisson settings with $1024$ initial points for advection. Both methods solve exactly the same problems and are implemented in the same JAX code base.
\paragraph{CS-PIKAN.} We use $R=32$ and Adam with peak learning rate $5\times10^{-3}$, decayed exponentially to $5\times10^{-5}$ over $6000$ epochs. Tensor grids have $48\times96$ nodes for the steady annular problems and $32\times\{32,64\}\times24$ nodes for transient problems. Table~\ref{tab:config} lists all settings.

\begin{table}[H]
\centering
\caption{Configurations of all test cases. CS-PIKAN: rank $32$, two Chebyshev-KAN layers of width $16$, $K=8$, $K_f=8$, $6000$ epochs, Adam $5\times10^{-3}\to5\times10^{-5}$, no loss weights. CT-PIKAN: RBF-KAN input layer ($16$ bases) and $\tanh$ trunk of width $64$, Adam $10^{-3}$, loss weights $\lambda_b=\lambda_0=10$.}
\label{tab:config}
\footnotesize
\begin{adjustbox}{max width=\linewidth}
\begin{tabular}{lcccccc}
\toprule
& \multicolumn{3}{c}{CS-PIKAN} & \multicolumn{3}{c}{CT-PIKAN}\\
\cmidrule(lr){2-4}\cmidrule(lr){5-7}
Case & grid $N_\xi\times N_\eta(\times N_t)$ & $\eta$-basis & params & depth / epochs & $N_r/N_b/N_0$ & params\\
\midrule
Advection, Cartesian & $32\times32\times24$ & Chebyshev & 18,816 & 4 / 4000 & 4096/1024/1024 & 15,777\\
Advection, polar ($\beta=2,8$) & $32\times32\times24$ & Chebyshev & 18,816 & 4 / 4000 & 4096/1024/1024 & 15,777\\
Advection, wavy & $32\times32\times24$ & Chebyshev & 18,816 & 4 / 4000 & 4096/1024/1024 & 15,777\\
Poisson, annulus & $48\times96$ & Fourier & 12,688 & 4 / 4000 & 4096/1024/-- & 14,721\\
Poisson, sampled boundary & $48\times96$ & Fourier & 12,688 & 4 / 4000 & 4096/1024/-- & 14,721\\
Heat, wavy channel & $32\times32\times24$ & Chebyshev & 18,816 & 3 / 3000 & 4096/1024/1024 & 11,617\\
Heat, star (both I.C.) & $32\times64\times24$ & Fourier & 18,960 & 4 / 3000 & 2500/1024/1024 & 15,777\\
\bottomrule
\end{tabular}
\end{adjustbox}
\end{table}

\paragraph{References and errors.} The advection problems have exact solutions. For the Poisson and heat problems we solve the transformed equations with a second-order conservative nine-point finite-difference scheme on the same body-fitted grid, with Crank--Nicolson in time (Appendix~\ref{app:fd}). Manufactured solutions verify the scheme at second order on all four curved maps (Table~\ref{tab:fdver}). We combine the $256$-node and $128$-node solutions by Richardson extrapolation; the extrapolated reference differs from the fine-grid solution by at most $3.91\times10^{-4}$. Errors are reported as relative $\sqrt g$-weighted (that is, physical) $L_2$ and relative $L_\infty$ norms, on $65\times128$ (periodic) or $65\times65$ nodes, pooled over the snapshot times, as mean $\pm$ standard deviation over three seeds.

\begin{table}[H]\centering\caption{Verification of the finite-difference reference scheme with manufactured solutions: relative $L_2$ error on $N\times N$ ($N\times 2N$ for periodic maps) grids and observed orders.}\label{tab:fdver}\footnotesize\begin{tabular}{lccccc}\toprule
Map & $N=16$ & $N=32$ & $N=64$ & order $16\to32$ & order $32\to64$\\\midrule
Annulus & $3.36\times10^{-3}$ & $8.37\times10^{-4}$ & $2.09\times10^{-4}$ & 2.00 & 2.00\\
Wavy channel & $9.43\times10^{-3}$ & $2.51\times10^{-3}$ & $6.39\times10^{-4}$ & 1.91 & 1.97\\
Star & $3.19\times10^{-2}$ & $9.57\times10^{-3}$ & $2.60\times10^{-3}$ & 1.74 & 1.88\\
Sampled boundary (TFI) & $4.87\times10^{-3}$ & $1.26\times10^{-3}$ & $3.15\times10^{-4}$ & 1.95 & 2.00\\
\bottomrule\end{tabular}\end{table}

\paragraph{Timing.} All runs use a single core of an Intel\textsuperscript{\textregistered} Xeon\textsuperscript{\textregistered} Processor @ 2.30GHz CPU, JAX 0.12.0.1, and single precision. Training time excludes compilation and was measured separately from the accuracy runs, on an otherwise idle machine.

\subsection{Geometries and metric fields}
\label{sec:geoms}
\begin{figure}[H]
\centering
\includegraphics[width=\textwidth]{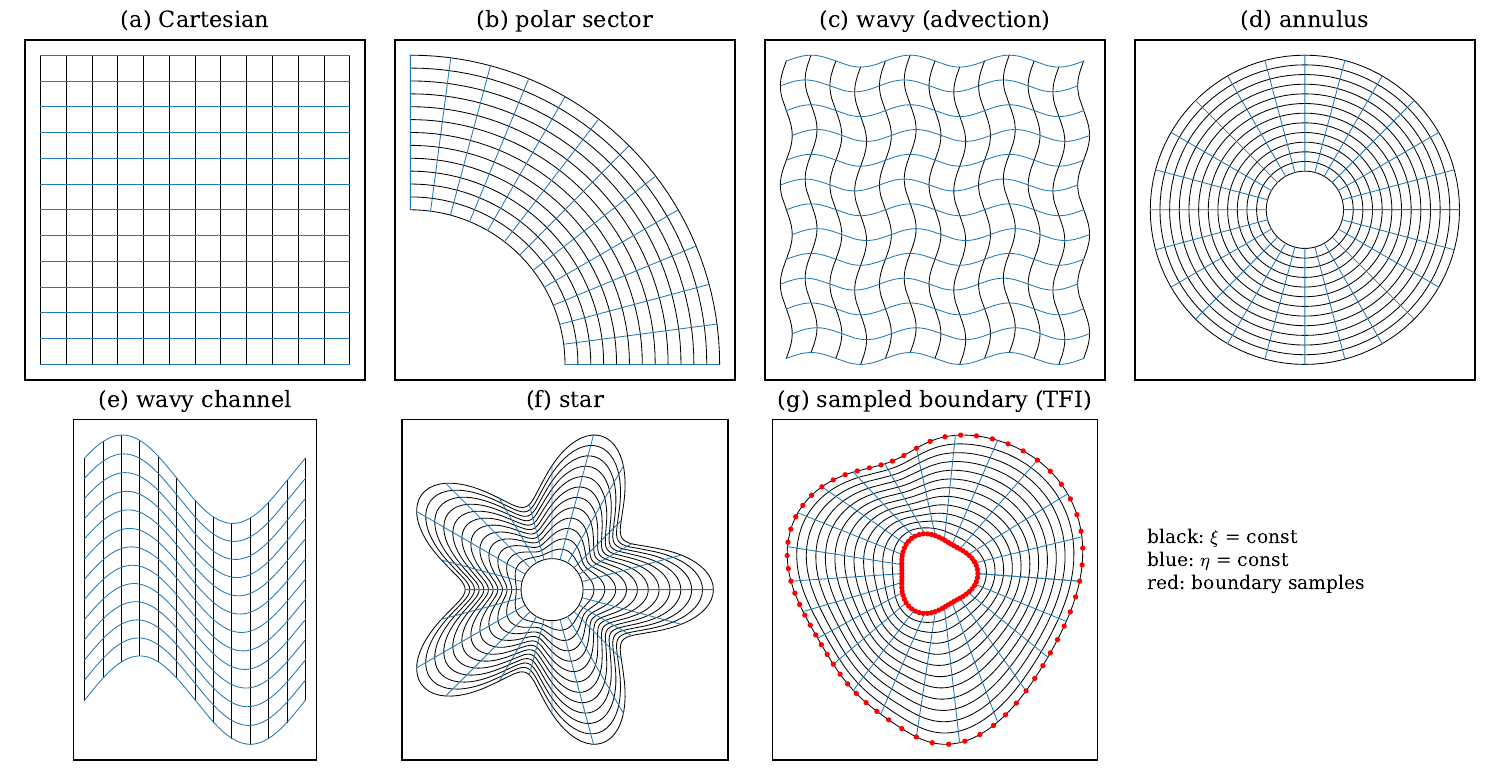}
\caption{Computational maps $\Phi(\Q)$ shown by their coordinate lines. In (g), the domain boundary is known only through the $64$ marked sample points per curve. The map is built by trigonometric interpolation of the samples and transfinite blending, and differentiated by autograd.}
\label{fig:geoms}
\end{figure}

\begin{figure}[H]
\centering
\includegraphics[width=0.95\textwidth]{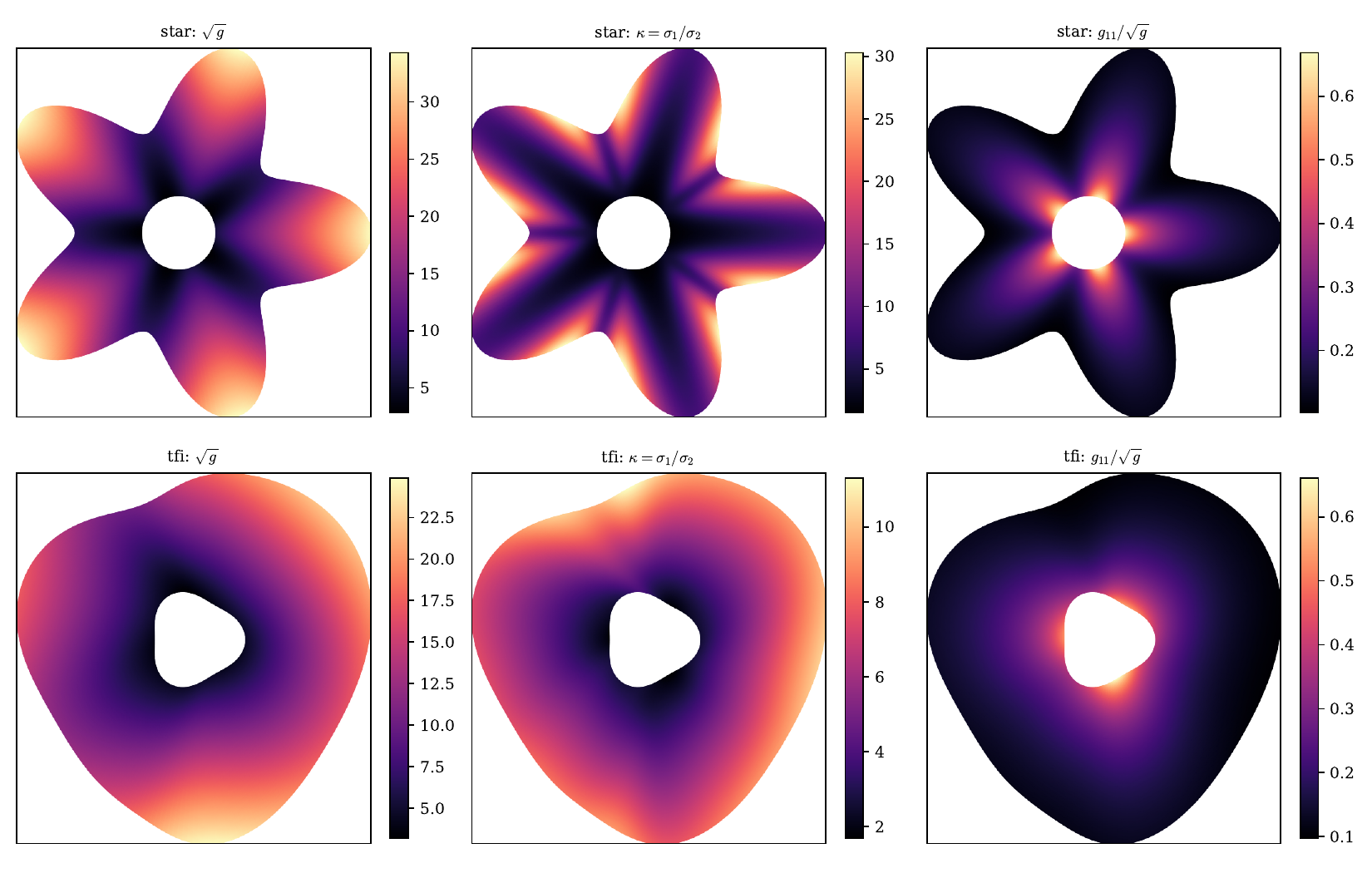}
\caption{Metric fields of the star map (top) and the sampled-boundary TFI map (bottom): area element $\sqrt g$, anisotropy $\kappa$, and the directional factor $g_{11}/\sqrt g$ that enters bound (ii) of Theorem~\ref{thm:poincare}. All fields are obtained by automatic differentiation of the maps.}
\label{fig:metric}
\end{figure}

Fig.~\ref{fig:geoms} shows the seven maps by their coordinate lines. Table~\ref{tab:geom} lists their metric extremes and the Poincar\'e constants of Theorem~\ref{thm:poincare}, and Fig.~\ref{fig:metric} shows the fields that enter the constants for the two most distorted maps.

\begin{table}[H]\centering\caption{Metric extremes of the maps (autograd) and Poincar\'e constants of Theorem~\ref{thm:poincare}: isotropic bound (i), the final constant $C_\Omega^2$ of \eqref{eq:CP}, and the sharp value $1/\lambda_1(\Omega)$ from the finite-difference eigenvalue problem (Richardson-extrapolated from $64$ and $128$ cells per direction).}\label{tab:geom}\footnotesize
\begin{adjustbox}{max width=\linewidth}\begin{tabular}{lcccccc}\toprule
Domain & $\min\sqrt g$ & $\max\sqrt g$ & $\max\kappa$ & $C_\Omega^2$ (i) & $C_\Omega^2$ \eqref{eq:CP} & $1/\lambda_1$\\\midrule
Polar sector & 0.393 & 0.785 & 3.14 & -- & -- & --\\
Wavy (advection) & 0.858 & 1.142 & 2.21 & -- & -- & --\\
Annulus & 4.712 & 18.850 & 8.38 & 16.00 & 0.912 & 0.238\\
Sampled boundary (TFI) & 3.183 & 24.862 & 11.31 & 28.48 & 1.664 & 0.219\\
Wavy channel & 0.901 & 1.099 & 3.47 & 0.19 & 0.122 & 0.040\\
Star & 2.827 & 34.306 & 30.30 & 105.31 & 2.323 & 0.162\\
\bottomrule\end{tabular}\end{adjustbox}\end{table}

\begin{figure}[H]
	\centering
	\subfloat [$\sqrt{g}$]{\includegraphics[width=0.2\textwidth]{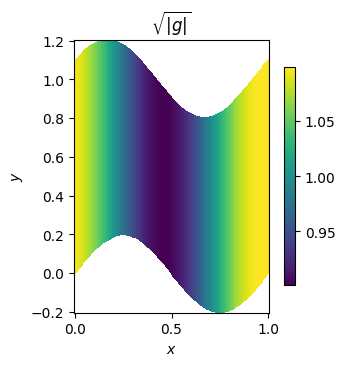}}
	\subfloat [$|\det(\mathbf J)|$]{\includegraphics[width=0.2\textwidth]{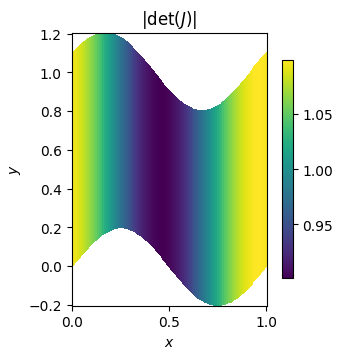}}
	\subfloat [$g_{11}$]{\includegraphics[width=0.2\textwidth]{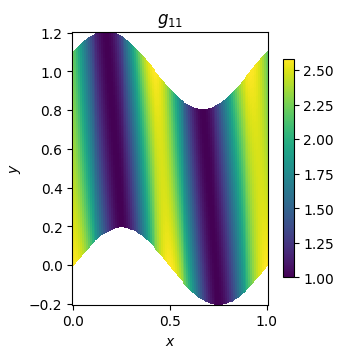}}
	\subfloat [$g_{12}$]{\includegraphics[width=0.2\textwidth]{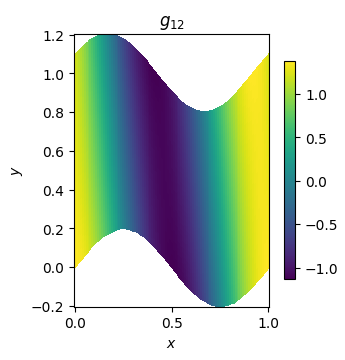}}
	\subfloat [$g_{22}$]{\includegraphics[width=0.2\textwidth]{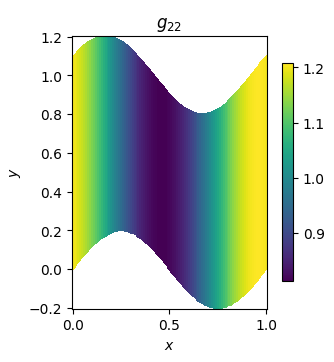}}
	\caption{Geometric quantities of the wavy-channel map computed by automatic differentiation (original CT-PIKAN implementation).}
	\label{fig:16}
\end{figure}

For the wavy heat channel, Fig.~\ref{fig:16} shows the geometric quantities as computed in the original CT-PIKAN implementation. They coincide with the closed forms of Appendix~\ref{app:metrics}: $\sqrt g=|\det\mathbf J|=y_\eta$, $g_{11}=1+y_\xi^2$, $g_{12}=y_\xi y_\eta$ and $g_{22}=y_\eta^2$.

\subsection{Hyperbolic problems: advection}
\label{sec:res_adv}
The two-dimensional advection equation describes the transport of a scalar quantity, such as temperature or concentration, by a prescribed velocity field. In applications such as meteorology, oceanography and chemical engineering, the geometry of the domain can strongly influence the advected quantity. We solve $u_t+\mathbf v\cdot\nabla u=0$ with inflow and initial data taken from the exact solution in three settings.
\begin{itemize}
\item \emph{Cartesian:} $\Omega=(0,1)^2$, $\mathbf v=(1,1)$ (Eq.~\eqref{eq:24} with $c=1$), $u=\sin 2\pi(x-t)\sin2\pi(y-t)$, $T=0.5$.
\item \emph{Polar sector:} $0.5<r<1$, $0<\theta<\pi/2$, $v_r=v_\theta=1$, with the initial condition \eqref{eq:27} for $\alpha=-1$, $\lambda=1$ and $\beta\in\{2,8\}$, and the exact solution \eqref{eq:exact_polar}. We use $T=0.45$, which keeps $r-t>0$ on the inflow circle. The map is $r=\tfrac12+\tfrac12\xi$, $\theta=\tfrac\pi2\eta$.
\item \emph{Wavy:} the map \eqref{eq:30} with $A=0.02$, $\mathbf v=(1,1)$, initial condition $u(x,y,0)=\sin2\pi(x+y)$ and exact solution $u=\sin2\pi(x+y-2t)$, $T=0.5$. These grids are standard test cases for curvilinear high-order schemes \cite{curve1,curve2,curve3}.
\end{itemize}
In the wavy case the contravariant velocity $\mathbf w=\mathbf J^{-1}\mathbf v$, Eq.~\eqref{eq:40}, varies in space and is evaluated pointwise by autograd. In all cases the inflow faces are $\xi=0$ and $\eta=0$. For the wavy map this holds because $6\pi A<1$, so that $\mathbf v\cdot\mathbf n<0$ on the left and lower curves. The outflow faces satisfy $\mathbf v\cdot\mathbf n\ge0$, as Theorem~\ref{thm:apost}(c) requires. The velocity is divergence-free except in the polar case, where $\nabla\cdot\mathbf v=1/r\le2$, so $\beta=2$ there.

\begin{figure}[H]
\centering
\includegraphics[width=0.85\textwidth]{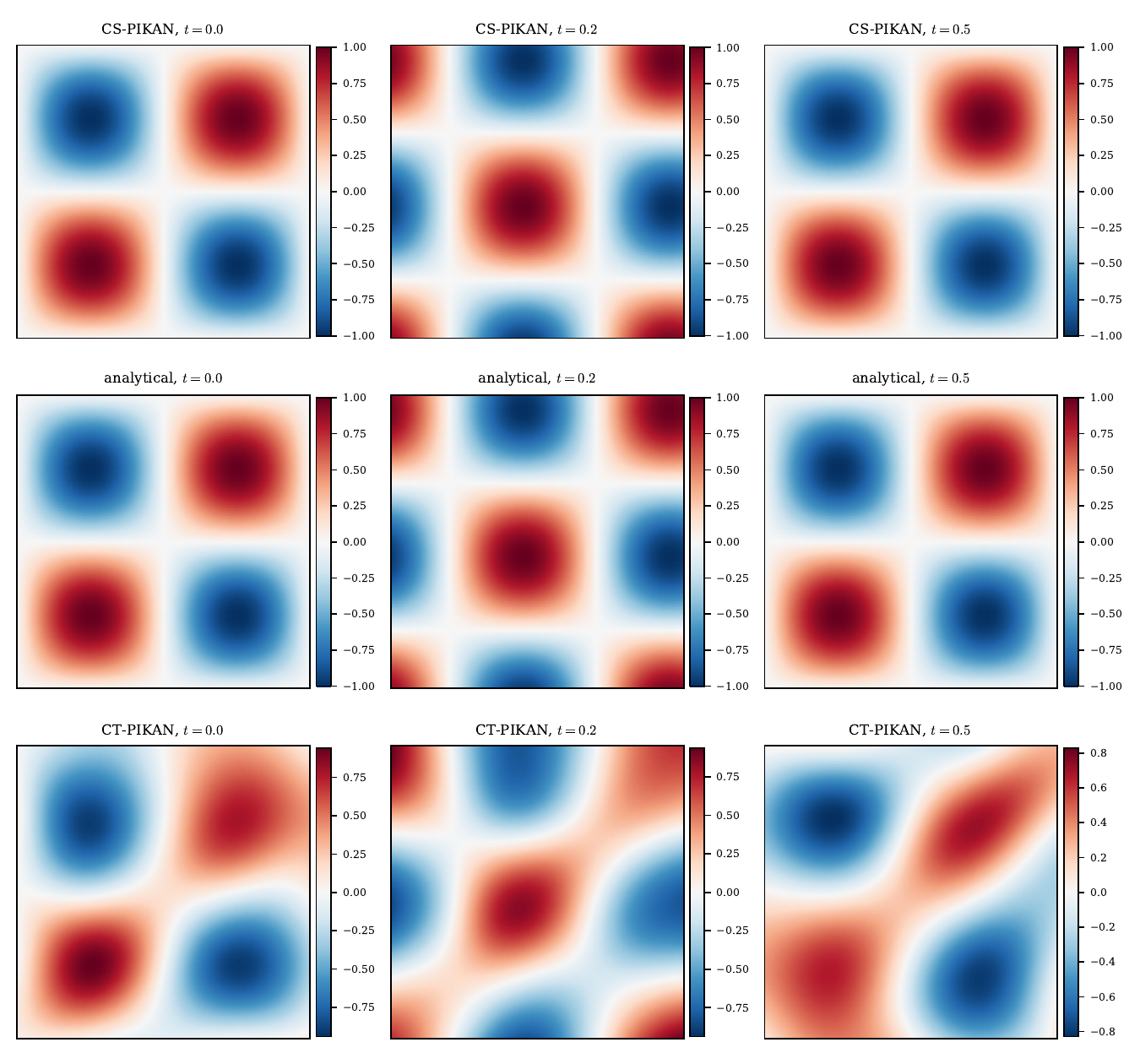}
\caption{Advection on the Cartesian domain at $t=0$, $0.2$ and $0.5$: CS-PIKAN (top), analytical solution (middle), CT-PIKAN (bottom); seed 0.}
\label{fig:adv_cart}
\end{figure}

\paragraph{Cartesian domain.} Fig.~\ref{fig:adv_cart} compares CS-PIKAN, the analytical solution and CT-PIKAN at three times. CS-PIKAN is visually indistinguishable from the analytical field, with a relative $L_2$ error of $2.88\times10^{-3}$ (Table~\ref{tab:main}). CT-PIKAN distorts and damps the pattern as it is transported (relative $L_2$ error $0.29$). The point-wise baseline thus fails already on the Cartesian problem once the inflow data have to be imposed by penalties over the whole time interval.

\begin{figure}[H]
\centering
\includegraphics[width=\textwidth]{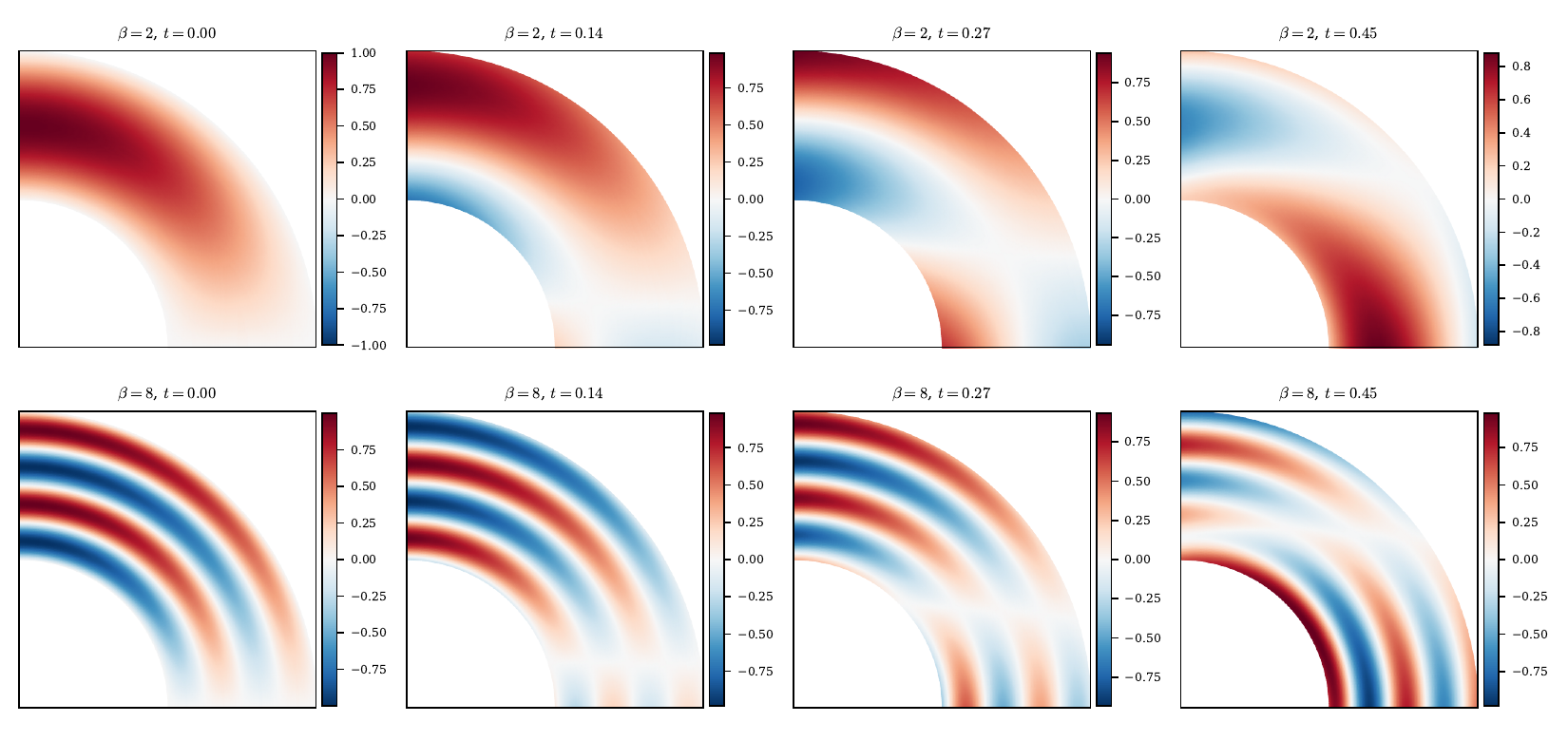}
\caption{CS-PIKAN solutions of the advection equation on the quarter-circle sector with the sinusoidal initial condition \eqref{eq:27}: $\beta=2$ (top) and $\beta=8$ (bottom), at $t=0$, $0.135$, $0.27$ and $0.45$ (seed 0).}
\label{fig:adv_polar}
\end{figure}

\paragraph{Polar sector.} Fig.~\ref{fig:adv_polar} shows the CS-PIKAN solutions for $\beta=2$ and $\beta=8$. The radial pattern is transported outwards while the angular phase is sheared by the factor $\ln\frac{r}{r-t}$. For $\beta=8$, four radial wavelengths are transported obliquely to the computational coordinates, which makes this the hardest transport benchmark for a separable representation (relative $L_2$ error $1.62\times10^{-2}$). CT-PIKAN fails on this case ($0.64$).

\begin{figure}[H]
\centering
\includegraphics[width=\textwidth]{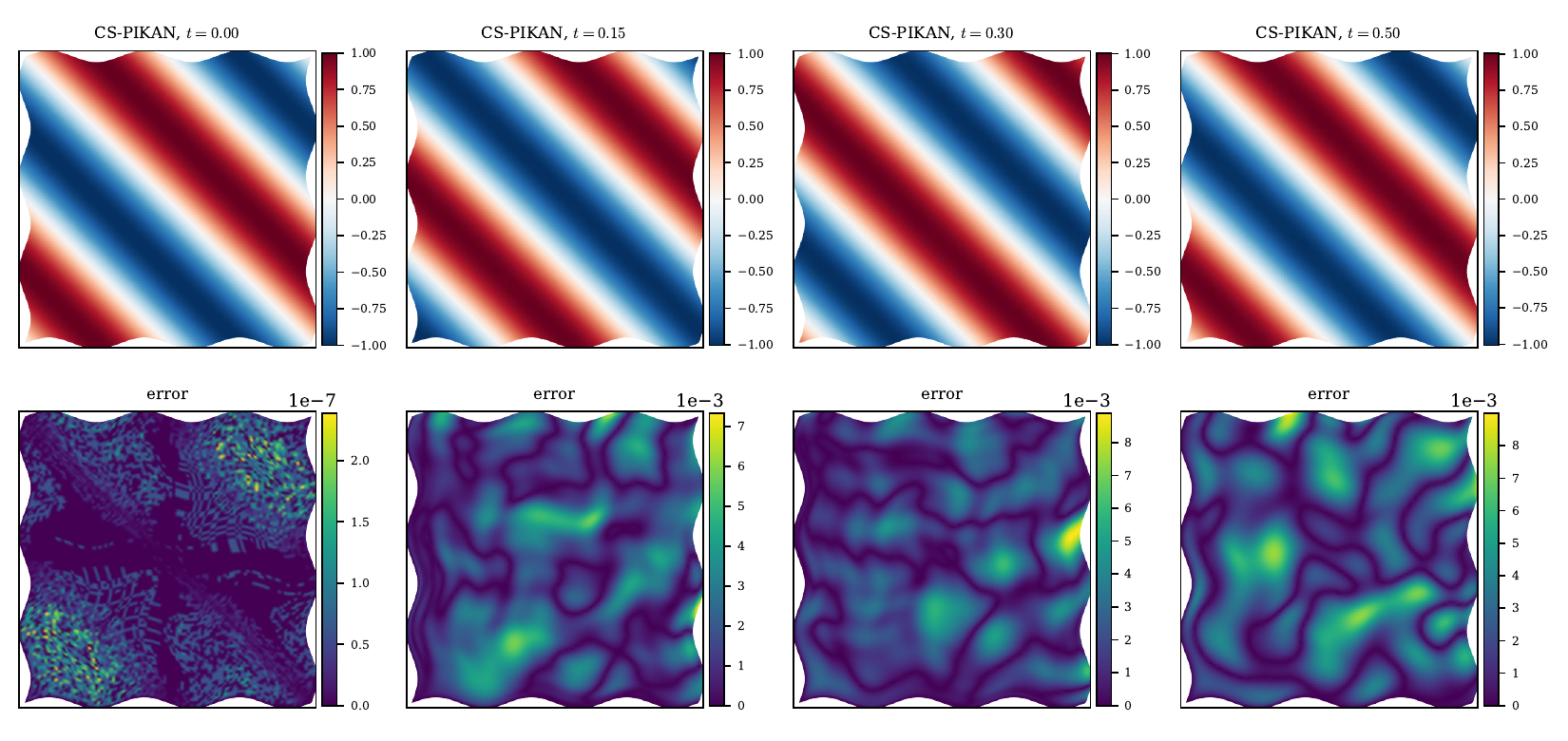}
\caption{CS-PIKAN solution (top) and pointwise error (bottom) of the advection equation on the wavy domain at $t=0$, $0.15$, $0.3$ and $0.5$ (seed 0).}
\label{fig:adv_wavy}
\end{figure}

\begin{figure}[H]
\centering
\includegraphics[width=\textwidth]{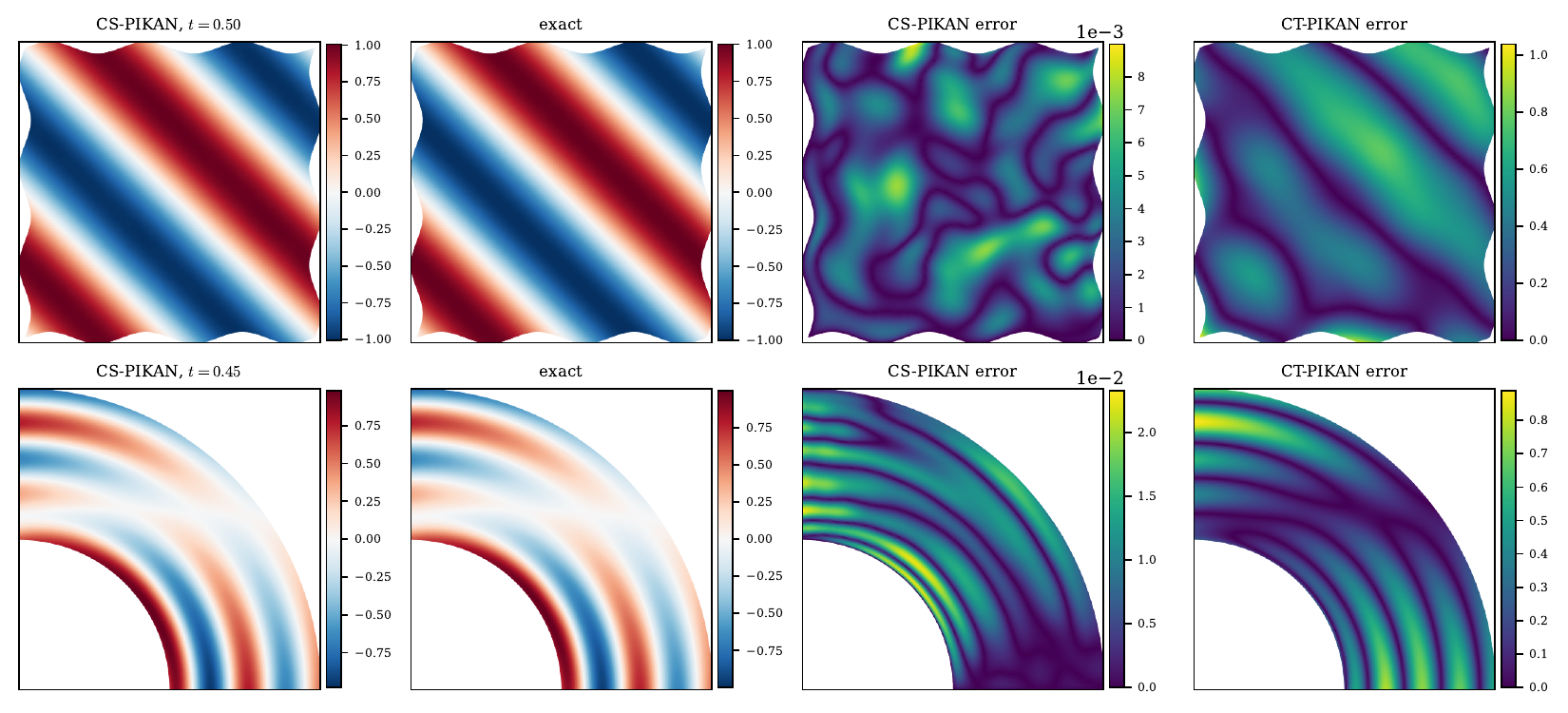}
\caption{Advection at the final time. Columns: CS-PIKAN solution, exact solution, CS-PIKAN error $|u_\theta-u|$, and CT-PIKAN error (same seed; note the different colour scales). Rows: wavy domain; polar sector with $\beta=8$.}
\label{fig:adv}
\end{figure}

\paragraph{Wavy domain.} Fig.~\ref{fig:adv_wavy} shows the transported pattern and the pointwise CS-PIKAN error on the non-orthogonal wavy domain. The maximum pointwise error remains below $1.2\%$ of the unit amplitude over all seeds and times. Fig.~\ref{fig:adv} compares the error maps of the two methods at the final time for the wavy and the polar ($\beta=8$) cases.

\subsection{Elliptic problems: Poisson on an annulus and on a sampled-boundary domain}
\label{sec:res_poisson}
To evaluate elliptic problems on non-Cartesian geometries, we consider the two-dimensional Poisson equation \cite{Poisson1,Poisson2}
\begin{equation}
	\nabla^{2}u+f(\mathbf{x})=0, \qquad \mathbf{x}=(x,y)\in\Omega,
	\label{eq:41}
\end{equation}
with homogeneous Dirichlet conditions on both boundary components and the Gaussian source
\begin{equation}
	f(x,y)= A \exp\left( -\frac{(x-x_0)^2+(y-y_0)^2}{2\sigma^2} \right),
	\label{eq:46}
\end{equation}
where $(x_0,y_0)=(-1,0)$ denotes the source location, $\sigma=1.0$ is the Gaussian width, and $A=2.0$ is the source amplitude. The localized forcing produces a smooth solution with sufficient spatial variation to test the approximation. Two domains are considered.
\begin{itemize}
\item \emph{Annulus:} $\Omega=\{(x,y):r_{\mathrm{inner}}\le\sqrt{x^{2}+y^{2}}\le r_{\mathrm{outer}}\}$ with $r_{\mathrm{inner}}=0.5$ and $r_{\mathrm{outer}}=2.0$, mapped by
\begin{equation}
	r(\xi) = r_{\mathrm{inner}} +\xi(r_{\mathrm{outer}}-r_{\mathrm{inner}}), \qquad \theta(\eta) = 2\pi\eta,\qquad x = r\cos\theta,\quad y = r\sin\theta .
	\label{eq:43}
\end{equation}
\item \emph{Sampled boundary:} a doubly connected domain whose inner and outer boundaries are known only through $64$ sample points each (Fig.~\ref{fig:geoms}g), mapped by the transfinite blend \eqref{eq:tfi_map} of the trigonometric interpolants of the samples.
\end{itemize}
In both cases $\eta$ is periodic, the Laplacian is evaluated with the Laplace--Beltrami operator \eqref{eq:12}, and the solver code is identical except for $\Phi$.

\begin{figure}[H]
	\centering
	\subfloat [CT-PIKAN solution]{\includegraphics[width=0.24\textwidth]{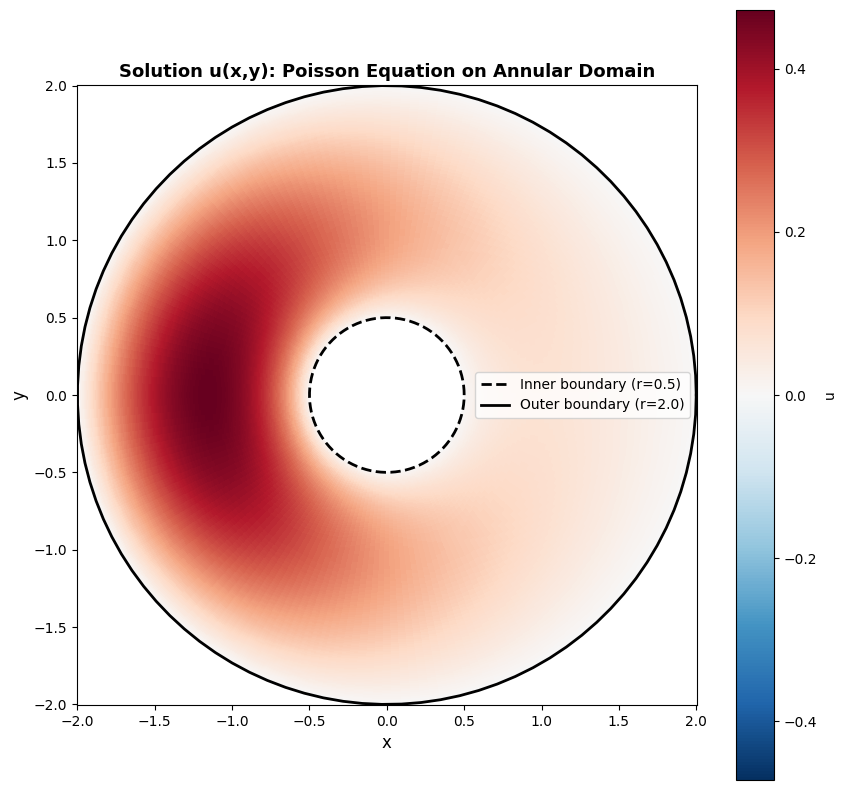}}
	\subfloat [Gaussian source term]{\includegraphics[width=0.24\textwidth]{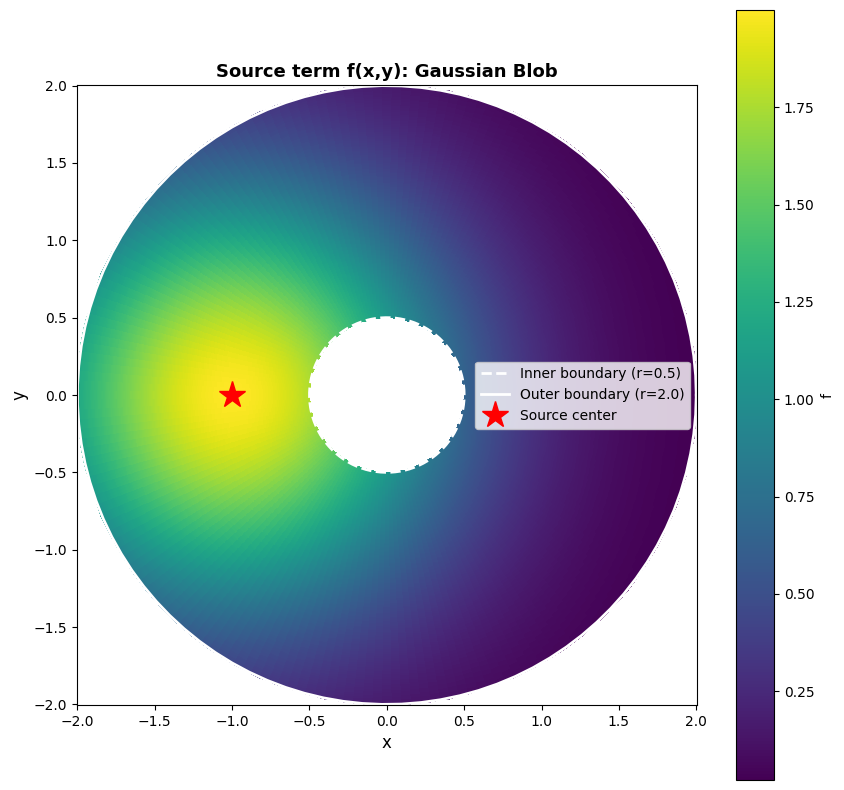}}
	\subfloat [CT-PIKAN residual]{\includegraphics[width=0.26\textwidth]{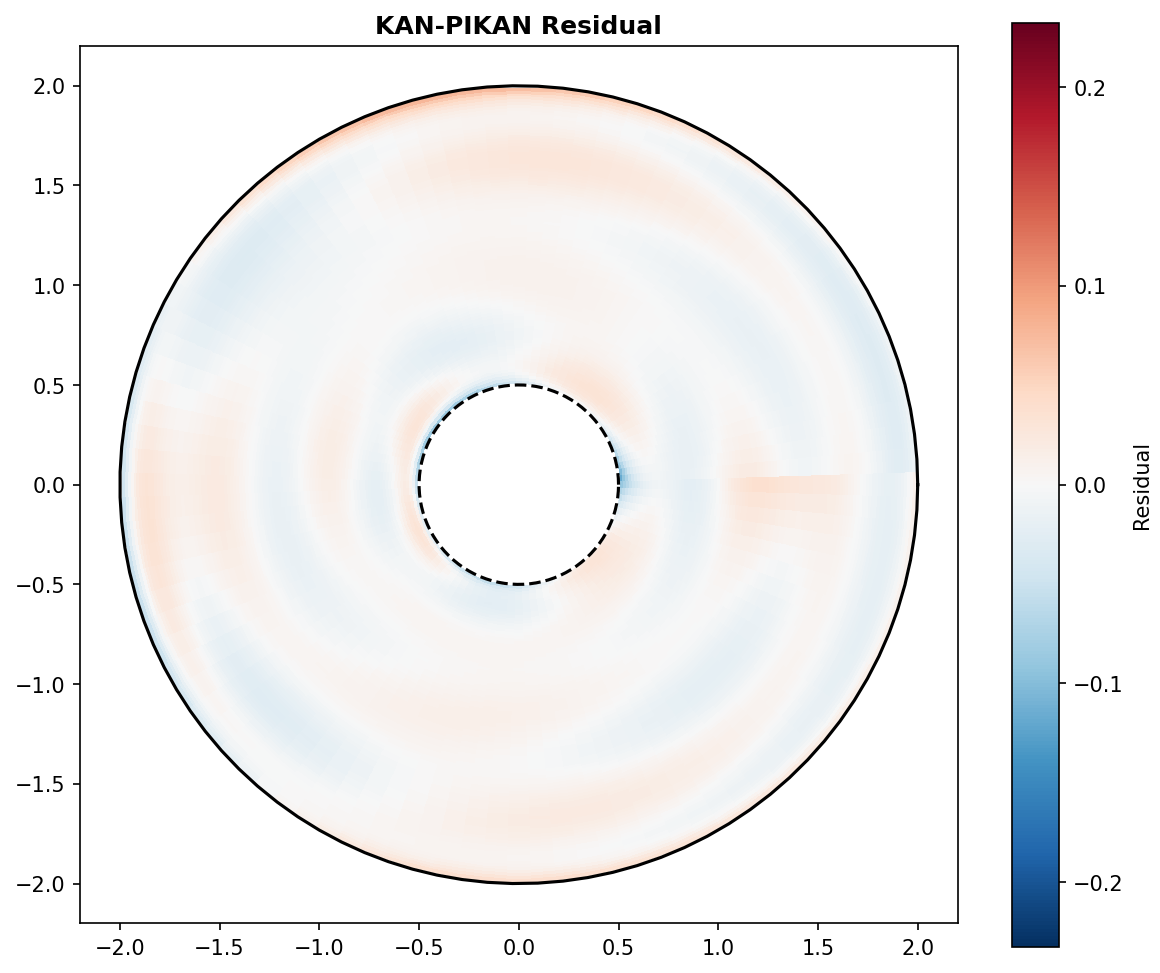}}
	\subfloat [Vanilla PINN residual]{\includegraphics[width=0.26\textwidth]{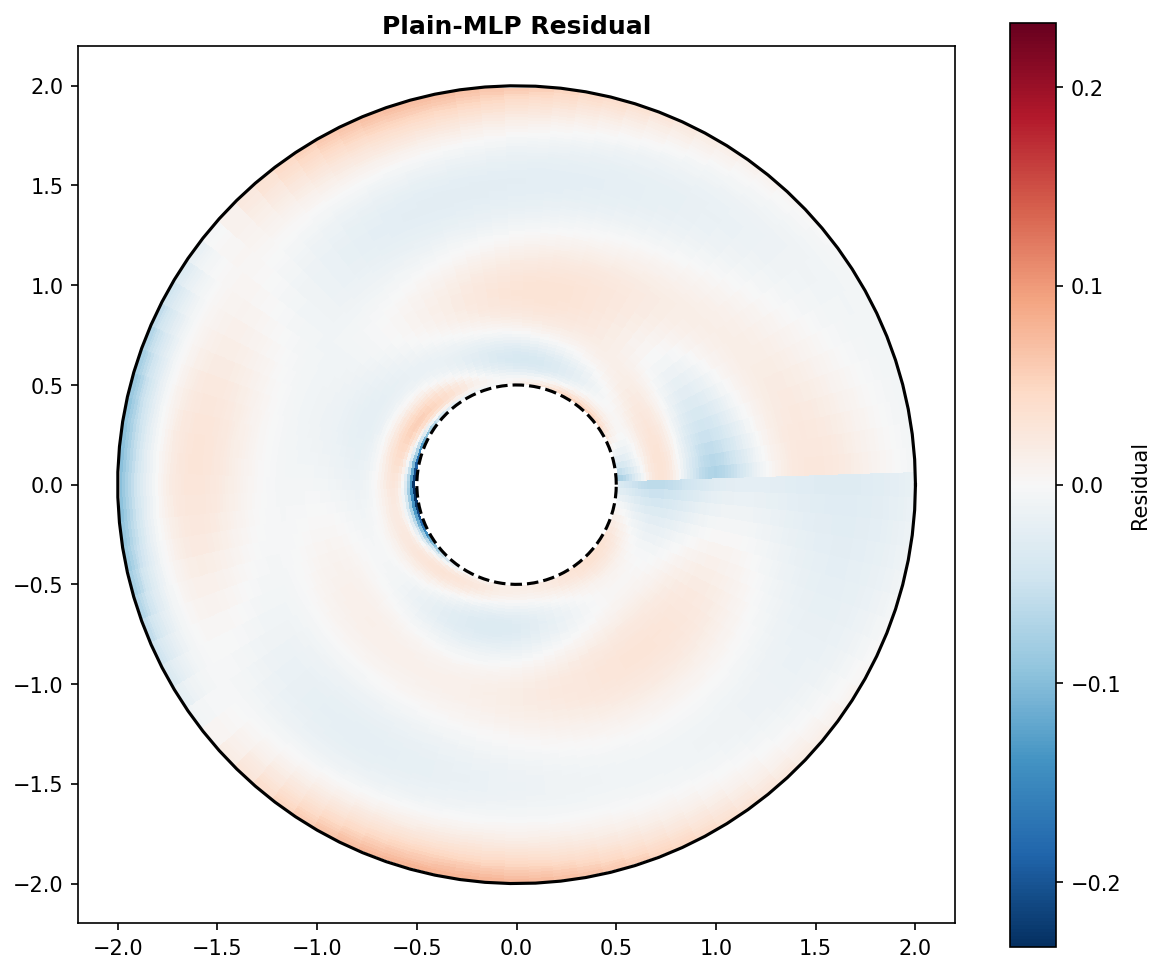}}
	\caption{Annular Poisson problem, original CT-PIKAN implementation: solution, Gaussian source and PDE residuals of CT-PIKAN and of a vanilla PINN.}
	\label{fig:8}
\end{figure}

\begin{figure}[H]
	\centering
	\subfloat [3D view of the CT-PIKAN solution]{\includegraphics[trim={0 0 0 0.65cm},clip,width=0.5\textwidth]{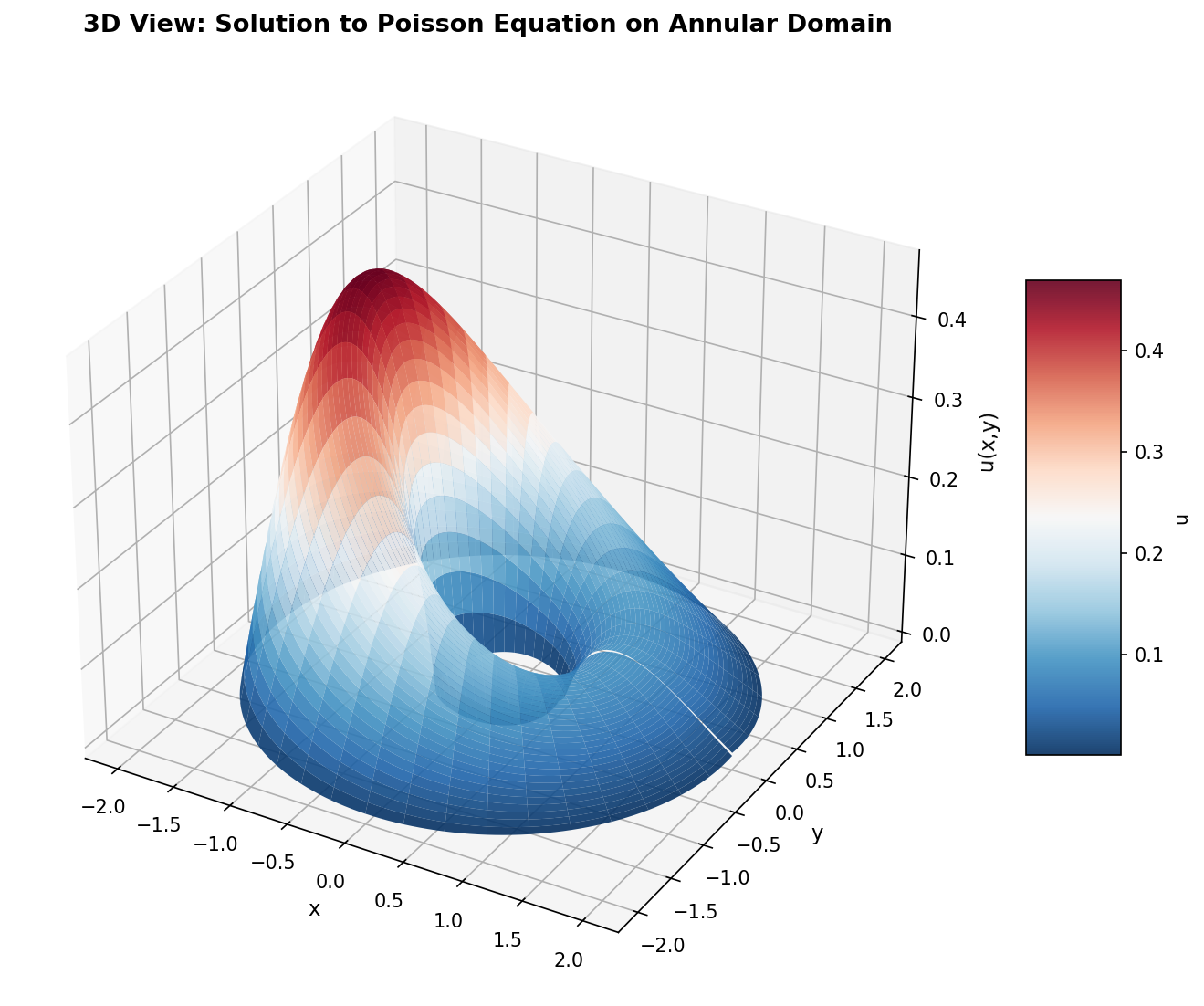}}
	\qquad \quad
	\subfloat [Area element $\sqrt{g}$]{\includegraphics[width=0.40\textwidth]{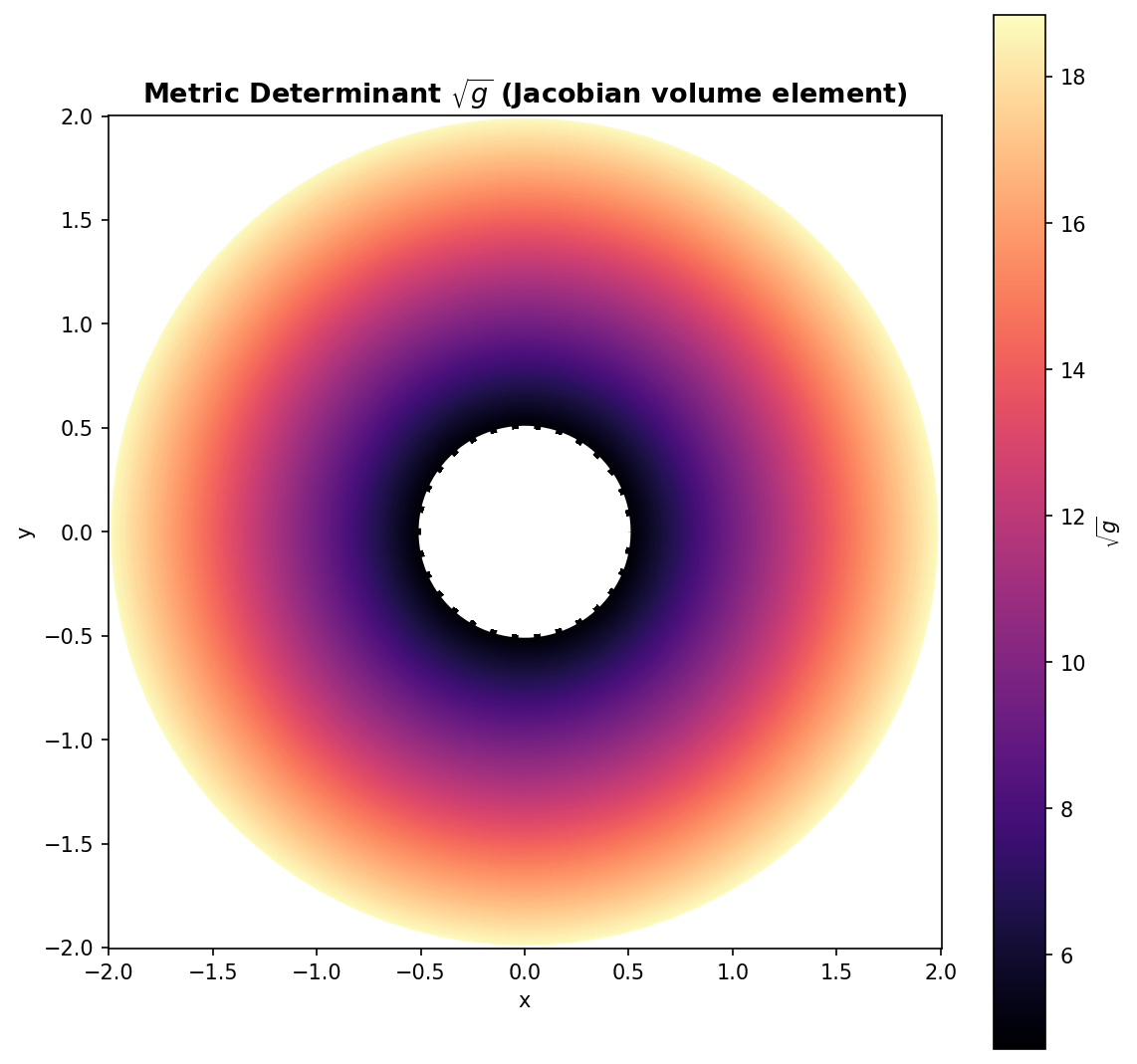}}
	\caption{Three-dimensional view of the CT-PIKAN solution and the area element of the annulus map.}
	\label{fig:9}
\end{figure}

\begin{figure}[H]
	\centering
	\subfloat [CT-PIKAN convergence]{\includegraphics[width=0.45\textwidth]{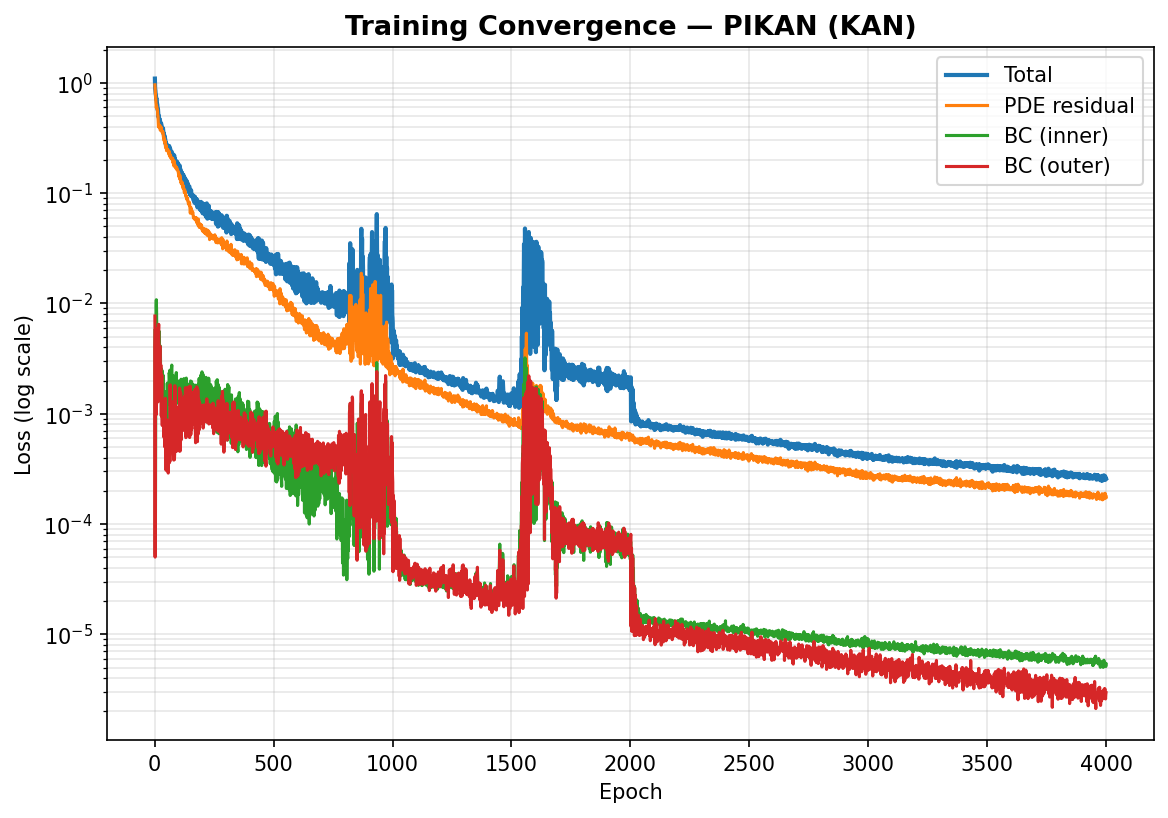}}
	\qquad \quad
	\subfloat [KAN versus MLP front end (loss)]{\includegraphics[width=0.45\textwidth]{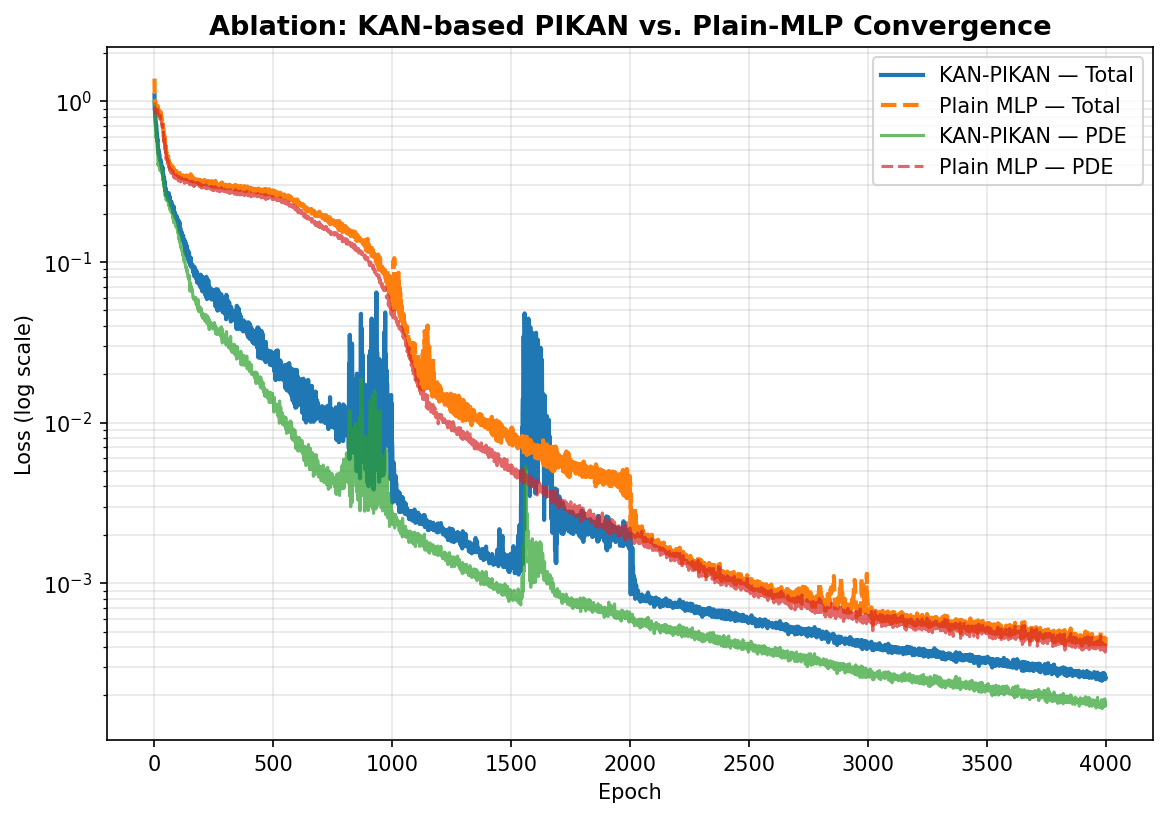}}
	\caption{Training convergence of CT-PIKAN and the loss-based KAN/MLP comparison of the original study. Losses of different formulations are not directly comparable; Table~\ref{tab:ablation} compares errors instead.}
	\label{fig:10}
\end{figure}

\paragraph{CT-PIKAN diagnostics (original implementation).} Figs.~\ref{fig:8}--\ref{fig:10} show the annular Poisson problem as solved by the original PyTorch implementation of CT-PIKAN: the solution, the source, the PDE residual of CT-PIKAN and of a vanilla PINN, a three-dimensional view, the area element $\sqrt g$, the training history, and the KAN-versus-MLP loss comparison of the original study. The solution peaks near the source and decays towards both circles. Residual maps alone, however, do not measure the error: the quantitative comparison below shows that CT-PIKAN's relative error is $2.07\times10^{-2}$, mostly concentrated along the seam.

\begin{figure}[H]
\centering
\includegraphics[width=\textwidth]{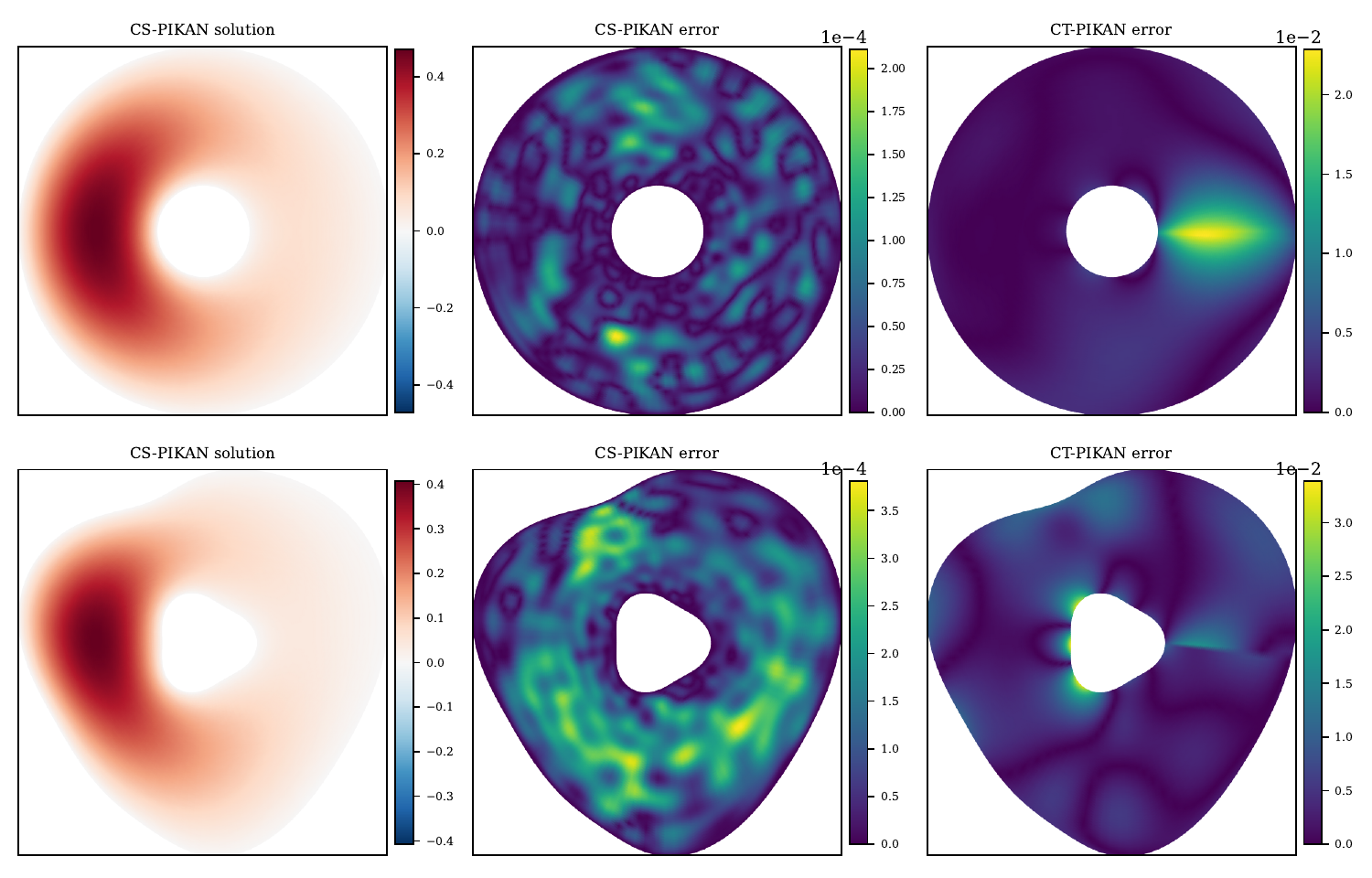}
\caption{Poisson problems. Columns: CS-PIKAN solution, CS-PIKAN error, CT-PIKAN error (same seed). Rows: annulus; sampled-boundary (TFI) domain. On the annulus the CT-PIKAN error concentrates along the ray $\theta=0$, i.e.\ along the seam $\eta=0\equiv1$ that its raw $\eta$ input leaves unconstrained.}
\label{fig:poisson}
\end{figure}

\paragraph{Error comparison.} Fig.~\ref{fig:poisson} compares the CS-PIKAN solution and error with the CT-PIKAN error on both domains. CS-PIKAN reaches relative $L_2$ errors of $6.30\times10^{-4}$ (annulus) and $9.32\times10^{-4}$ (sampled boundary), against $2.07\times10^{-2}$ and $3.91\times10^{-2}$ for CT-PIKAN. 

\begin{figure}[H]
\centering
\includegraphics[width=0.9\textwidth]{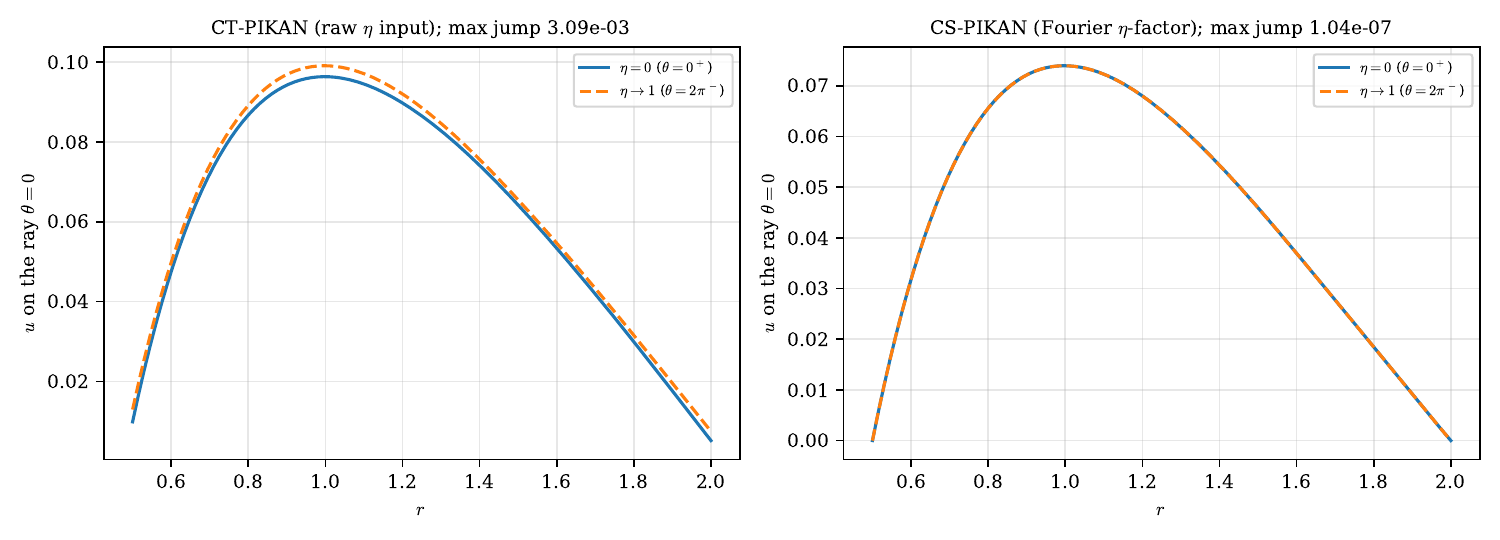}
\caption{Seam diagnostics for the annular Poisson problem (seed 0): the solution along the ray $\theta=0$ evaluated from both sides of the computational square, $\eta=0$ and $\eta\to1$. Left: CT-PIKAN, whose two traces differ. Right: CS-PIKAN, whose traces coincide.}
\label{fig:seam}
\end{figure}

On the annulus the CT-PIKAN error concentrates along the ray $\theta=0$, which is the image of the seam $\eta=0\equiv1$. Fig.~\ref{fig:seam} makes this explicit: along the seam, the two one-sided traces of the CT-PIKAN solution differ by up to $3.09\times10^{-3}$ and overshoot the reference by about $30\%$, whereas the Fourier factor of CS-PIKAN makes the two traces identical up to round-off ($1.04\times10^{-7}$).

\subsection{Parabolic problems: transient heat conduction}
\label{sec:res_heat}
We extend the comparison to the transient heat (diffusion) equation \cite{Transform1,Heat1}
\begin{equation}
	\frac{\partial u}{\partial t} = \alpha \, \nabla^2 u(x,y,t) + f(x,y,t),
	\qquad (x,y) \in \Omega,\; t \in [0,T],
	\label{eq:47}
\end{equation}
where $\alpha$ is the thermal diffusivity and $f$ is an optional source term, set to zero (pure diffusion) in all cases below. As before, $\Omega$ is the image of $\Q$ under a domain-specific map, and the metric is obtained by automatic differentiation of $\Phi$. In CT-PIKAN the network takes the augmented input $(\xi,\eta,t)$, the Laplace--Beltrami operator acts on the spatial sub-vector, and $u_t$ is obtained by a separate automatic-differentiation call. In CS-PIKAN the time factor $C_r(t)$ of Eq.~\eqref{eq:ansatz} provides $u_t$ through Proposition~\ref{prop:sep}.

\paragraph{Wavy channel.} The first case is a channel with sinusoidally corrugated top and bottom walls,
\begin{equation}
\begin{aligned}
&x = L_x\,\xi,\qquad y = (1-\eta)\, y_{\text{bottom}}(x) + \eta\, y_{\text{top}}(x),\\
&y_{\text{bottom}}(x) = A\sin(2\pi x), \qquad y_{\text{top}}(x) = L_y + A\sin(2\pi x + \phi),
\end{aligned}
	\label{eq:48}
\end{equation}
with amplitude $A = 0.2$, phase offset $\phi = 0.5$, and $L_x = L_y = 1$. The two walls are out-of-phase sinusoids, so the channel width varies along its length while remaining strictly positive ($\sqrt g\in[0.901,1.099]$). The diffusivity is $\alpha = 0.05$ and the horizon $T = 1$. Homogeneous Dirichlet conditions are imposed on all four edges of the computational square ($\xi = 0$, $\xi = 1$, $\eta = 0$, $\eta = 1$, corresponding to the inlet, the outlet and the two wavy walls). The initial condition is defined in computational coordinates,
\begin{equation}
	u_0\big(x(\xi,\eta),y(\xi,\eta)\big) = \sin(\pi\xi)\,\sin(\pi\eta),
	\label{eq:52}
\end{equation}
which vanishes on all four edges and therefore on both wavy walls, so the initial and boundary data are compatible.

\begin{figure}[H]
	\centering
	\subfloat [$t=0$]{\includegraphics[width=0.24\textwidth]{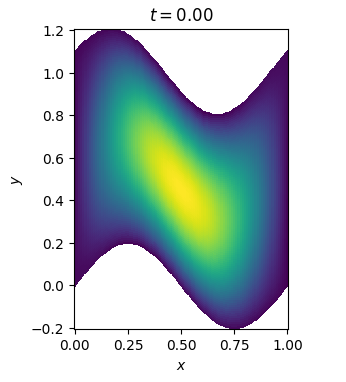}}
	\subfloat [$t=0.33$]{\includegraphics[width=0.238\textwidth]{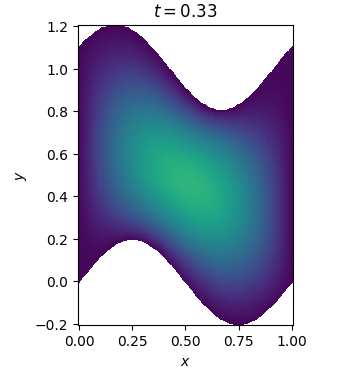}}
	\subfloat [$t=0.67$]{\includegraphics[width=0.24\textwidth]{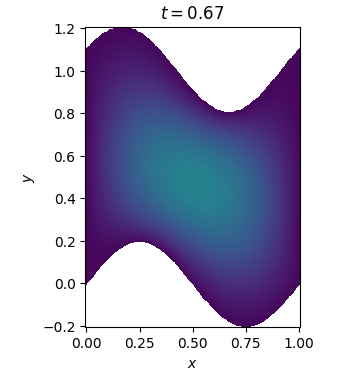}}
	\subfloat [$t=1$]{\includegraphics[width=0.27\textwidth]{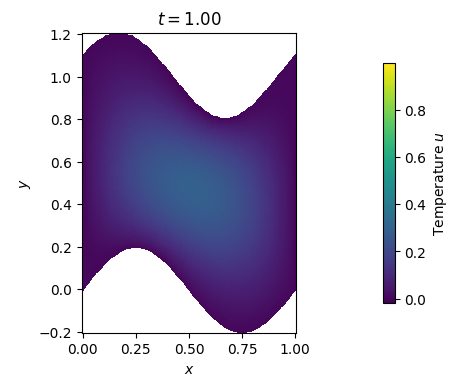}}
	\caption{Transient temperature evolution on the wavy channel, original CT-PIKAN implementation.}
	\label{fig:12}
\end{figure}

\begin{figure}[H]
	\centering
	\subfloat [Interior space--time collocation points]{\includegraphics[trim={0 0 0 0.75cm},clip,width=0.45\textwidth]{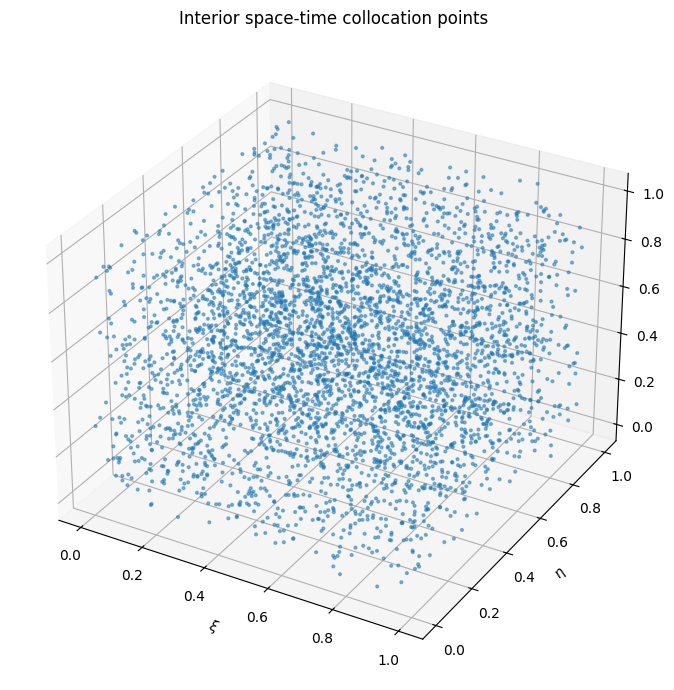}}
	\subfloat [Temperature field at $t=0.5$]{\includegraphics[trim={0 0 0 0.75cm},clip,width=0.5\textwidth]{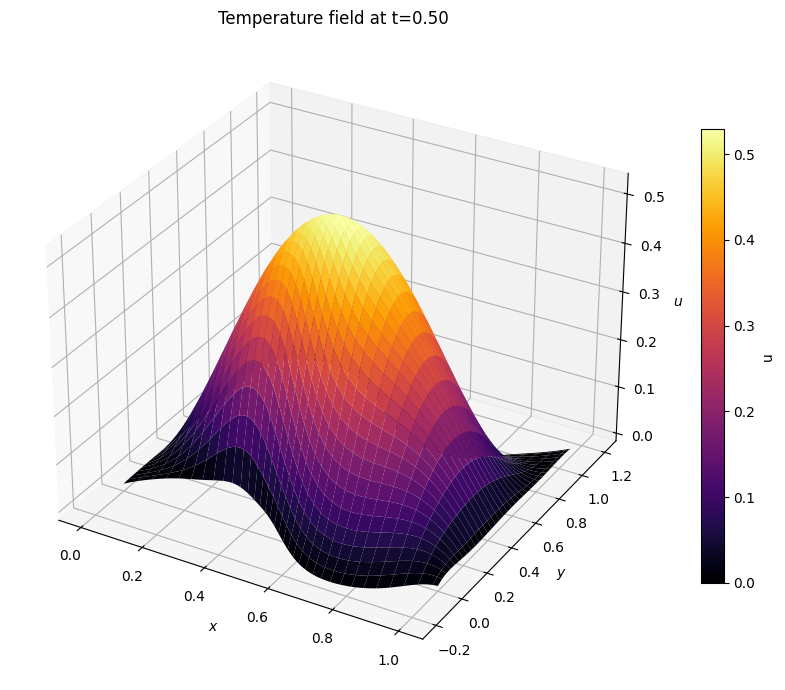}}
	\caption{Three-dimensional view of the CT-PIKAN space--time collocation points and of the CT-PIKAN solution on the wavy channel.}
	\label{fig:13}
\end{figure}

\begin{figure}[H]
	\centering
	\subfloat [Training convergence]{\includegraphics[width=0.41\textwidth]{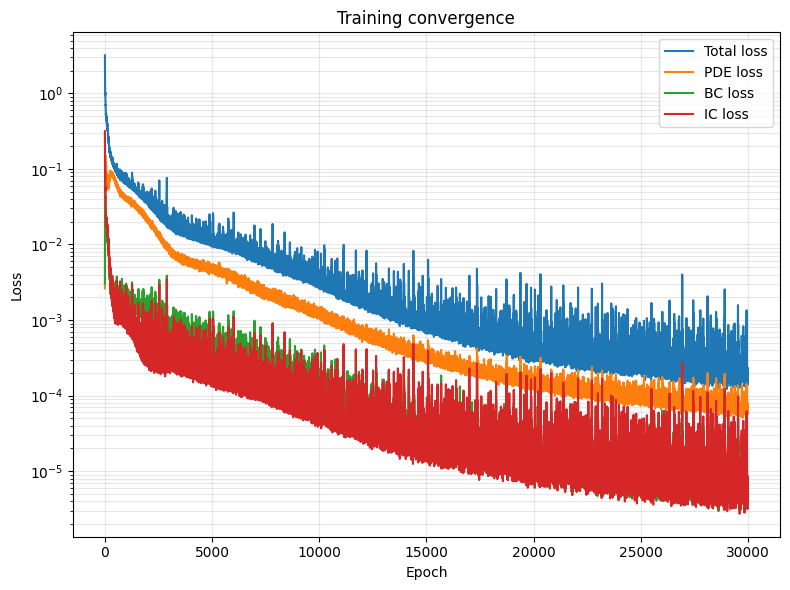}}
	\subfloat [Thermal energy decay $E(t)=\int_{\Omega} u^{2}\,d\Omega$]{\includegraphics[width=0.50\textwidth]{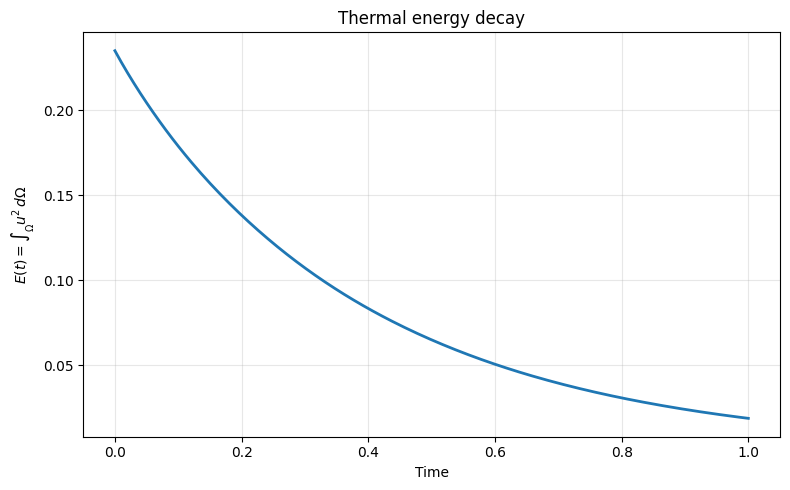}}
	\caption{Training convergence \& thermal-energy decay of CT-PIKAN on the wavy channel.}
	\label{fig:14}
\end{figure}

\begin{figure}[H]
	\centering
	\subfloat [Predicted I.C.]{\includegraphics[width=0.25\textwidth]{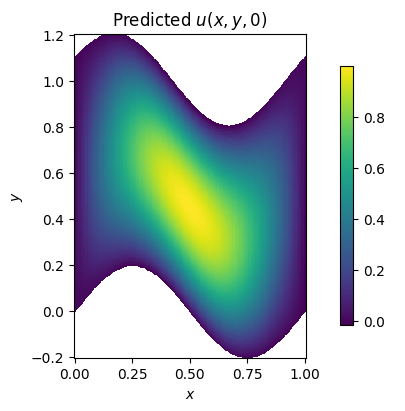}}
	\subfloat [Prescribed I.C.]{\includegraphics[width=0.25\textwidth]{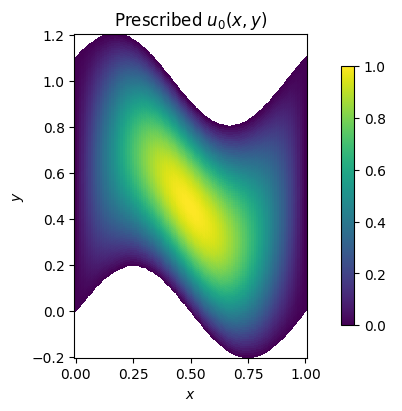}}
	\subfloat [I.C.~error]{\includegraphics[width=0.25\textwidth]{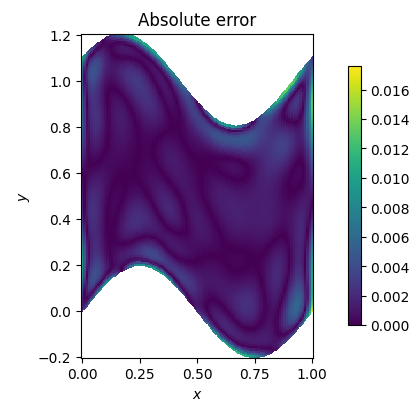}}
	\\
	\subfloat [PDE residual ($t=0$)]{\includegraphics[width=0.2\textwidth]{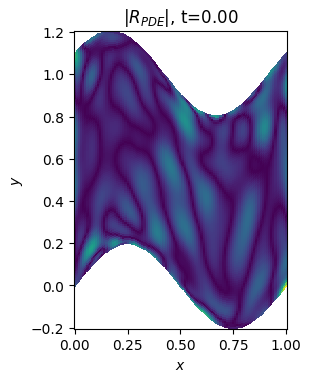}}
	\subfloat [$t=0.33$]{\includegraphics[width=0.2\textwidth]{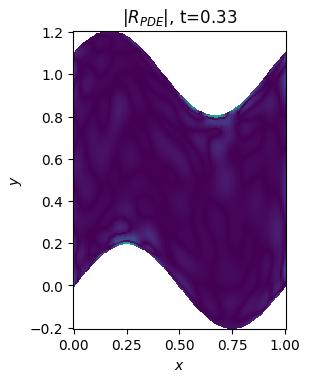}}
	\subfloat [$t=0.67$]{\includegraphics[width=0.2\textwidth]{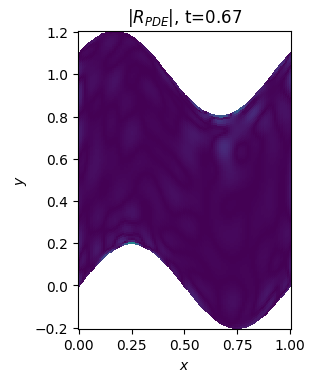}}
	\subfloat [$t=1$]{\includegraphics[width=0.28\textwidth]{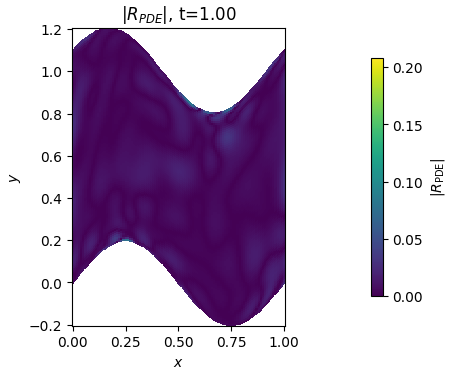}}
	\\
	\subfloat [Left B.C.~error]{\includegraphics[width=0.25\textwidth]{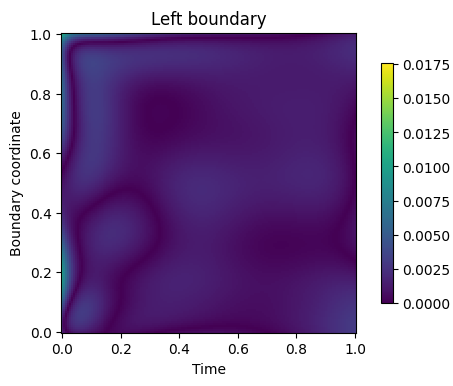}}
	\subfloat [Right B.C.~error]{\includegraphics[width=0.25\textwidth]{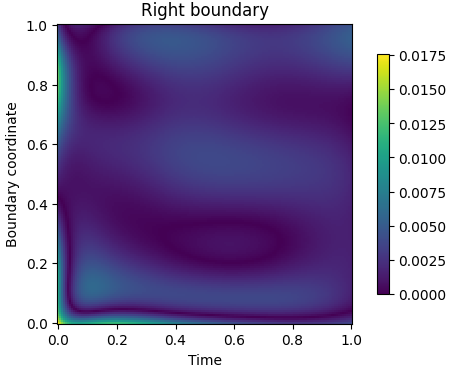}}
	\subfloat [Bottom B.C.~error]{\includegraphics[width=0.25\textwidth]{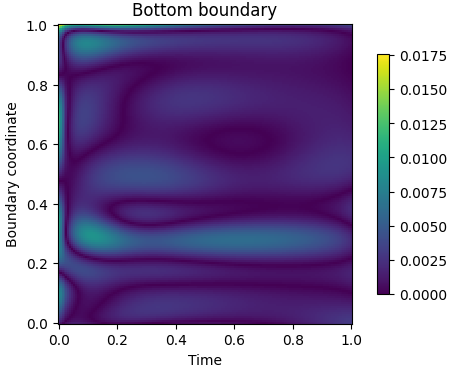}}
	\subfloat [Top B.C.~error]{\includegraphics[width=0.25\textwidth]{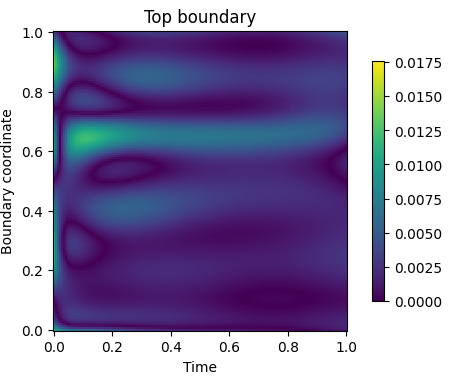}}
	\caption{Penalty-based constraints in CT-PIKAN.}
	\label{fig:15}
\end{figure}

\paragraph{CT-PIKAN diagnostics (baseline implementation).} Figs.~\ref{fig:12}--\ref{fig:15} show the wavy-channel problem as solved by the original PyTorch implementation of CT-PIKAN: the temperature evolution, the space--time collocation points, the training history and thermal-energy decay, and maps of the initial-condition error, the PDE residual and the four boundary-condition errors. The solution shows the expected smooth diffusion, and the thermal energy decreases monotonically. Fig.~\ref{fig:15} also shows the price of penalty-based constraints: the initial and boundary conditions are satisfied only approximately, with nonzero error maps on all four edges and at $t=0$. These are exactly the errors that the exact constraints of CS-PIKAN remove.

\begin{figure}[H]
\centering
\includegraphics[width=\textwidth]{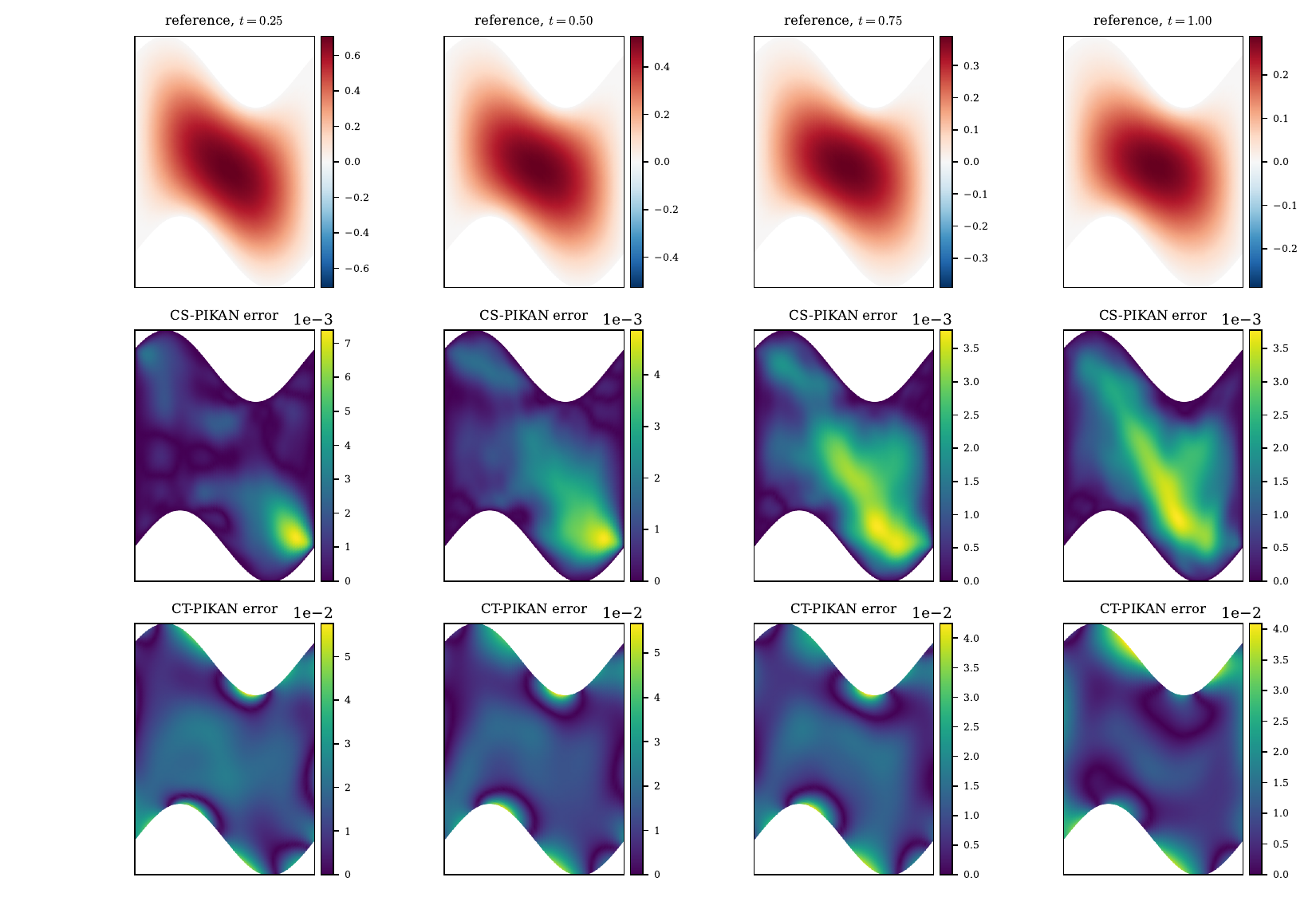}
\caption{Wavy-channel heat problem: reference solution (top), CS-PIKAN error (middle) and CT-PIKAN error (bottom) at $t=0.25$, $0.5$, $0.75$ and $1$ (seed 0; note the different colour scales).}
\label{fig:channel_err}
\end{figure}

In the above contours, initial-condition error, PDE residual at four times, and boundary-condition errors on the four edges of the wavy channel have been presented. In CS-PIKAN the initial and boundary errors vanish identically.
\paragraph{Error comparison on the channel.} Fig.~\ref{fig:channel_err} compares the pointwise errors of the two methods against the reference at four times. The CS-PIKAN error is smooth and of order $10^{-3}$; the CT-PIKAN error is about an order of magnitude larger and concentrated near the walls, where the penalised boundary data are imposed.

\paragraph{Star-shaped domain.} The second and third cases use a five-pointed star obtained by angularly modulating the outer radius of a polar map:
\begin{equation}
\begin{aligned}
&\theta(\eta) = 2\pi\eta,\qquad r_{\text{outer}}(\theta) = r_{\text{outer,base}} \left(1 + \gamma \cos(n_p\,\theta)\right),\\
&r(\xi,\theta) = r_{\text{inner}} + \xi\left(r_{\text{outer}}(\theta) - r_{\text{inner}}\right),
\end{aligned}
\label{eq:53}
\end{equation}
with $x = r\cos\theta$, $y = r\sin\theta$, inner radius $r_{\text{inner}} = 0.5$, base outer radius $r_{\text{outer,base}} = 2.0$, $n_p = 5$ star points and modulation amplitude $\gamma = 0.3$. As for the other maps, the metric tensor is obtained by automatic differentiation, with no manual derivation of curvature or Jacobian terms. The domain is doubly connected, so a well-posed heat problem requires a condition on \emph{both} boundary components. 

\begin{figure}[H]
\centering
\includegraphics[width=\textwidth]{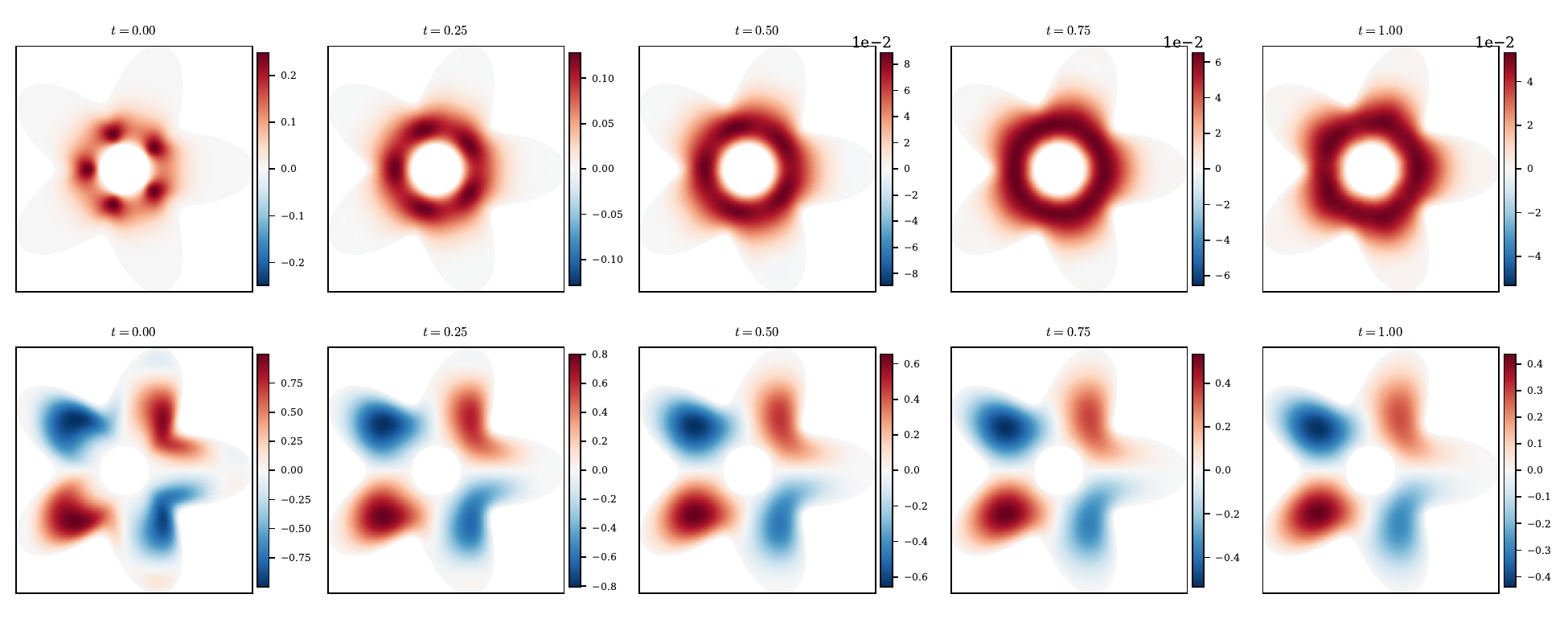}
\caption{CS-PIKAN temperature evolution on the star-shaped domain at $t=0$, $0.25$, $0.5$, $0.75$ and $1$: Gaussian initial condition (top) and sinusoidal initial condition (bottom); seed 0.}
\label{fig:star_snap}
\end{figure}

\begin{figure}[H]
\centering
\includegraphics[width=\textwidth]{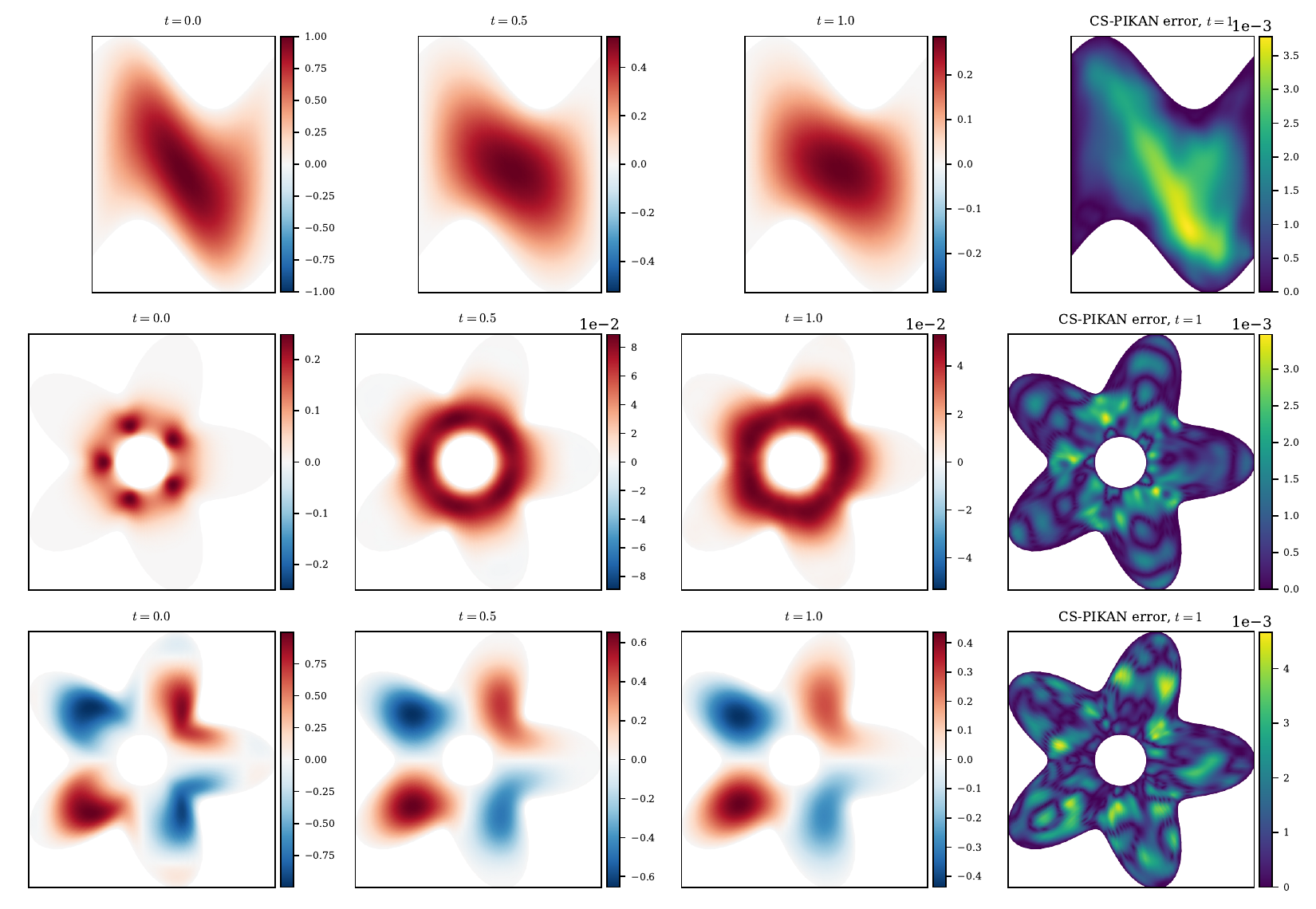}
\caption{Heat conduction. CS-PIKAN temperature at $t=0$, $0.5$ and $1$, and the CS-PIKAN error at $t=1$. Rows: wavy channel; star with the Gaussian initial condition; star with the sinusoidal initial condition.}
\label{fig:heat}
\end{figure}
\begin{figure}[H]
\centering
\includegraphics[width=\textwidth]{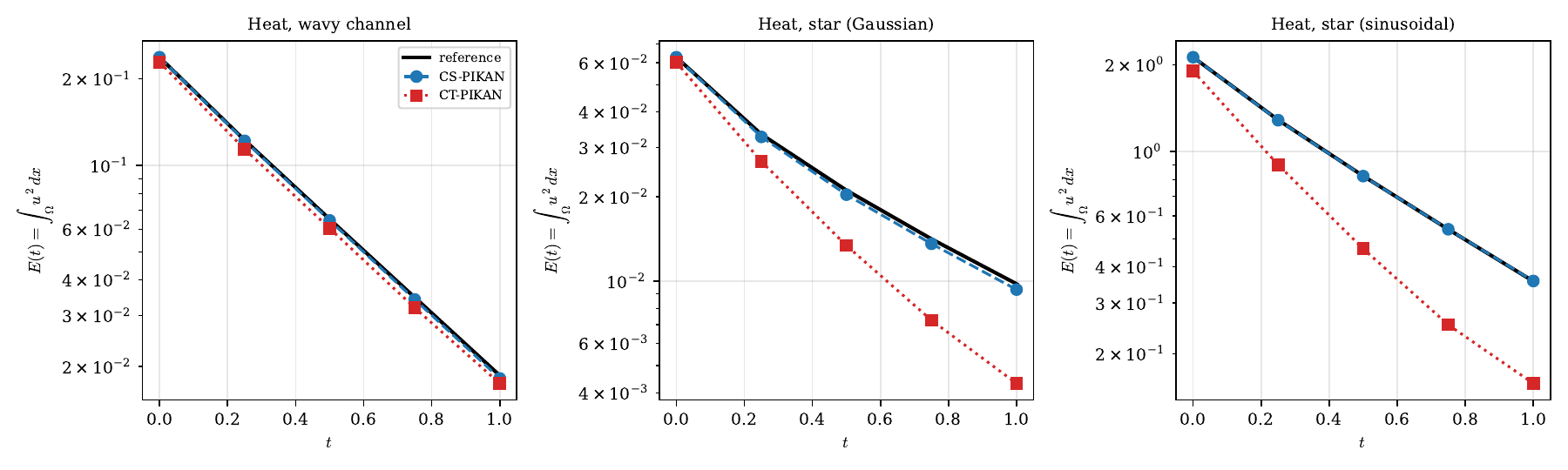}
\caption{Thermal energy $E(t)=\int_\Omega u^2\,dx$ of the reference, CS-PIKAN and CT-PIKAN solutions (seed 0).}
\label{fig:energy}
\end{figure}

We hold both the inner circle ($\xi=0$) and the outer star boundary ($\xi=1$) at $u=0$. Both cases solve the pure diffusion equation with $\alpha = 0.1$ over $T = 1$, with two boundary-compatible initial conditions:
\begin{equation}
	\begin{aligned}
		&u_0 = \sin(\pi\xi)\exp\!\left(-\frac{x^2+y^2}{2\sigma^2}\right),\quad \sigma=0.5 &&\text{(Gaussian)},\\
		&u_0 = \sin(\pi\xi)\sin(k_x x)\sin(k_y y),\quad k_x = k_y = \frac{2\pi}{4} &&\text{(sinusoidal)}.
	\end{aligned}
	\label{eq:58}
\end{equation}
The factor $\sin(\pi\xi)$ makes both initial conditions vanish on the two boundary components. Since the origin lies inside the hole, the Gaussian initial field is a ring profile that rises steeply from the inner circle. Fig.~\ref{fig:star_snap} shows the CS-PIKAN temperature evolution for both initial conditions, and Fig.~\ref{fig:heat} adds the pointwise error at $t=1$.
\paragraph{Thermal energy.} A physical check of a diffusion solution is the decay of the thermal energy $E(t)=\int_\Omega u^2\,dx$, which must decrease monotonically. Fig.~\ref{fig:energy} and Table~\ref{tab:energy} compare the energy of both solutions with the reference. On the channel both methods follow the reference. On the star, CS-PIKAN reproduces the reference energy closely, whereas CT-PIKAN dissipates energy far too fast: for both initial conditions its energy at $t=1$ is only $44\%$ of the reference value.

\begin{table}[H]\centering\caption{Thermal energy $E(t)=\int_\Omega u^2\,dx$: reference / CS-PIKAN / CT-PIKAN (seed 0).}\label{tab:energy}\footnotesize\begin{adjustbox}{max width=\linewidth}\begin{tabular}{lccc}\toprule
Case & $t=0$ & $t=0.5$ & $t=1$\\\midrule
Heat, wavy channel & 0.2380 / 0.2380 / 0.2288 & 0.0651 / 0.0644 / 0.0605 & 0.0186 / 0.0182 / 0.0174\\
Heat, star (Gaussian) & 0.0628 / 0.0628 / 0.0602 & 0.0211 / 0.0203 / 0.0134 & 0.0098 / 0.0093 / 0.0043\\
Heat, star (sinusoidal) & 2.1212 / 2.1212 / 1.9008 & 0.8240 / 0.8230 / 0.4617 & 0.3564 / 0.3571 / 0.1579\\
\bottomrule\end{tabular}\end{adjustbox}\end{table}

\subsection{Summary of accuracy}
Table~\ref{tab:main} shows that CS-PIKAN is more accurate than CT-PIKAN on all nine benchmarks, by a factor of 7--126 in relative $L_2$ error. 
\label{sec:summary}
\begin{table}[H]\centering\caption{Accuracy and cost on all benchmarks: relative $L_2$ and $L_\infty$ errors (mean $\pm$ std over 3 seeds) and wall-clock training time on one CPU core.}\label{tab:main}\footnotesize\setlength{\tabcolsep}{4pt}
\begin{adjustbox}{max width=\linewidth}\begin{tabular}{llcccr}\toprule
Case & Method & rel.~$L_2$ & rel.~$L_\infty$ & Time (s) & Params\\\midrule
\multirow{2}{*}{Advection, Cartesian} & CS-PIKAN & $2.88\times10^{-3}\pm1.84\times10^{-4}$ & $6.82\times10^{-3}\pm1.22\times10^{-3}$ & 28 & 18,816 \\
 & CT-PIKAN & $2.90\times10^{-1}\pm6.04\times10^{-2}$ & $5.02\times10^{-1}\pm4.84\times10^{-2}$ & 64 & 15,777 \\ \midrule
\multirow{2}{*}{Advection, polar ($\beta=2$)} & CS-PIKAN & $5.06\times10^{-3}\pm1.14\times10^{-3}$ & $1.26\times10^{-2}\pm3.67\times10^{-3}$ & 67 & 18,816 \\
 & CT-PIKAN & $3.84\times10^{-2}\pm4.26\times10^{-3}$ & $2.16\times10^{-1}\pm2.22\times10^{-2}$ & 68 & 15,777 \\ \midrule
\multirow{2}{*}{Advection, polar ($\beta=8$)} & CS-PIKAN & $1.62\times10^{-2}\pm3.31\times10^{-3}$ & $3.88\times10^{-2}\pm8.18\times10^{-3}$ & 63 & 18,816 \\
 & CT-PIKAN & $6.41\times10^{-1}\pm9.96\times10^{-2}$ & $9.87\times10^{-1}\pm4.41\times10^{-2}$ & 62 & 15,777 \\ \midrule
\multirow{2}{*}{Advection, wavy} & CS-PIKAN & $3.29\times10^{-3}\pm5.71\times10^{-4}$ & $9.82\times10^{-3}\pm1.04\times10^{-3}$ & 52 & 18,816 \\
 & CT-PIKAN & $4.16\times10^{-1}\pm3.50\times10^{-2}$ & $1.43\times10^{0}\pm4.88\times10^{-2}$ & 58 & 15,777 \\ \midrule
\multirow{2}{*}{Poisson, annulus} & CS-PIKAN & $6.30\times10^{-4}\pm4.67\times10^{-4}$ & $7.21\times10^{-4}\pm3.57\times10^{-4}$ & 32 & 12,688 \\
 & CT-PIKAN & $2.07\times10^{-2}\pm5.76\times10^{-3}$ & $4.22\times10^{-2}\pm1.02\times10^{-2}$ & 328 & 14,721 \\ \midrule
\multirow{2}{*}{Poisson, sampled boundary} & CS-PIKAN & $9.32\times10^{-4}\pm1.14\times10^{-4}$ & $1.20\times10^{-3}\pm2.04\times10^{-4}$ & 163 & 12,688 \\
 & CT-PIKAN & $3.91\times10^{-2}\pm3.69\times10^{-3}$ & $8.75\times10^{-2}\pm4.88\times10^{-3}$ & 365 & 14,721 \\ \midrule
\multirow{2}{*}{Heat, wavy channel} & CS-PIKAN & $6.17\times10^{-3}\pm1.64\times10^{-3}$ & $7.86\times10^{-3}\pm5.21\times10^{-4}$ & 34 & 18,816 \\
 & CT-PIKAN & $4.43\times10^{-2}\pm3.74\times10^{-3}$ & $7.87\times10^{-2}\pm8.61\times10^{-3}$ & 205 & 11,617 \\ \midrule
\multirow{2}{*}{Heat, star (Gaussian)} & CS-PIKAN & $2.37\times10^{-2}\pm3.43\times10^{-3}$ & $2.63\times10^{-2}\pm3.10\times10^{-3}$ & 90 & 18,960 \\
 & CT-PIKAN & $4.88\times10^{-1}\pm5.22\times10^{-3}$ & $4.78\times10^{-1}\pm1.38\times10^{-3}$ & 147 & 15,777 \\ \midrule
\multirow{2}{*}{Heat, star (sinusoidal)} & CS-PIKAN & $5.26\times10^{-3}\pm6.22\times10^{-4}$ & $8.26\times10^{-3}\pm1.78\times10^{-3}$ & 94 & 18,960 \\
 & CT-PIKAN & $2.82\times10^{-1}\pm1.84\times10^{-2}$ & $2.75\times10^{-1}\pm2.09\times10^{-2}$ & 143 & 15,777 \\
\bottomrule\end{tabular}\end{adjustbox}\end{table}

Its training time is comparable on the advection problems, where the point-wise baseline needs only first derivatives (ratio 0.99--2.31), and 1.5--10.3$\times$ lower on the elliptic and parabolic problems, which need second derivatives. The parameter counts are similar ($1.3$--$1.9\times10^4$ against $1.2$--$1.6\times10^4$).\\ 

CT-PIKAN does not reach a useful accuracy (relative $L_2$ error above $0.1$) on 5 of the nine benchmarks: Cartesian, polar ($\beta=8$) and wavy advection, and both star heat problems. On these transport and star-domain problems, its penalty terms and scattered points leave the boundary and initial data (and, on the star, the seam $\eta=0\equiv1$) poorly enforced.\\ 

The hardest case for CS-PIKAN is the star heat problem with the Gaussian initial condition ($2.37\times10^{-2}$). There the solution amplitude is small (at most $0.13$ for $t\ge0.25$), and the error concentrates in the outer lobes of the star, where $\sqrt g$ is up to twelve times its minimum. Increasing the grid to $48\times96\times32$ and training for $10^4$ epochs reduced this error only to $1.8\times10^{-2}$ (seed 0), so we report the uniform configuration.\\

\begin{figure}[H]
\centering
\includegraphics[width=\textwidth]{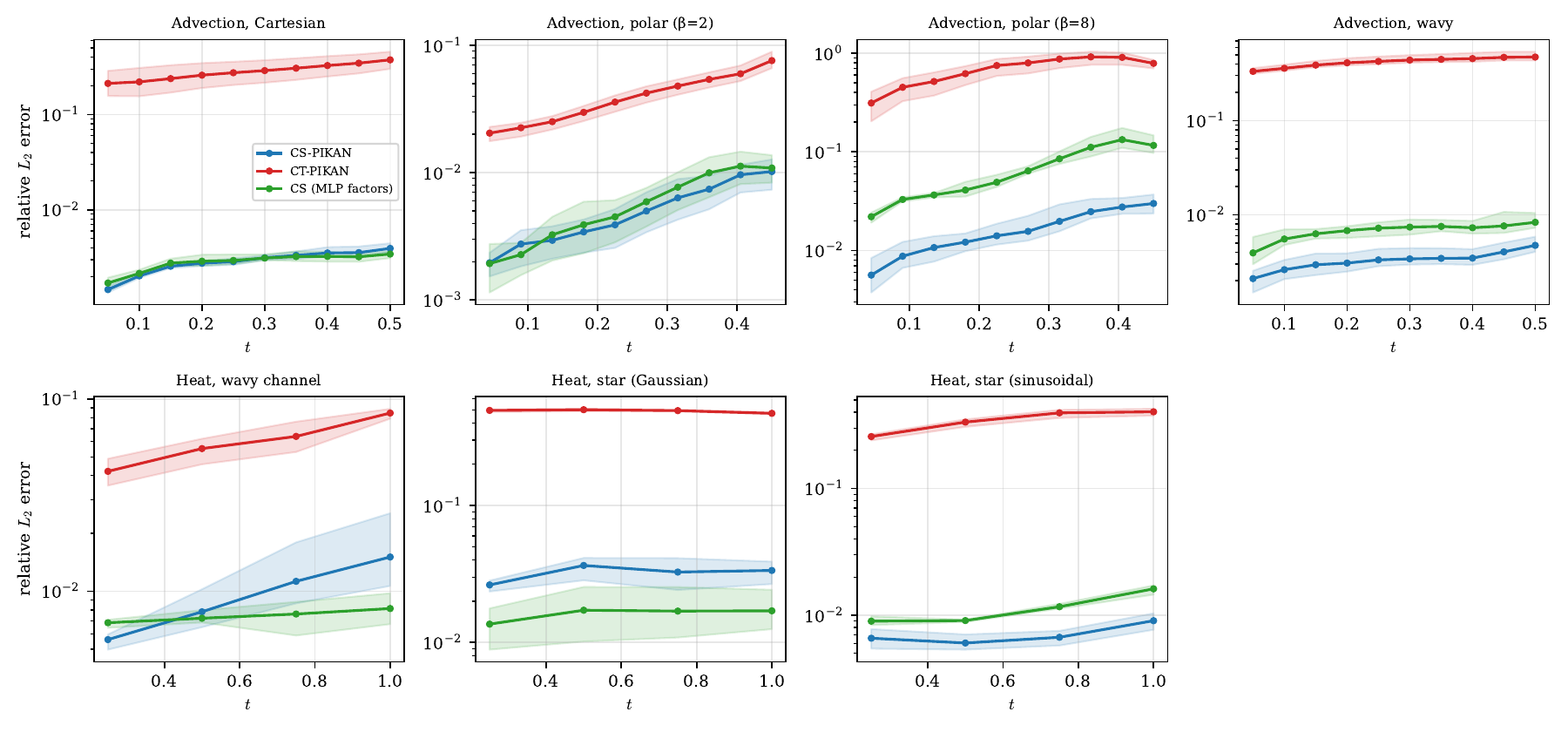}
\caption{Relative $L_2$ error versus time for all transient benchmarks: CS-PIKAN, CS-PIKAN with MLP factors, and CT-PIKAN. Lines show the geometric mean over three seeds.}
\label{fig:err_time}
\end{figure}

\begin{figure}[H]
	\centering
	\includegraphics[width=\textwidth]{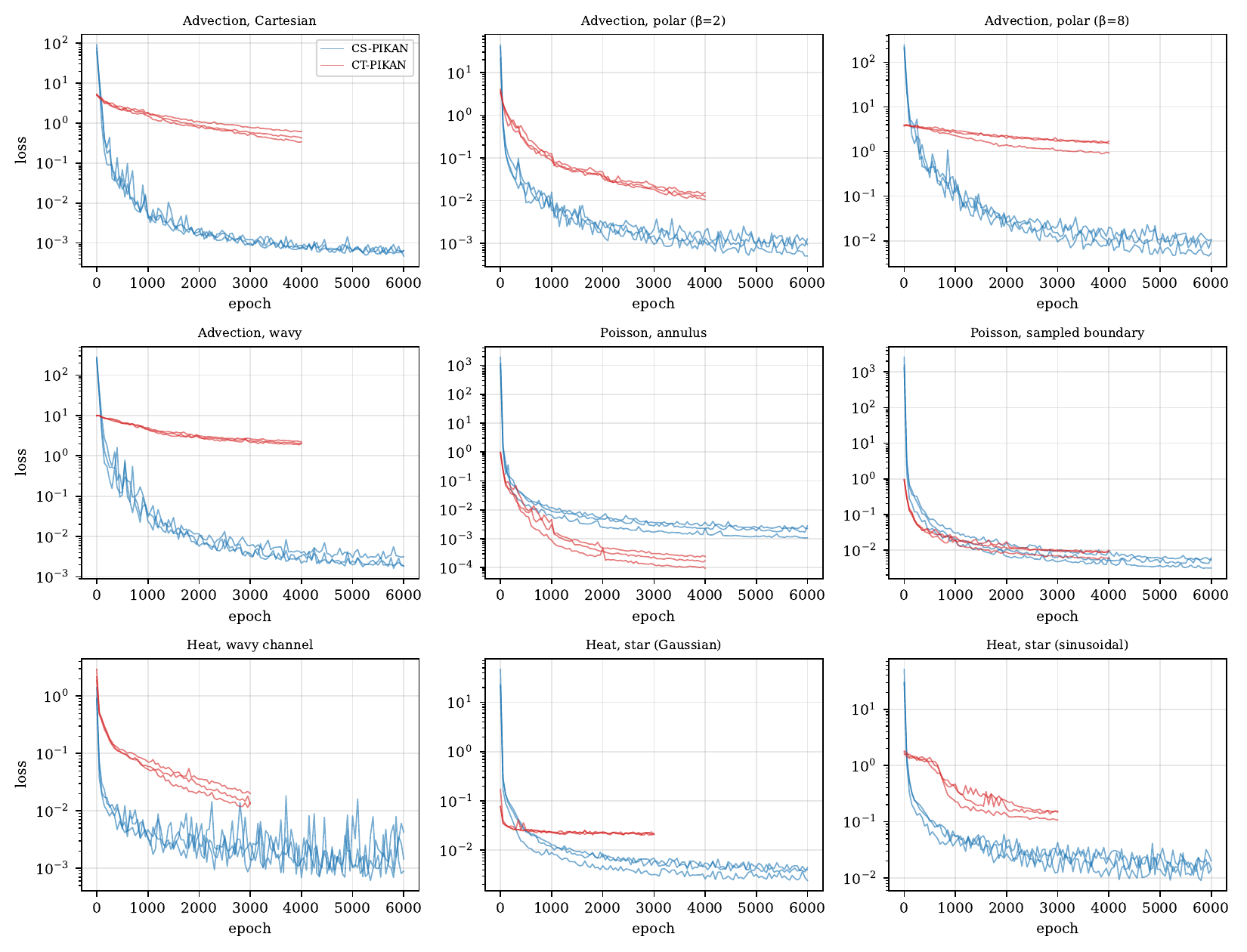}
	\caption{Training histories of all runs (three seeds per method).}
	\label{fig:histories}
\end{figure}

Fig.~\ref{fig:err_time} shows how the errors evolve in time for all transient problems, with the spread over seeds, and Table~\ref{tab:time} lists the errors at individual snapshot times. For the advection problems the CS-PIKAN error grows by a factor of $2$--$5$ over the time interval, as expected for transport of the residual, but it stays at least a factor of $6$ below that of CT-PIKAN at every snapshot (and two orders of magnitude below on the Cartesian and wavy domains). For the heat problems the CS-PIKAN error varies by less than a factor of three in time and stays at least a factor of $5$ below CT-PIKAN. Fig.~\ref{fig:histories} shows the training histories of all runs.

\noindent The CS-PIKAN loss is the single $\sqrt g$-weighted residual term; the CT-PIKAN loss includes the penalty terms, so the two are not directly comparable.
\begin{table}[H]\centering\caption{Relative $L_2$ error at individual snapshot times (mean over three seeds).}\label{tab:time}\footnotesize\begin{adjustbox}{max width=\linewidth}\begin{tabular}{llccc}\toprule
Case & Method & \multicolumn{3}{c}{error at snapshot times}\\\midrule
Advection, Cartesian & CS-PIKAN & $t=0.25$: $2.87\times10^{-3}$ & $t=0.50$: $3.96\times10^{-3}$\\
 & CT-PIKAN & $t=0.25$: $2.81\times10^{-1}$ & $t=0.50$: $3.80\times10^{-1}$\\
Advection, polar ($\beta=2$) & CS-PIKAN & $t=0.22$: $4.06\times10^{-3}$ & $t=0.45$: $1.05\times10^{-2}$\\
 & CT-PIKAN & $t=0.22$: $3.61\times10^{-2}$ & $t=0.45$: $7.64\times10^{-2}$\\
Advection, polar ($\beta=8$) & CS-PIKAN & $t=0.22$: $1.44\times10^{-2}$ & $t=0.45$: $3.04\times10^{-2}$\\
 & CT-PIKAN & $t=0.22$: $7.60\times10^{-1}$ & $t=0.45$: $7.91\times10^{-1}$\\
Advection, wavy & CS-PIKAN & $t=0.25$: $3.39\times10^{-3}$ & $t=0.50$: $4.80\times10^{-3}$\\
 & CT-PIKAN & $t=0.25$: $4.27\times10^{-1}$ & $t=0.50$: $4.75\times10^{-1}$\\
Heat, wavy channel & CS-PIKAN & $t=0.25$: $5.63\times10^{-3}$ & $t=0.50$: $7.98\times10^{-3}$ & $t=1.00$: $1.62\times10^{-2}$\\
 & CT-PIKAN & $t=0.25$: $4.24\times10^{-2}$ & $t=0.50$: $5.57\times10^{-2}$ & $t=1.00$: $8.46\times10^{-2}$\\
Heat, star (Gaussian) & CS-PIKAN & $t=0.25$: $2.64\times10^{-2}$ & $t=0.50$: $3.70\times10^{-2}$ & $t=1.00$: $3.40\times10^{-2}$\\
 & CT-PIKAN & $t=0.25$: $4.94\times10^{-1}$ & $t=0.50$: $5.00\times10^{-1}$ & $t=1.00$: $4.70\times10^{-1}$\\
Heat, star (sinusoidal) & CS-PIKAN & $t=0.25$: $6.65\times10^{-3}$ & $t=0.50$: $6.07\times10^{-3}$ & $t=1.00$: $9.12\times10^{-3}$\\
 & CT-PIKAN & $t=0.25$: $2.57\times10^{-1}$ & $t=0.50$: $3.35\times10^{-1}$ & $t=1.00$: $4.04\times10^{-1}$\\
\bottomrule\end{tabular}\end{adjustbox}\end{table}

\subsection{Ablation}
\label{sec:ablation}
\begin{figure}[H]
\centering
\includegraphics[width=\textwidth]{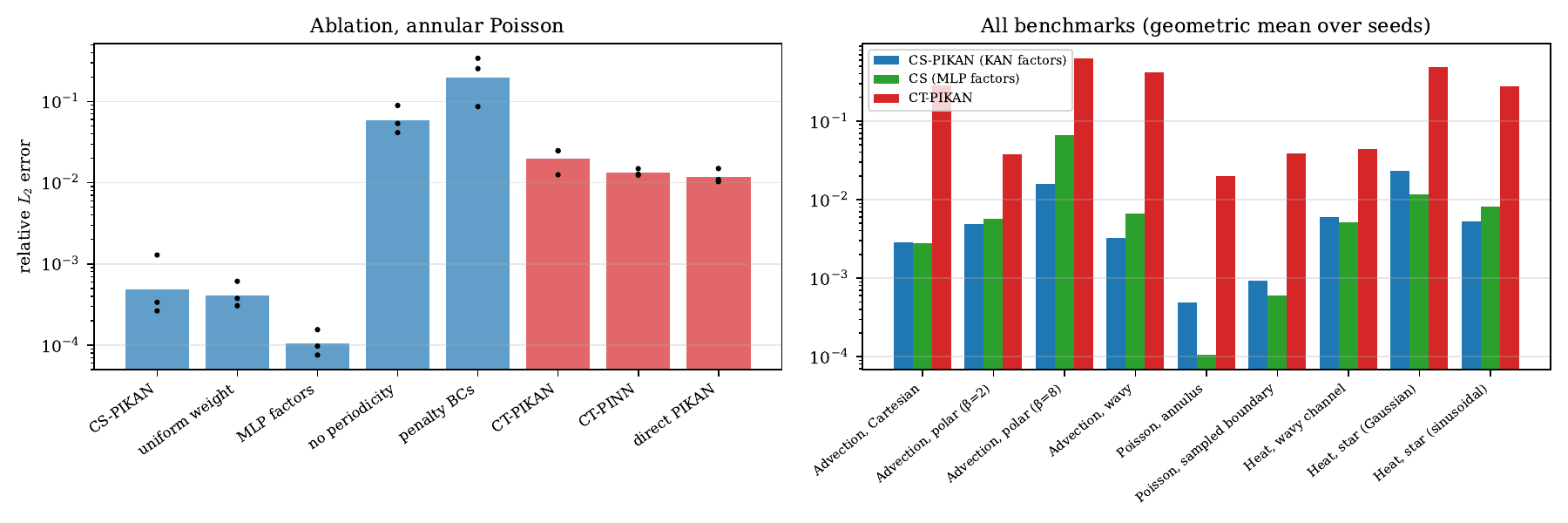}
\caption{Left: ablation on the annular Poisson problem (bars: geometric mean; dots: individual seeds). Right: relative $L_2$ error of CS-PIKAN with KAN factors, CS-PIKAN with MLP factors, and CT-PIKAN on all benchmarks.}
\label{fig:bars}
\end{figure}

Table~\ref{tab:ablation} isolates each component on the annular Poisson problem. All variants use the same seeds and budgets. The point-wise rows use the CT-PIKAN configuration; ``direct'' trains the point-wise RBF-KAN on physical coordinates $(x,y)$ with the Cartesian Laplacian. Fig.~\ref{fig:bars} (left) shows the same data with the individual seeds.
\begin{table}[H]\centering\caption{Ablation on the annular Poisson problem (3 seeds).}\label{tab:ablation}\footnotesize\setlength{\tabcolsep}{3pt}
\begin{adjustbox}{max width=\linewidth}\begin{tabular}{lcccccccc}\toprule
Variant & Repr. & Constraints & $\eta$-basis & Loss weight & Params & rel.~$L_2$ & Time (s)\\\midrule
CS-PIKAN (full) & sep. & exact & Fourier & $\sqrt g$ & 12,688 & $6.30\times10^{-4}\pm4.67\times10^{-4}$ & 32\\
\quad Chebyshev in $\eta$ (seam) & sep. & exact & Cheb. & $\sqrt g$ & 12,544 & $6.14\times10^{-2}\pm2.03\times10^{-2}$ & 49\\
\quad penalty BCs & sep. & penalty & Fourier & $\sqrt g$ & 12,688 & $2.27\times10^{-1}\pm1.05\times10^{-1}$ & 33\\
\quad unweighted loss & sep. & exact & Fourier & uniform & 12,688 & $4.32\times10^{-4}\pm1.31\times10^{-4}$ & 42\\
\quad MLP factors (SPINN-type) & sep. & exact & Fourier & $\sqrt g$ & 2,032 & $1.10\times10^{-4}\pm3.36\times10^{-5}$ & 27\\
CT-PIKAN (RBF-KAN) & point & penalty & -- & uniform & 14,721 & $2.07\times10^{-2}\pm5.76\times10^{-3}$ & 328\\
CT-PINN (MLP) & point & penalty & -- & uniform & 12,737 & $1.33\times10^{-2}\pm1.11\times10^{-3}$ & 270\\
Direct PIKAN, $(x,y)$ inputs & point & penalty & -- & uniform & 14,721 & $1.21\times10^{-2}\pm2.07\times10^{-3}$ & 289\\
\bottomrule\end{tabular}\end{adjustbox}\end{table}
Exact periodicity and exact constraints are the decisive ingredients. Replacing the Fourier $\eta$-layer by a Chebyshev one, i.e.\ dropping exact periodicity, raises the error from $6.30\times10^{-4}$ to $6.14\times10^{-2}$, and penalty boundary conditions (weight $10$, as in CT-PIKAN) raise it to $2.27\times10^{-1}$. The $\sqrt g$ weight is not needed for accuracy on this problem (uniform weighting gives $4.32\times10^{-4}$); its role is to make the loss an estimator of the physical residual norm and hence of the certificate (Proposition~\ref{prop:unbiased}). MLP factors give $1.10\times10^{-4}$, which is lower on this smooth problem; this trade-off is examined on all benchmarks below. Among the point-wise solvers, the coordinate transformation on its own does not help: the direct PIKAN on physical coordinates ($1.21\times10^{-2}$) is more accurate than CT-PIKAN ($2.07\times10^{-2}$), whose raw $\eta$ input leaves the seam unconstrained, and an MLP front end performs similarly (CT-PINN $1.33\times10^{-2}$). The transformation pays off only together with the separable, exactly constrained and exactly periodic representation.

\subsection{Choice of factor network}
\label{sec:factor}
The ablation shows that, on the smooth elliptic annulus problem, $\tanh$-MLP factors (a SPINN-type factor transplanted into our curvilinear framework) are more accurate than KAN factors. Table~\ref{tab:factor} therefore compares the two factor families on all nine benchmarks, keeping every other ingredient fixed.
\begin{table}[H]\centering\caption{Factor family inside the same separable curvilinear framework: Chebyshev-KAN factors (CS-PIKAN, $1.3$--$1.9\times10^4$ parameters) versus $\tanh$-MLP factors ($2.0$--$3.0\times10^3$ parameters); relative $L_2$ error, 3 seeds. Ratio $>1$ favours KAN.}\label{tab:factor}\footnotesize
\begin{adjustbox}{max width=\linewidth}\begin{tabular}{lcccc}\toprule
Case & KAN factors & MLP factors & ratio MLP/KAN & better\\\midrule
Advection, Cartesian & $2.88\times10^{-3}\pm1.84\times10^{-4}$ & $2.80\times10^{-3}\pm2.29\times10^{-4}$ & 0.97 & MLP\\
Advection, polar ($\beta=2$) & $5.06\times10^{-3}\pm1.14\times10^{-3}$ & $5.85\times10^{-3}\pm1.41\times10^{-3}$ & 1.16 & KAN\\
Advection, polar ($\beta=8$) & $1.62\times10^{-2}\pm3.31\times10^{-3}$ & $6.69\times10^{-2}\pm1.04\times10^{-2}$ & 4.12 & KAN\\
Advection, wavy & $3.29\times10^{-3}\pm5.71\times10^{-4}$ & $6.69\times10^{-3}\pm9.67\times10^{-4}$ & 2.03 & KAN\\
Poisson, annulus & $6.30\times10^{-4}\pm4.67\times10^{-4}$ & $1.10\times10^{-4}\pm3.36\times10^{-5}$ & 0.17 & MLP\\
Poisson, sampled boundary & $9.32\times10^{-4}\pm1.14\times10^{-4}$ & $6.18\times10^{-4}\pm1.50\times10^{-4}$ & 0.66 & MLP\\
Heat, wavy channel & $6.17\times10^{-3}\pm1.64\times10^{-3}$ & $5.13\times10^{-3}\pm3.86\times10^{-4}$ & 0.83 & MLP\\
Heat, star (Gaussian) & $2.37\times10^{-2}\pm3.43\times10^{-3}$ & $1.24\times10^{-2}\pm3.74\times10^{-3}$ & 0.52 & MLP\\
Heat, star (sinusoidal) & $5.26\times10^{-3}\pm6.22\times10^{-4}$ & $8.14\times10^{-3}\pm2.32\times10^{-4}$ & 1.55 & KAN\\
\bottomrule\end{tabular}\end{adjustbox}\end{table}
KAN factors are more accurate on 4 of the 9 benchmarks (polar advection with $\beta=2$ and $\beta=8$, wavy advection, and the sinusoidal star heat problem), and MLP factors on the remaining 5 (Cartesian advection, both Poisson problems, the wavy-channel heat problem and the Gaussian star heat problem). The KAN advantage is largest for transport of high-frequency patterns (polar sector, $\beta=8$: $1.62\times10^{-2}$ against $6.69\times10^{-2}$), which is consistent with the reduced spectral bias reported for KANs \cite{Mos1}. MLP factors are preferable for smooth elliptic solutions, where they reach a given accuracy with fewer parameters. Both factor families outperform the point-wise CT-PIKAN by roughly one to two orders of magnitude. The decisive ingredients are therefore the separable curvilinear representation, the exact constraints and the physical quadrature, while the factor family is a secondary, problem-dependent choice. We keep KAN factors as the default for two reasons. Their largest disadvantage (a factor $5.8$, annulus) occurs on smooth problems that both families already solve to well below $10^{-2}$, whereas their largest advantage occurs on the hardest transport problem, where it decides whether a $10^{-2}$ accuracy level is reached.

\subsection{Certified bounds}
\label{sec:res_bounds}
Table~\ref{tab:bounds} compares the certificates of Theorem~\ref{thm:apost} with the true errors of the trained CS-PIKAN models. For the transient problems both quantities are taken at the final time $T$.
\begin{table}[H]\centering\caption{Certified bounds (Theorem~\ref{thm:apost}) versus true errors for CS-PIKAN, both relative to $\|u\|_{L_2(\Omega)}$ (at $t=T$ for transient problems); effectivity = bound/true error ($\ge1$ certifies the bound).}\label{tab:bounds}\footnotesize
\begin{adjustbox}{max width=\linewidth}\begin{tabular}{llccc}\toprule
Case & Constant & true rel.~error & certified bound & effectivity\\\midrule
Advection, Cartesian & $\beta=0$ & $3.93\times10^{-3}\pm3.50\times10^{-4}$ & $2.40\times10^{-2}\pm9.21\times10^{-4}$ & 6.2 \\
Advection, polar ($\beta=2$) & $\beta=2$ & $1.05\times10^{-2}\pm2.27\times10^{-3}$ & $5.50\times10^{-2}\pm6.44\times10^{-3}$ & 5.4 \\
Advection, polar ($\beta=8$) & $\beta=2$ & $3.08\times10^{-2}\pm5.48\times10^{-3}$ & $1.58\times10^{-1}\pm2.17\times10^{-2}$ & 5.2 \\
Advection, wavy & $\beta=0$ & $4.82\times10^{-3}\pm7.76\times10^{-4}$ & $3.42\times10^{-2}\pm3.96\times10^{-3}$ & 7.1 \\
Poisson, annulus & $C_\Omega^2=0.912$ & $6.30\times10^{-4}\pm4.67\times10^{-4}$ & $5.80\times10^{-2}\pm8.40\times10^{-3}$ & 140.8 \\
Poisson, sampled boundary & $C_\Omega^2=1.664$ & $9.32\times10^{-4}\pm1.14\times10^{-4}$ & $2.27\times10^{-1}\pm2.45\times10^{-2}$ & 250.1 \\
Heat, wavy channel & $C_\Omega^2=0.122$ & $1.62\times10^{-2}\pm6.57\times10^{-3}$ & $2.26\times10^{-1}\pm5.05\times10^{-3}$ & 15.9 \\
Heat, star (Gaussian) & $C_\Omega^2=2.323$ & $3.40\times10^{-2}\pm5.27\times10^{-3}$ & $5.56\times10^{-1}\pm5.62\times10^{-2}$ & 16.5 \\
Heat, star (sinusoidal) & $C_\Omega^2=2.323$ & $9.12\times10^{-3}\pm1.10\times10^{-3}$ & $1.89\times10^{-1}\pm1.78\times10^{-2}$ & 20.8 \\
\bottomrule\end{tabular}\end{adjustbox}\end{table}
All 27 certificates hold (effectivity $\ge1$). They are tightest for advection (effectivity 4.9--7.8), where the bound involves no Poincar\'e constant, intermediate for heat conduction (9.2--21.4), and loosest for the Poisson problems (47--332). For Poisson the bound $C_\Omega^2\|\mathcal R\|$ uses the Poincar\'e constant twice, and Table~\ref{tab:geom} shows that $C_\Omega^2$ overestimates the sharp $1/\lambda_1$ by a factor of 4--8 on these domains. It also bounds the error by the $L_2$ norm of the residual, whereas the sharp dual norm $\|\mathcal R\|_{H^{-1}}$ is much smaller for oscillatory residuals. Because the loss is an unbiased estimate of the squared residual norm (Proposition~\ref{prop:unbiased}), the certificate decreases together with the training loss. 

\begin{figure}[H]
\centering
\includegraphics[width=\textwidth]{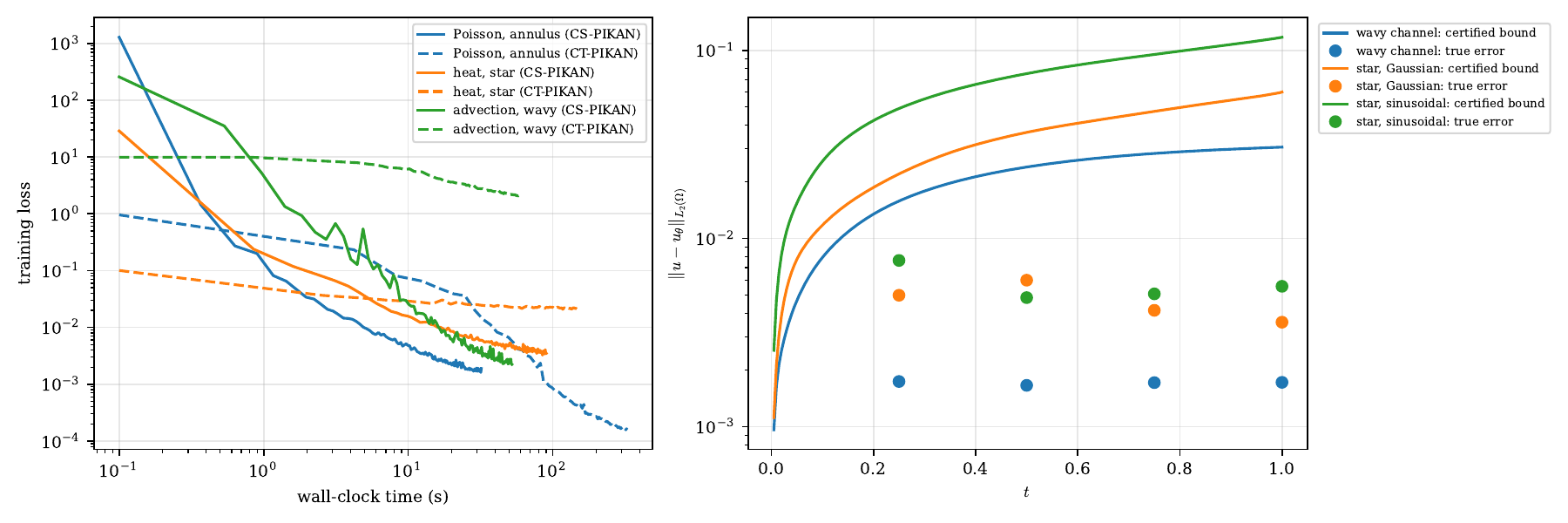}
\caption{Left: training histories (mean over seeds) of CS-PIKAN and CT-PIKAN on three benchmarks, plotted against wall-clock time. The two losses are not directly comparable: CT-PIKAN's includes penalty terms. Right: certified bound and true $L_2$ error over time for the heat problems (CS-PIKAN, seed 0).}
\label{fig:bounds}
\end{figure}

\begin{figure}[H]
\centering
\includegraphics[width=0.85\textwidth]{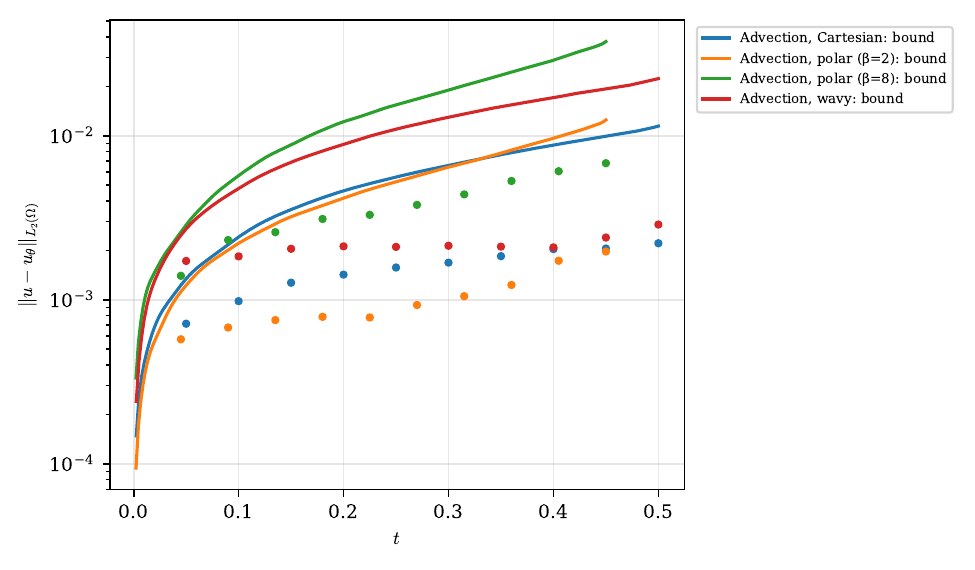}
\caption{Certified bound of Theorem~\ref{thm:apost}(c) (lines) and true $L_2$ error (markers) over time for the four advection problems (CS-PIKAN, seed 0).}
\label{fig:bound_adv}
\end{figure}

Figs.~\ref{fig:bounds} and \ref{fig:bound_adv} show the bounds and the true errors as functions of time.

\subsection{Computational cost and comparison with a classical body-fitted solver}
\label{sec:cost}
Table~\ref{tab:timing} breaks the training time down into time per epoch and number of epochs. Per epoch, CS-PIKAN is $3.0$--$15.4\times$ cheaper on the problems with second derivatives and $1.5$--$3.5\times$ cheaper on advection. Because it runs $6000$ epochs instead of $3000$--$4000$, the total training times on advection are comparable.

\begin{table}[H]\centering\caption{Wall-clock cost on one CPU core (compilation excluded): time per epoch, number of epochs and total training time.}\label{tab:timing}\footnotesize\begin{adjustbox}{max width=\linewidth}\begin{tabular}{lccccccc}\toprule
& \multicolumn{3}{c}{CS-PIKAN} & \multicolumn{3}{c}{CT-PIKAN} & \\\cmidrule(lr){2-4}\cmidrule(lr){5-7}
Case & ms/epoch & epochs & total (s) & ms/epoch & epochs & total (s) & ratio CT/CS\\\midrule
Advection, Cartesian & 4.6 & 6000 & 28 & 16.1 & 4000 & 64 & 2.31\\
Advection, polar ($\beta=2$) & 11.1 & 6000 & 67 & 17.1 & 4000 & 68 & 1.03\\
Advection, polar ($\beta=8$) & 10.5 & 6000 & 63 & 15.5 & 4000 & 62 & 0.99\\
Advection, wavy & 8.7 & 6000 & 52 & 14.5 & 4000 & 58 & 1.11\\
Poisson, annulus & 5.3 & 6000 & 32 & 82.1 & 4000 & 328 & 10.27\\
Poisson, sampled boundary & 27.1 & 6000 & 163 & 91.1 & 4000 & 365 & 2.24\\
Heat, wavy channel & 5.6 & 6000 & 34 & 68.4 & 3000 & 205 & 6.09\\
Heat, star (Gaussian) & 15.0 & 6000 & 90 & 48.9 & 3000 & 147 & 1.63\\
Heat, star (sinusoidal) & 15.6 & 6000 & 94 & 47.5 & 3000 & 143 & 1.52\\
\bottomrule\end{tabular}\end{adjustbox}\end{table}

\begin{figure}[H]
\centering
\includegraphics[width=\textwidth]{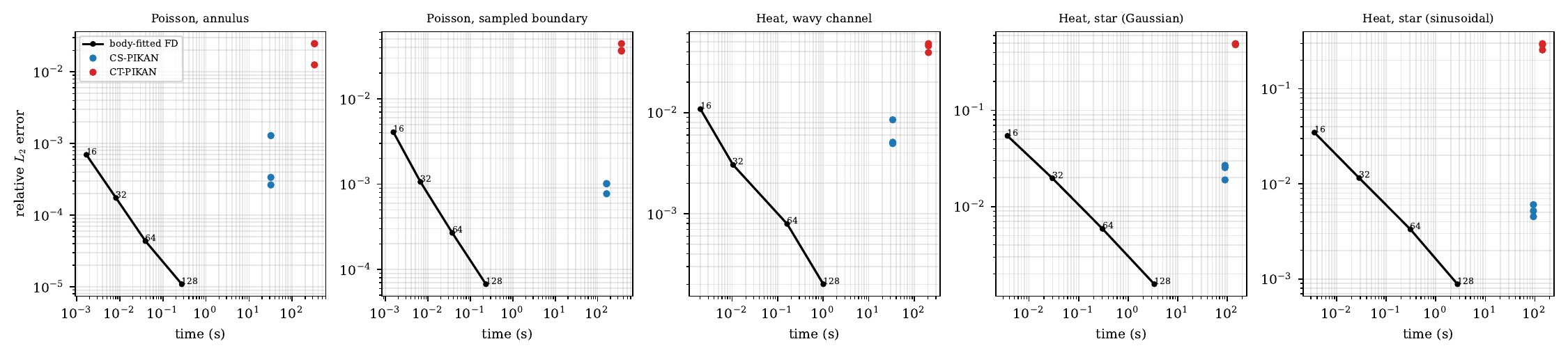}
\caption{Error versus cost. Black: body-fitted finite differences on successively refined grids (labels: cells in $\xi$). Blue: CS-PIKAN; red: CT-PIKAN (three seeds each, plotted at the training time).}
\label{fig:fdcost}
\end{figure}

For the problems with a finite-difference reference, we also measured the cost of the classical body-fitted solver (sparse factorisation and Crank--Nicolson time stepping; matrix assembly, done by an unoptimised Python loop, is excluded) on a sequence of grids (Table~\ref{tab:fdcost}, Fig.~\ref{fig:fdcost}). At the coarsest grid whose error is no larger than the CS-PIKAN error, the results are as follows. For the annulus Poisson problem, FD $32\times64$ reaches $1.74\times10^{-4}$ in 0.01~s, against $6.30\times10^{-4}$ in 32~s for CS-PIKAN. For the sampled-boundary Poisson problem, FD $64\times128$ reaches $2.68\times10^{-4}$ in 0.04~s, against $9.32\times10^{-4}$ in 163~s. For the wavy-channel heat problem, FD $32\times32$ reaches $3.04\times10^{-3}$ in 0.01~s, against $6.17\times10^{-3}$ in 34~s. For the star heat problem (Gaussian), FD $32\times64$ reaches $1.97\times10^{-2}$ in 0.03~s, against $2.37\times10^{-2}$ in 90~s. For these two-dimensional linear forward problems the classical solver is therefore three to four orders of magnitude cheaper at equal accuracy, and we do not claim otherwise. Compared with CT-PIKAN, CS-PIKAN narrows the gap by being at least as fast and $7$--$126$ times more accurate, while it keeps the mesh-free, continuous and differentiable representation and adds a certificate.

\begin{table}[H]\centering\caption{Body-fitted finite-difference solver: error against the Richardson reference and solve time (factorisation and time stepping) per grid.}\label{tab:fdcost}\footnotesize\begin{tabular}{lccc}\toprule
Case & grid & rel.~$L_2$ & solve time (s)\\\midrule
Poisson, annulus & $16\times32$ & $6.96\times10^{-4}$ & 0.002\\
 & $32\times64$ & $1.74\times10^{-4}$ & 0.008\\
 & $64\times128$ & $4.35\times10^{-5}$ & 0.039\\
 & $128\times256$ & $1.09\times10^{-5}$ & 0.275\\\midrule
Poisson, sampled boundary & $16\times32$ & $4.06\times10^{-3}$ & 0.002\\
 & $32\times64$ & $1.06\times10^{-3}$ & 0.007\\
 & $64\times128$ & $2.68\times10^{-4}$ & 0.038\\
 & $128\times256$ & $6.72\times10^{-5}$ & 0.233\\\midrule
Heat, wavy channel & $16\times16$ & $1.08\times10^{-2}$ & 0.002\\
 & $32\times32$ & $3.04\times10^{-3}$ & 0.011\\
 & $64\times64$ & $7.92\times10^{-4}$ & 0.162\\
 & $128\times128$ & $2.01\times10^{-4}$ & 1.011\\\midrule
Heat, star (Gaussian) & $16\times32$ & $5.44\times10^{-2}$ & 0.004\\
 & $32\times64$ & $1.97\times10^{-2}$ & 0.030\\
 & $64\times128$ & $5.88\times10^{-3}$ & 0.304\\
 & $128\times256$ & $1.56\times10^{-3}$ & 3.384\\\midrule
Heat, star (sinusoidal) & $16\times32$ & $3.46\times10^{-2}$ & 0.004\\
 & $32\times64$ & $1.15\times10^{-2}$ & 0.028\\
 & $64\times128$ & $3.33\times10^{-3}$ & 0.307\\
 & $128\times256$ & $8.82\times10^{-4}$ & 2.744\\
\bottomrule\end{tabular}\end{table}

\section{Discussion}
\label{sec:discussion}

\paragraph{What the coordinate transformation does and does not do.} The ablation separates two roles of the transformation. As a change of variables in the operator, which is all that CT-PIKAN uses, it is not beneficial on its own: on the annulus the point-wise network trained directly on $(x,y)$ is more accurate than the point-wise network trained on $(\xi,\eta)$. The reason is the seam: a closed physical curve is cut open, and a network with raw inputs has no reason to be continuous across the cut (Fig.~\ref{fig:seam}). As a change of the \emph{domain of the representation}, which is what CS-PIKAN exploits, the transformation is decisive: it is what makes the solution expressible as a sum of products, the data imposable by an envelope and a lift, the loss computable by tensor quadrature, and the periodicity exact. Every one of these properties requires a box.

\paragraph{Why exact constraints matter more than the architecture.} With penalties, the boundary and initial conditions are satisfied only to the extent that their loss terms are small relative to the residual. For transport problems an inflow error is advected into the whole domain, and for diffusion on the star a boundary error at the inner circle feeds the large outer lobes. This explains why CT-PIKAN fails precisely on the advection and star problems (Table~\ref{tab:main}), and why its thermal energy on the star decays far too fast (Fig.~\ref{fig:energy}). With exact constraints, the only remaining error source is the residual, which the certificates of Theorem~\ref{thm:apost} turn into a guaranteed bound.

\paragraph{KAN versus MLP factors.} Inside the separable framework the factor family matters much less than the framework itself (Table~\ref{tab:factor}). KAN factors are more accurate for transported high-frequency patterns, consistent with the reduced spectral bias reported for KANs \cite{Mos1}; MLP factors are more accurate, with fewer parameters, for smooth elliptic solutions. A practical recommendation is to use KAN factors for transport-dominated problems and to treat the factor family as a tunable choice for diffusion-dominated ones.

\paragraph{The role of the $\sqrt g$ weight.} The weight does not improve accuracy on the annulus (Table~\ref{tab:ablation}). Its justification is theoretical: it makes the training loss an unbiased estimate of the physical residual norm (Proposition~\ref{prop:unbiased}), so that the quantity being minimised is, up to the constants of Theorem~\ref{thm:apost}, a guaranteed bound on the error (Corollary~\ref{cor:energy}).

\paragraph{When to use a network solver.} For two-dimensional linear forward problems on a single smooth map, a body-fitted finite-difference solver is far cheaper at equal accuracy (Section~\ref{sec:cost}). The value of CS-PIKAN lies elsewhere: it delivers a continuous and differentiable solution, imposes the data exactly, provides a certificate without a reference solution, and keeps the same code for every map. The separable structure also makes higher-dimensional extensions attractive, since the network cost grows with the sum of the grid sizes. Inverse, parametric and design problems, in which the solution must be differentiated with respect to parameters or geometry, are natural applications.

\section{Limitations}
\label{sec:limits}
\begin{itemize}
\item \textbf{Single-block maps.} CS-PIKAN needs one smooth diffeomorphism from $\Q$ (or from $(0,1)\times\mathbb T$) onto $\Omega$. Domains with several holes or corners require multi-block maps with interface conditions, which would break global separability, or a learned map \cite{JacobiNet,Khairullin2025}. Neither has been tested here.
\item \textbf{Separability.} The rank needed grows with the anisotropy of the solution relative to the computational coordinates. The two hardest benchmarks for CS-PIKAN are the transported high-frequency pattern of the polar case with $\beta=8$ and the small-amplitude star heat problem with the Gaussian initial condition, both with errors around $2\times10^{-2}$ (Table~\ref{tab:main}).
\item \textbf{Certificates.} The bounds are guaranteed up to quadrature error, but they can be pessimistic (effectivity indices are reported in Table~\ref{tab:bounds}). They cover linear problems with exact data; nonlinear problems would need stability estimates specific to each equation.
\item \textbf{Problem class.} All benchmarks are linear, scalar and two-dimensional. Nonlinear, multi-physics and three-dimensional problems are not covered.
\item \textbf{Cost.} For these two-dimensional linear forward problems, a body-fitted finite-difference solver remains faster at equal accuracy. Network solvers are attractive when a continuous, differentiable representation is needed, or for inverse and parametric extensions.
\end{itemize}

\section{Conclusion}
\label{sec:conc}
This study developed physics-informed Kolmogorov--Arnold networks for partial differential equations on curvilinear domains in two steps. The first step, CT-PIKAN, integrates coordinate-transformation theory with PIKANs: the PDE is solved on a regular computational domain, and all geometric quantities (the Jacobian, its inverse, the metric tensor, its determinant and the derivatives that enter the Laplace--Beltrami operator) are evaluated by automatic differentiation of the map. This removes the manual derivation of metric coefficients and makes the formulation independent of the geometry; closed-form, numerically constructed and sampled-boundary maps are handled by the same code.\\

The second step showed that a coordinate transformation is more than a way to rewrite the operator: it changes what the solver can be. Pulling a curved domain back to a box makes it possible to use a separable, topology-adapted KAN representation, to impose all boundary and initial data exactly through an envelope and a transfinite lift, and to compute the residual with a $\sqrt g$-weighted tensor quadrature. The same ingredients provide what is needed for certified error bounds, whose constants follow from the metric of the map. On nine benchmarks covering hyperbolic, elliptic and parabolic problems on seven maps, CS-PIKAN was more accurate than the point-wise CT-PIKAN by 7--126 at equal or lower training cost, and its training loss is an unbiased estimate of a guaranteed error bound. The ablation identified exact periodicity and exact constraints as the decisive ingredients; KAN factors are preferable for transport-dominated problems and MLP factors for smooth elliptic ones.\\

Future work will address multi-block and learned maps for multiply connected and strongly non-convex domains, three-dimensional curvilinear geometries, moving and deformable domains, nonlinear and coupled multiphysics systems with problem-specific stability estimates, and inverse and parametric problems, where a continuous, differentiable and certified surrogate is most valuable.

\paragraph{Code and data availability.} All code (JAX), including the scripts that produce every table and figure, is available at \url{https://github.com/m-heravifard/CT-PIKAN}.

\appendix
\section*{Appendix A: Closed-form metrics of the maps}
\label{app:metrics}
The following closed forms were used only to verify the autograd metrics; the solvers never use them. With $\sqrt g=x_\xi y_\eta-x_\eta y_\xi$, $g_{11}=x_\xi^2+y_\xi^2$, $g_{12}=x_\xi x_\eta+y_\xi y_\eta$ and $g_{22}=x_\eta^2+y_\eta^2$:
\begin{itemize}
\item \emph{Polar sector} ($r=\tfrac12+\tfrac12\xi$, $\theta=\tfrac\pi2\eta$): $g_{11}=\tfrac14$, $g_{12}=0$, $g_{22}=\tfrac{\pi^2}{4}r^2$, $\sqrt g=\tfrac\pi4 r\in[0.393,0.785]$.
\item \emph{Annulus} ($r=\tfrac12+\tfrac32\xi$, $\theta=2\pi\eta$): $g_{11}=\tfrac94$, $g_{12}=0$, $g_{22}=4\pi^2r^2$, $\sqrt g=3\pi r\in[4.712,18.850]$, $\kappa=\tfrac{4\pi r}{3}\le8.38$.
\item \emph{Wavy (advection)}: Eq.~\eqref{eq:wavy_metrics}.
\item \emph{Wavy channel} ($x=\xi$, $y=(1-\eta)y_b(\xi)+\eta y_t(\xi)$): $y_\xi=(1-\eta)y_b'+\eta y_t'$, $y_\eta=y_t-y_b=1+2A\sin\tfrac\phi2\cos(2\pi\xi+\tfrac\phi2)$, $\sqrt g=y_\eta\in[0.901,1.099]$, $g_{11}=1+y_\xi^2$, $g_{12}=y_\xi y_\eta$, $g_{22}=y_\eta^2$.
\item \emph{Star} ($\theta=2\pi\eta$, $r=r_i+\xi(r_o(\theta)-r_i)$): with $r_\xi=r_o-r_i$ and $r_\eta=2\pi\xi\, r_o'(\theta)$, $r_o'=-r_{o,\rm base}\gamma n_p\sin n_p\theta$,
\[\sqrt g=2\pi\,r\,r_\xi\in[2.827,34.306],\qquad g_{11}=r_\xi^2,\qquad g_{12}=r_\xi r_\eta,\qquad g_{22}=r_\eta^2+4\pi^2r^2.\]
\item \emph{Sampled boundary (TFI)}: $\mathbf x_\xi=\mathbf c_{\rm out}-\mathbf c_{\rm in}$, $\mathbf x_\eta=(1-\xi)\mathbf c_{\rm in}'+\xi\mathbf c_{\rm out}'$, so that $\sqrt g$ is affine in $\xi$ for each $\eta$.
\end{itemize}
On $2000$ random points per map, the autograd values agree with these expressions to within $2.3\times10^{-13}$ in double precision.

\section*{Appendix B: Finite-difference reference scheme}
\label{app:fd}
The reference solutions solve the transformed equations in conservative form, $\Delta_x u=\frac{1}{\sqrt g}\big[\partial_\xi(M_{11}u_\xi+M_{12}u_\eta)+\partial_\eta(M_{12}u_\xi+M_{22}u_\eta)\big]$ with $M=\sqrt g\,G^{-1}$, on the uniform $(\xi,\eta)$ grid $\xi_i=ih_\xi$, $\eta_j=jh_\eta$ ($\eta$ periodic for the annular, star and TFI maps). The diagonal fluxes use the metric at the half points,
\[\partial_\xi(M_{11}u_\xi)\big|_{ij}\approx\frac{M_{11}^{i+\frac12,j}(u_{i+1,j}-u_{ij})-M_{11}^{i-\frac12,j}(u_{ij}-u_{i-1,j})}{h_\xi^2},\]
and analogously in $\eta$. The cross terms use central differences of the products at the nodes,
\[\partial_\xi(M_{12}u_\eta)\big|_{ij}\approx\frac{M_{12}^{i+1,j}\,\delta_\eta u_{i+1,j}-M_{12}^{i-1,j}\,\delta_\eta u_{i-1,j}}{2h_\xi},\qquad \delta_\eta u_{ij}=\frac{u_{i,j+1}-u_{i,j-1}}{2h_\eta},\]
which gives a nine-point stencil. All metric values are computed by automatic differentiation of the map at the required (half) points. The Poisson problems are solved by a sparse direct solver; the heat problems use Crank-Nicolson time stepping with an LU factorisation computed once. References use $256\times512$ (periodic) or $256\times256$ grids with $2000$ time steps, combined with the $128$-grid solution (and $1000$ time steps) by Richardson extrapolation. Table~\ref{tab:fdver} verifies second-order convergence, and the finite-difference eigenvalue problem for $-\Delta_x$ on the same grids provides the sharp values $1/\lambda_1$ in Table~\ref{tab:geom}.

\section*{Appendix C: Results of all individual runs}
\label{app:seeds}
\begin{table}[H]\centering\caption{Relative $L_2$ error of every run (seeds 0, 1, 2).}\label{tab:seeds}\footnotesize\begin{adjustbox}{max width=\linewidth}\begin{tabular}{llccc}\toprule
Case & Method & seed 0 & seed 1 & seed 2\\\midrule
Advection, Cartesian & CS-PIKAN & $3.13\times10^{-3}$ & $2.81\times10^{-3}$ & $2.69\times10^{-3}$ \\
 & CS (MLP) & $3.10\times10^{-3}$ & $2.55\times10^{-3}$ & $2.74\times10^{-3}$ \\
 & CT-PIKAN & $2.90\times10^{-1}$ & $3.63\times10^{-1}$ & $2.15\times10^{-1}$ \\\midrule
Advection, polar ($\beta=2$) & CS-PIKAN & $3.45\times10^{-3}$ & $5.75\times10^{-3}$ & $5.98\times10^{-3}$ \\
 & CS (MLP) & $6.51\times10^{-3}$ & $7.15\times10^{-3}$ & $3.89\times10^{-3}$ \\
 & CT-PIKAN & $3.89\times10^{-2}$ & $4.33\times10^{-2}$ & $3.29\times10^{-2}$ \\\midrule
Advection, polar ($\beta=8$) & CS-PIKAN & $1.30\times10^{-2}$ & $2.08\times10^{-2}$ & $1.50\times10^{-2}$ \\
 & CS (MLP) & $5.88\times10^{-2}$ & $6.02\times10^{-2}$ & $8.15\times10^{-2}$ \\
 & CT-PIKAN & $6.88\times10^{-1}$ & $7.32\times10^{-1}$ & $5.02\times10^{-1}$ \\\midrule
Advection, wavy & CS-PIKAN & $2.91\times10^{-3}$ & $4.10\times10^{-3}$ & $2.87\times10^{-3}$ \\
 & CS (MLP) & $5.71\times10^{-3}$ & $6.37\times10^{-3}$ & $8.01\times10^{-3}$ \\
 & CT-PIKAN & $4.65\times10^{-1}$ & $3.92\times10^{-1}$ & $3.90\times10^{-1}$ \\\midrule
Poisson, annulus & CS-PIKAN & $2.65\times10^{-4}$ & $1.29\times10^{-3}$ & $3.37\times10^{-4}$ \\
 & CS (MLP) & $9.76\times10^{-5}$ & $7.57\times10^{-5}$ & $1.55\times10^{-4}$ \\
 & CT-PIKAN & $2.47\times10^{-2}$ & $2.48\times10^{-2}$ & $1.25\times10^{-2}$ \\\midrule
Poisson, sampled boundary & CS-PIKAN & $1.02\times10^{-3}$ & $1.01\times10^{-3}$ & $7.72\times10^{-4}$ \\
 & CS (MLP) & $6.85\times10^{-4}$ & $7.59\times10^{-4}$ & $4.11\times10^{-4}$ \\
 & CT-PIKAN & $3.67\times10^{-2}$ & $4.43\times10^{-2}$ & $3.62\times10^{-2}$ \\\midrule
Heat, wavy channel & CS-PIKAN & $4.94\times10^{-3}$ & $8.48\times10^{-3}$ & $5.09\times10^{-3}$ \\
 & CS (MLP) & $4.62\times10^{-3}$ & $5.56\times10^{-3}$ & $5.20\times10^{-3}$ \\
 & CT-PIKAN & $4.80\times10^{-2}$ & $4.58\times10^{-2}$ & $3.92\times10^{-2}$ \\\midrule
Heat, star (Gaussian) & CS-PIKAN & $2.54\times10^{-2}$ & $2.68\times10^{-2}$ & $1.89\times10^{-2}$ \\
 & CS (MLP) & $1.30\times10^{-2}$ & $1.66\times10^{-2}$ & $7.52\times10^{-3}$ \\
 & CT-PIKAN & $4.81\times10^{-1}$ & $4.89\times10^{-1}$ & $4.94\times10^{-1}$ \\\midrule
Heat, star (sinusoidal) & CS-PIKAN & $5.21\times10^{-3}$ & $6.05\times10^{-3}$ & $4.52\times10^{-3}$ \\
 & CS (MLP) & $8.21\times10^{-3}$ & $8.38\times10^{-3}$ & $7.83\times10^{-3}$ \\
 & CT-PIKAN & $2.56\times10^{-1}$ & $2.98\times10^{-1}$ & $2.91\times10^{-1}$ \\
\bottomrule\end{tabular}\end{adjustbox}\end{table}

\end{document}